\documentclass[a4paper,reqno]{amsart}

\usepackage{iftex}
\ifpTeX
  \PassOptionsToPackage{dvipdfmx}{graphicx}
  \PassOptionsToPackage{dvipdfmx}{xcolor}
  \PassOptionsToPackage{dvipdfmx}{hyperref}
  \AtBeginDvi{}
\fi

\usepackage{amsmath}
\usepackage{amssymb}
\usepackage{amscd}
\usepackage{amsthm}
\usepackage{type1cm}
\usepackage{tcolorbox}
\usepackage{mathrsfs}
\usepackage{needspace}

\usepackage[colorlinks,citecolor=darkgreen,linkcolor=red]{hyperref}
\definecolor{darkgreen}{rgb}{0.0, 0.6, 0.0}

\usepackage[all]{xy}
\input xy
\xyoption{all}
\tcbuselibrary{breakable, skins, theorems}

\def\A{\mathcal{A}}
\def\B{\mathcal{B}}
\def\C{\mathcal{C}}
\def\E{\mathcal{E}}
\def\F{\mathcal{F}}
\def\H{\mathcal{H}}
\def\I{\mathcal{I}}

\def\K{\mathcal{K}}
\def\L{\mathcal{L}}

\def\O{\mathcal{O}}
\def\P{\mathcal{P}}
\def\S{\mathcal{S}}
\def\T{\mathcal{T}}
\def\U{\mathcal{U}}
\def\V{\mathcal{V}}
\def\W{\mathcal{W}}

\def\m{\mathfrak{m}}
\def\n{\mathfrak{n}}
\def\p{\mathfrak{p}}

\def\La{\Lambda}

\DeclareMathOperator{\md}{\mathsf{mod}}
\renewcommand{\mod}{\md}
\DeclareMathOperator{\proj}{\mathsf{proj}}
\DeclareMathOperator{\fl}{\mathsf{fl}}
\DeclareMathOperator{\CM}{\mathsf{CM}}
\DeclareMathOperator{\refl}{\mathsf{ref}}
\DeclareMathOperator{\add}{\mathsf{add}}

\DeclareMathOperator{\Fac}{\mathsf{Fac}}
\DeclareMathOperator{\Coh}{\mathsf{Coh}}
\DeclareMathOperator{\vect}{\mathsf{vect}}
\DeclareMathOperator{\per}{\mathsf{per}}
\DeclareMathOperator{\thick}{\mathsf{thick}}

\DeclareMathOperator{\Hom}{Hom}
\DeclareMathOperator{\End}{End}
\DeclareMathOperator{\Ext}{Ext}
\DeclareMathOperator{\op}{op}
\DeclareMathOperator{\Frac}{Frac}
\DeclareMathOperator{\Ker}{Ker}
\DeclareMathOperator{\im}{Im}

\DeclareMathOperator{\Cl}{Cl}
\DeclareMathOperator{\Spec}{Spec}
\DeclareMathOperator{\Sing}{Sing}
\DeclareMathOperator{\Ass}{Ass}
\DeclareMathOperator{\height}{ht}
\DeclareMathOperator{\Pic}{Pic}

\DeclareMathOperator{\RHom}{\mathbb{R}Hom}

\DeclareMathOperator{\tilt}{tilt}
\DeclareMathOperator{\ptilt}{ptilt}
\DeclareMathOperator{\modif}{modif}
\DeclareMathOperator{\MM}{MM}
\DeclareMathOperator{\rk}{rk}

\DeclareMathOperator{\Conv}{conv}
\DeclareMathOperator{\depth}{depth}

\def\gl{\mathop{\rm gl.dim}\nolimits}

\def\pd{\mathop{\rm proj.dim}\nolimits}

\theoremstyle{definition}
\newtheorem{Thm}{Theorem}[section]
\newtheorem{IntroThm}{Theorem}

\newtheorem{Lem}[Thm]{Lemma}
\newtheorem{Prop}[Thm]{Proposition}
\newtheorem{Cor}[Thm]{Corollary}
\newtheorem{Def}[Thm]{Definition}
\newtheorem{Ex}[Thm]{Example}
\newtheorem{Rem}[Thm]{Remark}
\newtheorem*{Rem*}{Remark}

\newcommand{\FRAC}[2]{\leavevmode\kern.1em\raise.5ex\hbox{\the\scriptfont0 #1}\kern-.1em/\kern-.15em\lower.25ex\hbox{\the\scriptfont0 #2}}

\title{Tilting bundles and modifying modules in relative dimension one}
\author{Ryu Tomonaga}
\address{Graduate School of Mathematical Sciences, The University of Tokyo, 3-8-1 Komaba, Meguro-ku, Tokyo, 153-8914, Japan}
\email{ryu-tomonaga@g.ecc.u-tokyo.ac.jp}

\begin{document}
\begin{abstract}
For projective morphisms $\pi\colon Y\to X=\Spec R$ over noetherian complete local rings, with fibres of dimension at most one and $\mathbb R\pi_*\O_Y\cong\O_X$, we introduce a tilting fan formed by cones of pretilting bundles. Maximal cones parametrize tilting bundles and facet adjacency encodes irreducible mutation, translating both into polyhedral geometry. Completeness of this fan in positive rank is equivalent to finiteness of tilting bundles modulo line-bundle twists. For closed fibres that are split ADE highest-root cycles or their Stein contractions, we identify tilting bundles with chambers of restricted affine Tits cones and determine the fan explicitly.

Over a normal Gorenstein complete local domain, reflexive pushforward takes pretilting bundles on normal projective birational models $\pi\colon Y\to X$ with fibres of dimension at most one and $\pi^!\O_X\cong\O_Y$ to modifying modules. The pushforwards of canonical tilting bundles recover each model as a module blowup and detect contractions by additive containment, giving an Auslander--McKay embedding.

The hypersurfaces $R=S[[u,v]]/(uv-f)$, with $S$ regular complete local and $f\in S$ a nonzero nonunit, generalize type $A$ Du Val surfaces and $cA_n$ threefolds to arbitrary dimension. We classify modifying modules by flags of divisors of $f$ and rank-one reflexive twists; flag blowups make the Auslander--McKay embedding an equivalence. If $f$ has $r\ge2$ prime factors counted with multiplicity, reflexive pushforward on any fixed maximal-flag model $Y$ gives a mutation-preserving surjection
\[
\operatorname{Cham}(\widetilde A_{r-1})\cong\tilt Y
\twoheadrightarrow\MM R.
\]
Its fibres are affine Weyl group orbits determined by repeated factors of $f$. All maximal-flag models and maximal modification algebras are derived equivalent. These algebras are noncommutative crepant resolutions precisely when $S/(p)$ is regular for every prime factor $p$ of $f$.
\end{abstract}

\maketitle
\enlargethispage{3\baselineskip}
\tableofcontents

\section*{Introduction}

Tilting bundles translate geometry into the representation theory of endomorphism algebras. In relative dimension one, we study three related problems: classifying these bundles from the closed fibre, recovering models and contractions from modules, and describing the ambiguity of realizing a module on a fixed model.

The derived McKay correspondence for quotient threefolds \cite{BKR} and Bridgeland's equivalences for threefold flops \cite{Bri02} made this relationship concrete. Van den Bergh gave an algebraic approach to flops using tilting bundles \cite{VdB04b} and introduced noncommutative crepant resolutions (NCCRs) as algebraic counterparts of crepant resolutions \cite{VdB04a}. Smooth crepant resolutions need not exist, and crepant models may retain singularities. To study such models, Iyama--Wemyss developed modifying and maximal modifying modules, whose endomorphism algebras are Cohen--Macaulay but may have infinite global dimension \cite{IW}. They related maximal modification algebras to three-dimensional $\mathbb Q$-factorial terminalizations \cite{IW14}. Related work includes constructions of tilting bundles \cite{Bei,Kap,BH,TU10,HP14,Har,Tom25a} and NCCRs \cite{Aus86,Bro,HN,SVdB17,SVdB20a,SVdB20b,Tom24}, as well as derived equivalences and mutation for NCCRs and modifying modules \cite{IR,IW13,IW,IW14,IW18,Wem18,HH24}.

Classifying all such generators calls for a picture of how they fit together. For finite-dimensional algebras, Asai described two-term silting cones using numerical torsion pairs and wall--chamber structures \cite{Asa21}. Broomhead--Pauksztello--Ploog--Woolf's heart fans place such cones in the wider setting of hearts of bounded t-structures \cite{BPPW24}. This picture also appears in birational geometry: Shimpi used heart fans to classify t-structures intermediate with respect to the perverse coherent heart on the local derived category of a complex cDV flopping contraction \cite{Shi26}.

Line-bundle twists can produce infinitely many tilting bundles, while mutation changes their summands. Our starting point is that the closed fibre controls both the bundles and the fan of their numerical cones, without imposing birational hypotheses.

In Part \ref{part:tilting bundles}, let $(R,\m,k)$ be a noetherian complete local ring and let $\pi\colon Y\to X=\Spec R$ be projective, with
\[
\mathbb R\pi_*\O_Y\cong\O_X,\qquad
C:=Y\times_Rk,\qquad \dim C\le1.
\]
The residue field is arbitrary, and mixed characteristic is allowed. A vector bundle is \emph{pretilting} if it has no positive self-extensions, and \emph{tilting} if it also generates $\per Y$. Write $\ptilt Y$ and $\tilt Y$ for the isomorphism classes of basic pretilting and tilting bundles, respectively; \emph{basic} means that no indecomposable summand is repeated. For $U=\bigoplus_iU_i\in\ptilt Y$ with $U_i$ indecomposable, put
\[
V:=K_0(\per Y)_{\mathbb R},\qquad
C(U):=\sum_i\mathbb R_{\ge0}[U_i],\qquad
\Sigma(Y):=\{C(U)\mid U\in\ptilt Y\}.
\]
Let $\Coh_0Y$ denote sheaves supported on the closed fibre. For a vector bundle $\V$ and $M\in\Coh_0Y$, use the Euler pairing
\[
\langle[\V],[M]\rangle:=\sum_{i\ge0}(-1)^i\ell_R\Ext_Y^i(\V,M),
\qquad \theta(M):=\langle\theta,[M]\rangle\quad(\theta\in V),
\]
extended linearly to $V$. A nonzero $M$ is $\theta$-semistable if $\theta(M)=0$ and $\theta(N)\le0$ for every subsheaf $N\subseteq M$. In $V_{>0}:=\{\theta\in V\mid\rk\theta>0\}$, let $\mathrm{Wall}$ be the locus admitting such a sheaf. Chambers are the connected components of $V_{>0}\setminus\overline{\mathrm{Wall}}$, with closure taken in $V_{>0}$. Sections \ref{sec:tilting-torsion}--\ref{sec:completeness} establish the following.

\begin{IntroThm}[The tilting fan and closed-fibre reduction]\label{intro:tilting-fan}
The cones $\Sigma(Y)$ form a rational simplicial fan, with maximal cones the tilting cones, and
\[
\tilt Y\xrightarrow{\ \sim\ }\{\text{chambers in }V_{>0}\},
\qquad \E\longmapsto C(\E)^\circ.
\]
Two distinct basic tilting bundles are related by an irreducible mutation if and only if their cones share a facet.
Restriction preserves summands and identifies the bundles and their lattice fans:
\[
\ptilt Y\xrightarrow{\ \sim\ }\ptilt C,\qquad
\bigl(K_0(\per Y),\Sigma(Y)\bigr)\xrightarrow{\ \sim\ }
\bigl(K_0(\per C),\Sigma(C)\bigr).
\]
\end{IntroThm}

Restriction reduces explicit classification to the scheme-theoretic closed fibre $C$ (Proposition \ref{prop:closed-fibre-fan}). The chamber correspondence uses the reconstruction of tilting bundles in Theorem \ref{thm:tilting-torsion}, based on Fushimi's silting--simple-minded correspondence \cite[Theorem 3.6]{Fus26}.

The fan need not fill the positive-rank region. Theorem \ref{thm:complete-finite} measures exactly when it does:
\[
\begin{aligned}
\Sigma(Y)\text{ is complete in }V_{>0}
&\quad\Longleftrightarrow\quad \#\bigl(\tilt Y/\Pic(Y)\bigr)<\infty\\
&\quad\Longleftrightarrow\quad
\sup_{E\in\ptilt Y\colon\text{indecomposable}}\rk E<\infty.
\end{aligned}
\]
Here $\Pic(Y)$ acts by tensoring. Under these conditions, every basic tilting bundle can be reached from any other by mutation (Corollary \ref{cor:complete-mutation-connectivity}).

For ADE fibres, the reduction leads to an explicit root arrangement. Let $Z$ be a split ADE highest-root cycle of type $\Delta$ (Definition \ref{def:ade-cycle-scheme}). Here \emph{split} refers to the components and nodes of $Z_{\mathrm{red}}$. For $J\subseteq\Delta$, let $C_J$ be its Stein contraction along $J$. Write $\widetilde\Delta$ for the affine Dynkin diagram and $\operatorname{Cham}(\widetilde\Delta,J)$ for the chambers of its restricted Tits cone. Rank-one slicing identifies them with the $J$-alcoves of the affine ADE root arrangement restricted to $z_j=0$ for $j\in J$, as in \eqref{eq:restricted-affine-ade-arrangement}.

\begin{IntroThm}[ADE fans and Coxeter combinatorics]\label{intro:ade-alcoves}
Assume $C\cong C_J$ over $k$. Rank-and-degree coordinates on $\{\rk=1\}$ give the correspondences
\[
\begin{array}{c@{\quad\longleftrightarrow\quad}l}
\text{basic tilting bundles on }Y & \operatorname{Cham}(\widetilde\Delta,J)\text{ ($J$-alcoves)}\\
\text{nonzero pretilting cones} & \text{cones over nonempty faces of closed $J$-alcoves}\\
\text{irreducible mutation} & \text{facet adjacency}.
\end{array}
\]
For $J=\varnothing$, the underlying undirected mutation graph is isomorphic to
\[
\operatorname{Cay}\bigl(W_{\mathrm{aff}}(\Delta),S_{\mathrm{aff}}\bigr),
\]
where $S_{\mathrm{aff}}$ denotes the affine simple reflections.
\end{IntroThm}

These fans are complete, giving finiteness modulo line bundles and mutation connectivity. For $k$ algebraically closed and $R$ a complete local $k$-algebra, the result applies to projective crepant partial resolutions of Du Val singularities and, in characteristic zero, cDV singularities (Corollary \ref{cor:duval-cdv-tilting}).

Theorem \ref{thm:ade-cycle-fan} derives this description intrinsically from the curve; compare \cite[Theorems 5.2 and 5.26]{IWTC} for contracted preprojective algebras. General parabolic restrictions need not be Coxeter arrangements. On a reduced split nodal chain of $m\ge1$ projective lines, every indecomposable pretilting bundle is a line bundle, and there are $m!$ basic tilting bundles modulo line-bundle twists (Theorem \ref{thm:type-a-classification}).

We next ask what a module of sections remembers about the model. A module blowup makes the pullback of a module locally free modulo torsion \cite{Vil06,BCZ12,GI18}. This suggests using the globally generated tilting bundle containing $\O_Y$: its sections contain $R$, and its summands should encode contractions.

In Part \ref{part:auslander-mckay}, let $R$ be a normal Gorenstein complete local domain. A reflexive module $M$ is \emph{modifying} if $\End_R(M)$ is Cohen--Macaulay, and \emph{maximal modifying} if $M\oplus N$ modifying implies $N\in\add M$. Write $\MM R$ for the basic isomorphism classes of maximal modifying modules. An NCCR is $\End_R(M)$ for a nonzero modifying $M$ with finite global dimension. Let $\mathsf{ModGen}(R)$ be the poset of basic modifying generators ($R\in\add M$), ordered by additive inclusion. Let $\mathsf{Crep}_1(R)$ consist of normal integral projective birational models $\pi\colon Y\to X=\Spec R$ satisfying
\[
\dim\pi^{-1}(x)\le1\quad(x\in X),\qquad \pi^!\O_X\cong\O_Y,
\]
with $X$-morphisms as arrows. The duality condition is equivalent to $Y$ being Cohen--Macaulay with $\omega_Y\cong\O_Y$, and hence implies that $Y$ is Gorenstein (Lemma \ref{lem:crepancy-and-duality}). Crepant vanishing gives $\mathbb R\pi_*\O_Y\cong\O_X$ (Proposition \ref{prop:crepant-vanishing}), so $Y$ carries the unique basic globally generated tilting bundle $\P_Y$ containing $\O_Y$. Sections \ref{sec:general-reflexive-pushforward}--\ref{sec:general-auslander-mckay} prove the following.

\begin{IntroThm}[Reconstruction from canonical generators]\label{intro:reconstruction}
The canonical generators define a fully faithful functor recovering models and morphisms:
\[
\begin{gathered}
\mathsf{Crep}_1(R)^{\op}\lhook\joinrel\longrightarrow\mathsf{ModGen}(R),
\qquad Y\longmapsto M_Y:=\pi_*\P_Y,\\
Y\cong_X\operatorname{Bl}_{M_Y}X,\qquad
\Hom_X(Y,Z)\ne\varnothing\ \Longleftrightarrow\ \add M_Z\subseteq\add M_Y.
\end{gathered}
\]
\end{IntroThm}

Contractions from $Y$ correspond to summands $N$ of $M_Y$ containing $R$, with target $\operatorname{Bl}_N X$ up to $X$-isomorphism (Theorem \ref{thm:general-auslander-mckay-embedding}). They remove summands, so the faces of the canonical tilting cone containing $[\O_Y]$ encode contractions. Compare \cite[Corollary 4.4]{Wem18}.

For a contraction $h$, an idempotent quotient of $\End_Y(\P_Y)$ describes the coherent sheaves with $\mathbb Rh_*F=0$; its support records the exceptional image in $X$ (Theorem \ref{thm:contraction-null-category}); compare Donovan--Wemyss \cite{DW19contractions}.

For a characteristic-zero cDV singularity over an algebraically closed field and any fixed minimal model $\pi\colon Y\to\Spec R$, Corollaries \ref{cor:duval-cdv-tilting} and \ref{cor:cdv-modifying-realization} give
\[
\tilt Y\cong\operatorname{Cham}(\widetilde\Delta,J)
\overset{\Psi_Y}{\twoheadrightarrow}\MM R,
\qquad \Psi_Y(\E):=(\pi_*\E)^{**}.
\]
This is a bijection when $R$ is isolated. Over $\mathbb C$, compare \cite[Theorem 0.13]{IWTC} and the discussion after Corollary \ref{cor:cdv-modifying-realization}.

Theorem \ref{intro:reconstruction} leaves the essential image of the embedding to be determined. For the family $uv=f$, divisors $h\mid f$ give explicit rank-one ideals $(u,h)$, and chains of divisors give candidate generators and models. Part \ref{part:MM} shows that they exhaust the correspondence. We work with
\[
R=S[[u,v]]/(uv-f),\qquad
f=\varepsilon\prod_{\alpha=1}^{\ell}p_\alpha^{a_\alpha}\in\n\setminus\{0\},
\]
where $(S,\n)$ is a regular complete local ring of arbitrary dimension $d\ge1$, $\varepsilon\in S^\times$, and the $p_\alpha$ are pairwise nonassociated primes with $a_\alpha\ge1$. Since $\dim R=d+1$ has no upper bound, Part \ref{part:MM} treats higher-dimensional analogues of type $A$ Du Val surfaces and $cA_n$ threefolds, also in mixed characteristic. For $S=k[[x,y]]$, Iyama--Wemyss classified modifying generators by flags and, over $\mathbb C$, all modifying modules by rank-one reflexive twists \cite[Theorem 5.1, Lemma 5.22 and Proposition 5.24]{IW18}.

Put $a=(a_1,\ldots,a_\ell)$ and $n=\sum_\alpha a_\alpha$. A flag is a chain of integer exponent vectors
\[
\F=(b_1,\ldots,b_m),\qquad
0=b_0<b_1<\cdots<b_m<b_{m+1}=a,
\]
where $b<c$ means $b\le c$ coordinatewise and $b\ne c$. Let $\mathsf{Flag}(a)$ be their inclusion poset, including the empty flag. A flag is maximal when $m=n-1$. Define
\[
h_j:=\prod_\alpha p_\alpha^{b_{j,\alpha}}\quad(1\le j\le m),\qquad
T^\F:=R\oplus\bigoplus_{j=1}^{m}(u,h_j),
\]
and write $Y^\F=\operatorname{Bl}_{T^\F}X$ for the flag model. Theorem \ref{thm:auslander-mckay} gives the geometric classification.

\begin{IntroThm}[Flags and crepant models]\label{intro:flags}
The embedding of Theorem \ref{intro:reconstruction} is an equivalence, giving
\[
\xymatrix@C=3em@R=.8em{
\mathsf{Crep}_1(R)^{\op}\ar[r]^-{\sim}
 &\mathsf{ModGen}(R)\ar[r]^-{\sim}&\mathsf{Flag}(a)\\
Y^\F\ar@{|->}[r]&T^\F\ar@{|->}[r]&\F.
}
\]
\end{IntroThm}

Thus every model in $\mathsf{Crep}_1(R)$ is a flag blowup, and contractions are exactly flag deletions. A nonempty flag of length $m$ has a reduced split nodal chain of $m$ projective lines as its scheme-theoretic closed fibre; the empty flag gives $X$. The classification supplies the split-chain structure needed for the type $A$ tilting fan.

Rank-one summands provide the normalization that extends the classification from generators to all modifying modules over the same base $S$ (Theorems \ref{modrk1} and \ref{thm:modifying-classification}, and Corollary \ref{cor:mm-classification}).

\begin{IntroThm}[Modifying and maximal modifying modules]\label{intro:modifying-modules}
Every indecomposable modifying module has rank one. For a nonzero basic modifying module $M$ and an indecomposable summand $L$,
\[
M\cong(L\otimes_RT^\F)^{**}
\quad\text{for a unique flag }\F=(b_1,\ldots,b_m).
\]
Conversely, every rank-one reflexive $L$ and every flag give a basic modifying module. Moreover,
\[
M\text{ is maximal modifying}
\quad\Longleftrightarrow\quad m=n-1.
\]
\end{IntroThm}

The flag is unique relative to the chosen summand $L$; for generators, choosing $L=R$ recovers $T^\F$.

The NCCR criterion is geometric: an NCCR exists exactly when every factor hypersurface $S/(p_\alpha)$ is regular, equivalently $p_\alpha\notin\n^2$ for all $\alpha$. When this holds, precisely the maximal modifying modules give NCCRs (Theorem \ref{thm:nccr-criterion}). For $S=k[[x,y]]$, this recovers \cite[Theorem 5.1(3)]{IW18} via the CT--NCCR equivalence \cite[Corollary 5.9]{IW}.

The preceding classifications allow the flag and rank-one twist to vary. We now realize all these modules within one tilting fan. Fix any maximal-flag model $\pi\colon Y\to\Spec R$ and put
\[
\Psi_Y(U):=(\pi_*U)^{**}.
\]
Write $\modif R$ for the basic isomorphism classes of modifying modules, including $0$.

For $n\ge2$, Theorem \ref{thm:ade-cycle-fan} identifies $\tilt Y\cong\operatorname{Cham}(\widetilde A_{n-1},\varnothing)$.

\begin{IntroThm}[Realization on a fixed model]\label{intro:pushforward}
Every basic modifying module has a basic pretilting lift on this fixed $Y$, and $\Psi_Y^{-1}(\MM R)=\tilt Y$. The fibres of $\Psi_Y$ are exactly the orbits of an action of
\[
\Gamma\cong\prod_{\alpha:a_\alpha>1}W_{\mathrm{aff}}(A_{a_\alpha-1})
\]
on $\ptilt Y$, preserving $\tilt Y$. Thus $\Psi_Y$ induces the bijections
\[
\xymatrix@C=4em@R=1.4em{
\ptilt Y/\Gamma\ar[r]^-{\overline\Psi_Y}_-{\sim}&\modif R\\
\tilt Y/\Gamma\ar@{^{(}->}[u]\ar[r]^-{\overline\Psi_Y}_-{\sim}
 &\MM R\ar@{^{(}->}[u].
}
\]
For every $M\in\MM R$, the fibre $\Psi_Y^{-1}(M)$ is a $\Gamma$-torsor.
\end{IntroThm}

In Theorem \ref{intro:ade-alcoves}, the Weyl group describes alcove adjacency; here $\Gamma$ records the choices of lifts of one module. Its translations record line-bundle lifts of a rank-one summand, and its permutations record labels of repeated prime factors (Theorem \ref{thm:pushforward-fibres}). Both pushforward maps are bijective exactly when $f$ is squarefree (Corollary \ref{cor:pushforward-bijections}). For $n\ge2$, $\Psi_Y$ intertwines irreducible bundle mutation with both left and right modifying mutation (Theorem \ref{thm:pushforward-mutation}).

A single maximal-flag model, which may be singular, therefore supplies all maximal modification algebras through its tilting bundles, up to Morita equivalence. For any maximal-flag model $Y'$ and any maximal modification algebra $\Lambda$, there are $R$-linear equivalences
\[
D^b(\Coh Y)\simeq D^b(\Coh Y')\simeq D^b(\mod\Lambda),
\]
restricting to perfect complexes (Corollary \ref{cor:mma-derived-equivalence}). The three-dimensional MMA statement is already contained in \cite[Theorem 4.16]{IW}; the fixed-model realization here works in arbitrary dimension.

Finally, the wall categories retain geometry that the numerical fan forgets. For $n\ge2$, the tilting fan of a maximal-flag model depends only on the chain of $n-1$ curves, while the intersections of the prime hypersurfaces can vary. For a contraction $c$ deleting an entry of any nonempty flag, let $g,h$ be the adjacent factor blocks joined by the deletion. The coherent sheaves with $\mathbb Rc_*F=0$ form a category equivalent to $\mod A_c$; its algebra and closed-fibre semistable category satisfy
\[
A_c\cong S/(g,h),\qquad \W_c\simeq\fl S/(g,h)
\]
(Theorem \ref{thm:elementary-contraction-algebra} and Corollary \ref{cor:explicit-wall-category}); here $\fl$ denotes finite-length modules. For example, with $S=k[[x,y]]$ and $f=x(x+y^r)$, $r\ge1$, the fan is always of affine type $A_1$, whereas the wall algebra is $k[[y]]/(y^r)$. Its length recovers the intersection multiplicity $r$.

\subsection*{Conventions}
Base rings and fields are specified in each part; in particular, Part \ref{part:tilting bundles} does not assume that $R$ contains a field. The residue field of a base local ring is denoted by $k$, with no assumption on its characteristic or perfection unless stated. For a scheme point $x$, we write $\kappa(x)$ for its residue field, also viewed as a skyscraper sheaf when $x$ is closed. For a scheme $X$, we write $\vect X$ for the category of vector bundles on $X$.

We denote the projective morphism under consideration by $\pi\colon Y\to X$; its hypotheses are specified separately in each part.

Let $\U,\V$ be full subcategories of a triangulated or abelian category $\C$. We write $\U*\V$ for the full subcategory of $\C$ consisting of objects $M$ that fit into a triangle
\[
U\longrightarrow M\longrightarrow V\longrightarrow U[1]
\]
in the triangulated case, or a short exact sequence
\[
0\longrightarrow U\longrightarrow M\longrightarrow V\longrightarrow0
\]
in the abelian case, with $U\in\U$ and $V\in\V$.

For a finite module $M$ over a commutative noetherian ring $B$, put $\dim_BM:=\dim\operatorname{Supp}_B M$, with $\dim_B0=-\infty$.

For an object $M$ of an additive category, $\add M$ denotes the full subcategory of direct summands of finite direct sums of copies of $M$. Let $\U$ be a full subcategory of an additive category $\C$. A morphism $u\colon U\to M$ with $U\in\U$ is a \emph{right $\U$-approximation} of $M$ if every morphism $U'\to M$ with $U'\in\U$ factors through $u$. Left $\U$-approximations are defined dually. If $\C$ is linear over a commutative ring $B$ and $\Hom_{\C}(M,N)$ is finitely generated over $B$, a right $(\add M)$-approximation of $N$ can be chosen with source $M^{\oplus a}$; the dual statement holds for left approximations.

For an object $M$ of an abelian category $\A$, we write $\Fac M$ for the full subcategory of quotients of finite direct sums of copies of $M$:
\[
\Fac M:=\{N\in\A\mid M^{\oplus a}\twoheadrightarrow N\text{ for some }a\ge0\}.
\]

For a collection $\U$ of objects of a triangulated category, $\thick(\U)$ denotes the smallest full triangulated subcategory containing $\U$ and closed under direct summands. We write $\thick(M)$ when $\U=\{M\}$. For a full subcategory $\U$ of an abelian category $\A$, put
\[
\U^{\perp_i}:=\{M\in\A\mid\Ext_{\A}^i(U,M)=0\text{ for every }U\in\U\}
\qquad(i=0,1),
\]
where $\Ext_{\A}^0=\Hom_{\A}$. For a single object $U$, we write $U^{\perp_i}$ in place of $\{U\}^{\perp_i}$.

\subsection*{Use of AI}
The author used ChatGPT and Codex (OpenAI) to generate proof strategies for Theorems~\ref{thm:ade-cycle-fan} and~\ref{thm:general-auslander-mckay-embedding}. These tools were also used to assist with checking mathematical arguments, simplifying proofs, and revising the structure, language, and exposition of the manuscript. The author reviewed and revised the AI-generated material and takes full responsibility for the mathematical arguments and the final manuscript.

\section*{Acknowledgements}
He is grateful to Riku Fushimi, Wahei Hara, Osamu Iyama and Shu Nimura for fruitful discussions. This work was supported by the WINGS-FMSP program at the Graduate School of Mathematical Sciences, the University of Tokyo, and JSPS KAKENHI Grant Number JP25KJ0818.

\part{Tilting bundles in relative dimension one}\label{part:tilting bundles}

Throughout Part \ref{part:tilting bundles}, let $(R,\mathfrak{m})$ be a commutative noetherian complete local ring, put $k:=R/\m$, and let
\begin{equation}\label{eq:geometric-setup}
\pi\colon Y\longrightarrow X:=\Spec R
\qquad\text{with}\qquad
\mathbb{R}\pi_*\O_Y=\O_X,\quad \dim \pi^{-1}(\m)\le1
\end{equation}
be a projective morphism.
We write $C:=Y\times_Rk$ for the scheme-theoretic closed fibre. We impose no birationality assumption on $\pi$, no perfection assumption on $k$ and no reducedness assumption on $C$ unless explicitly stated.

We write $\per Y$ for the full subcategory of perfect complexes in $D^b(\Coh Y)$. Let
\[
\Coh_0Y:=\{M\in\Coh Y\mid\operatorname{Supp}M\subseteq \pi^{-1}(\m)\},
\]
and let $D^b_0(\Coh Y)$ denote the full subcategory of $D^b(\Coh Y)$ consisting of complexes whose cohomology sheaves belong to $\Coh_0Y$. Thus $\Coh_0Y$ is a Serre subcategory of $\Coh Y$, and its objects are precisely the coherent sheaves annihilated by some power of $\m\O_Y$.

For $\E\in\vect Y$, the subcategories $\Fac\E$ and $\E^{\perp_i}$ are taken in $\Coh Y$; in particular,
\[
\E^{\perp_i}=\{M\in\Coh Y\mid\Ext_Y^i(\E,M)=0\}\qquad(i=0,1).
\]
We call $\E$ a \emph{pretilting bundle} if $\Ext_Y^i(\E,\E)=0$ for all $i>0$, and a \emph{tilting bundle} if, in addition, $\thick(\E)=\per Y$.

A bundle is \emph{basic} if its indecomposable direct summands are pairwise nonisomorphic. Such decompositions exist and are unique up to isomorphism and order: endomorphism rings of vector bundles are finite algebras over the complete local ring $R$, so the category $\vect Y$ is Krull--Schmidt. We use the notation
\[
\begin{split}
\ptilt Y&:=\{\text{isomorphism classes of basic pretilting bundles on $Y$}\},\\
\tilt Y&:=\{\text{isomorphism classes of basic tilting bundles on $Y$}\}
\subseteq\ptilt Y.
\end{split}
\]
Properness and upper semicontinuity of fibre dimension imply that every fibre of $\pi$ has dimension at most one. Thus Van den Bergh's construction \cite[Proposition 3.2.5 and Corollary 3.2.8]{VdB04b}, which is valid over a general noetherian affine base, gives $\tilt Y\ne\emptyset$.

All Grothendieck groups of perfect complexes in this part are denoted explicitly by $K_0(\per Y)$. Set
\[
V:=K_0(\per Y)\otimes_{\mathbb Z}\mathbb R,
\qquad V_{>0}:=\{\theta\in V\mid\rk\theta>0\},
\qquad H:=\{\theta\in V\mid\rk\theta=1\}.
\]
The rank is well defined because $Y$ is connected, as proved below. Corollary \ref{cor:grothendieck-groups} shows that $K_0(\per Y)$ is free of finite rank, denoted by $n$. For a basic bundle or module $U$, let $|U|$ be its number of indecomposable summands. For a module-finite $R$-algebra $\La$, $\mod\La$ denotes the category of finitely generated right modules and $\fl\La$ its subcategory of modules of finite $R$-length.

For a tilting bundle $\E$, we use throughout this part
\[
\T_\E:=\Fac\E,\qquad\F_\E:=\E^{\perp_0},\qquad
\T_{\E,0}:=\T_\E\cap\Coh_0Y,\qquad\F_{\E,0}:=\F_\E\cap\Coh_0Y.
\]
Their HRS hearts, constructed in Section \ref{sec:tilting-torsion}, are denoted by $\H_\E:=\F_\E[1]*\T_\E$ and $\H_{\E,0}:=\F_{\E,0}[1]*\T_{\E,0}$.

The argument proceeds from torsion pairs and canonical generators (Section \ref{sec:tilting-torsion}) to the stability fan (Section \ref{sec:stability-fan}), then to completeness, finiteness and mutation connectivity (Section \ref{sec:completeness}). Section \ref{sec:ade-fibres} identifies the cycle-fibre fans and their ADE symmetries, with explicit bundles in type $A$.

\section{Tilting bundles and torsion pairs}\label{sec:tilting-torsion}

Recall that a \emph{torsion pair} in an abelian category $\A$ is a pair $(\T,\F)$ of full subcategories such that
\[
\Hom_{\A}(\T,\F)=0
\qquad\text{and}\qquad
\A=\T*\F.
\]
Thus every $M\in\A$ has a short exact sequence, called its torsion decomposition,
\[
0\longrightarrow T_M\longrightarrow M\longrightarrow F_M\longrightarrow0,
\qquad T_M\in\T,\quad F_M\in\F.
\]
In particular, $\T$ is closed under quotients and extensions, and $\F=\T^{\perp_0}$. The associated Happel--Reiten--Smal\o\ heart is
\begin{equation}\label{eq:hrs-heart}
\H:=\F[1]*\T
=\{C\in D^b(\A)\mid H^{-1}(C)\in\F,\ H^0(C)\in\T,\ H^i(C)=0\ (i\ne-1,0)\}.
\end{equation}
This is the heart of a bounded $t$-structure on $D^b(\A)$; see \cite[Example 2.7(2)]{PS20}. We will use the natural identification
\begin{equation}\label{eq:heart-ext-one}
\Ext^1_{\H}(A,B)\cong\Hom_{D^b(\A)}(A,B[1])\qquad(A,B\in\H).
\end{equation}
Here a \emph{progenerator} of $\H$ means a projective object $P$ such that every object of $\H$ is a quotient, in $\H$, of a finite direct sum of copies of $P$. In this case the projective objects of $\H$ are precisely the objects of $\add P$.

An abelian category is called a \emph{length category} if every object has a finite composition series. For a torsion pair in $\Coh_0Y$, we always regard its HRS tilt as a heart in $D^b_0(\Coh Y)$, using the standard $t$-structure with heart $\Coh_0Y$.

The aim of this section is to prove the following correspondence.

\begin{Thm}\label{thm:tilting-torsion}
Under the assumptions \eqref{eq:geometric-setup}, there are bijections between the following sets:
\begin{enumerate}
\item isomorphism classes of basic tilting bundles on $Y$;
\item torsion pairs $(\T,\F)$ in $\Coh Y$ such that $\T\ne0$ and $\F[1]*\T$ has a progenerator;
\item torsion pairs $(\T_0,\F_0)$ in $\Coh_0Y$ such that $\T_0\ne0$ and $\F_0[1]*\T_0$ is a length category.
\end{enumerate}
The maps from (1) to (2) and (3) are, respectively,
\[
\E\longmapsto(\Fac\E,\E^{\perp_0}),
\qquad
\E\longmapsto(\Fac\E\cap\Coh_0Y,\E^{\perp_0}\cap\Coh_0Y).
\]
The map from (2) to (3) is restriction:
\[
(\T,\F)\longmapsto(\T\cap\Coh_0Y,\F\cap\Coh_0Y).
\]
Its inverse is given by
\begin{equation}\label{eq:recover-torsion-class}
\T=\{M\in\Coh Y\mid M/\m M\in\T_0\},
\qquad\F=\T^{\perp_0}.
\end{equation}
The inverse from (2) to (1) takes the basic progenerator of $\F[1]*\T$. For corresponding $\E$ and $(\T,\F)$, we have
\[
\T=\Fac\E=\E^{\perp_1},\qquad
\F=\E^{\perp_0},\qquad
\proj(\F[1]*\T)=\add\E.
\]
\end{Thm}

The proof has two steps. Right exactness and Nakayama's lemma show that a torsion pair in (2) is determined by its restriction to the closed fibre. Fushimi's correspondence then reconstructs a perfect silting object from a heart in (3), and its rank forces it to be a vector bundle. Together with the tilting construction, these steps prove all three bijections.

\subsection{Cohomology and derived realization}\label{sec:derived-restrictions}

\begin{Lem}\label{lem:cohomological-dimension-one}
The scheme $Y$ is connected, and all morphism spaces in $D^b(\Coh Y)$ are finitely generated $R$-modules. For every $\V\in\vect Y$ and $M\in\Coh Y$, we have
\[
\Ext_Y^i(\V,M)=0\qquad(i>1).
\]
\end{Lem}
\begin{proof}
Since $\Gamma(Y,\O_Y)=R$ is local, it has no nontrivial idempotents, so $Y$ is connected. Properness and the local-to-global Ext spectral sequence give the asserted finiteness. The theorem on formal functions gives
\[
H^i(Y,N)=0\qquad(N\in\Coh Y,\ i>1),
\]
since every infinitesimal neighbourhood of the closed fibre has dimension at most one; see \cite[Tag 02V7]{stacks-project}. Now use $\Ext_Y^i(\V,M)=H^i(Y,\V^\vee\otimes M)$.
\end{proof}

The following observation will be used both for vector-bundle cohomology and for HRS cohomology.

\begin{Lem}[Right exactness and Nakayama's lemma]\label{lem:right-exact-nakayama}
Let $Q\colon\Coh Y\to\mod R$ be an $R$-linear right exact functor. For every $M\in\Coh Y$, there is a natural isomorphism
\[
Q(M/\m M)\cong Q(M)/\m Q(M).
\]
In particular, $Q(M)=0$ if and only if $Q(M/\m M)=0$.
\end{Lem}
\begin{proof}
Choose generators $a_1,\ldots,a_s$ of $\m$ and apply $Q$ to
\[
M^{\oplus s}\xrightarrow{(a_1,\ldots,a_s)}M\longrightarrow M/\m M\longrightarrow0.
\]
The isomorphism follows from right exactness and $R$-linearity, and the last assertion is Nakayama's lemma.
\end{proof}

\begin{Cor}\label{lem:closed-fibre-detects-ext}
For $\V\in\vect Y$ and $M\in\Coh Y$, we have
\[
\Ext_Y^1(\V,M)=0
\quad\Longleftrightarrow\quad
\Ext_Y^1(\V,M/\m M)=0.
\]
\end{Cor}
\begin{proof}
By Lemma \ref{lem:cohomological-dimension-one}, the functor $\Ext_Y^1(\V,-)\colon\Coh Y\to\mod R$ is right exact. Apply Lemma \ref{lem:right-exact-nakayama}.
\end{proof}

We also record how line bundles are determined by the closed fibre. This gives the Picard lattice of translations in Section \ref{sec:completeness} and the multidegree description of $\Pic(Y)$ in Section \ref{sec:ade-fibres}.

\begin{Lem}[Picard restriction]\label{lem:picard-restriction}
Restriction gives an isomorphism and a surjection
\[
\Pic(Y)\xrightarrow{\ \sim\ }\Pic(C)\twoheadrightarrow\Pic(C_{\mathrm{red}}).
\]
\end{Lem}
\begin{proof}
Use the formal lifting argument of \cite[proof of Lemma 3.4.3]{VdB04b}. Put $I:=\m\O_Y$ and $C_q:=Y\times_R R/\m^{q+1}$. Lemma \ref{lem:cohomological-dimension-one} gives $H^1(C,\O_C)=0$. Each successive ideal $J_q:=I^{q+1}/I^{q+2}$ is a quotient of a finite sum of $\O_C$, so $H^1(C,J_q)=H^2(C,J_q)=0$. The square-zero unit sequence \cite[Tag 0C6R]{stacks-project} therefore identifies the Picard groups of all $C_q$ and makes their global-unit maps surjective. Grothendieck existence \cite[Tag 088E]{stacks-project} then gives $\Pic(Y)\cong\Pic(C)$. This argument uses neither a coefficient field nor a residue-field hypothesis. Finally, the same unit sequence for the nilradical filtration of $C$, with obstruction groups in $H^2$, gives the asserted surjection.
\end{proof}

A perfect complex $P$ is \emph{silting} if
\[
\Hom_{D^b(\Coh Y)}(P,P[i])=0\quad(i>0),
\qquad \thick(P)=\per Y.
\]
It is a \emph{tilting complex} if the vanishing holds for every $i\ne0$. Properness and completeness of $R$ make $\per Y$ Krull--Schmidt, so these objects admit basic representatives.

\begin{Lem}[Derived realization and restrictions]\label{lem:tilting-derived-realization}
Let $T\in\per Y$ be a tilting complex and put $\Lambda:=\End_{D^b(\Coh Y)}(T)$. Then
\[
\Phi_T:=\RHom_Y(T,-)\colon
D^b(\Coh Y)\xrightarrow{\ \sim\ }D^b(\mod\Lambda)
\]
is an $R$-linear triangle equivalence and restricts to equivalences
\begin{equation}\label{eq:tilting-derived-restrictions}
D^b_0(\Coh Y)\xrightarrow{\ \sim\ }D^b_{\fl}(\mod\Lambda),
\qquad
\per Y\xrightarrow{\ \sim\ }\per\Lambda.
\end{equation}
Here $\per\Lambda=K^b(\proj\Lambda)$ and
\[
D^b_{\fl}(\mod\Lambda):=
\{M\in D^b(\mod\Lambda)\mid H^i(M)\in\fl\Lambda\text{ for every }i\}.
\]
\end{Lem}
\begin{proof}
Derived Morita theory \cite{Kel94} applied directly to the compact generator $T$ gives an $R$-linear equivalence on unbounded derived categories. The perfect restriction follows from $\Phi_T(T)=\Lambda$. Properness makes $\Lambda$ finite over $R$, and $\Phi_T$ carries $D^b(\Coh Y)$ into $D^b(\mod\Lambda)$. For the inverse, bounded-above resolutions by finitely generated projectives give bounded-above coherent cohomology, using finite truncations and the fixed amplitude of $T$. Boundedness below is detected by $T$ \cite[Tag 0GEQ]{stacks-project}, proving the bounded equivalence.

Both support subcategories consist of bounded complexes annihilated, as derived objects, by some power of $\m$: this follows from the corresponding condition on cohomology by induction on truncation triangles. The $R$-linearity of $\Phi_T$ therefore gives the remaining restriction.
\end{proof}

\begin{Cor}[Grothendieck groups]\label{cor:grothendieck-groups}
The groups
\[
K_0(\per Y),\qquad
K_0(D^b_0(\Coh Y))\cong K_0(\Coh_0Y)
\]
are finitely generated free abelian groups of the same rank $n$. For any basic tilting complex $T=\bigoplus_{i=1}^nT_i$, the classes $[T_i]$ form a basis of the first group, and the inverse images of the simple $\End(T)$-modules form a basis of the second. Every length heart of a bounded $t$-structure on $D^b_0(\Coh Y)$ has exactly $n$ isomorphism classes of simple objects.
\end{Cor}
\begin{proof}
Choose a basic tilting complex $T$. The algebra $\Lambda=\End(T)$ is finite over the complete local ring $R$, hence semiperfect. The indecomposable projectives form a basis of $K_0(\per\Lambda)$, and their simple tops form a basis of $K_0(\fl\Lambda)\cong K_0(D^b_{\fl}(\mod\Lambda))$. Apply Lemma \ref{lem:tilting-derived-realization}. A bounded $t$-structure identifies the Grothendieck group of its heart with that of the ambient category; for a length heart, the simple classes form a basis. This proves the last assertion as well.
\end{proof}

\subsection{Progenerators and restriction to the closed fibre}

We first prove that restriction determines every torsion pair whose HRS heart has a progenerator. No perfectness assumption on that progenerator is needed.

\begin{Lem}\label{lem:progenerator-restriction}
Let $(\T,\F)$ be a torsion pair in $\Coh Y$ with $\T\ne0$, and suppose that $\H:=\F[1]*\T$ has a progenerator. Then
\[
(\T_0,\F_0):=(\T\cap\Coh_0Y,\F\cap\Coh_0Y)
\]
is a torsion pair in $\Coh_0Y$ with $\T_0\ne0$, and $\F_0[1]*\T_0$ is a length category. Moreover, $(\T,\F)$ is recovered by \eqref{eq:recover-torsion-class}. In particular, restriction is injective on the torsion pairs in Theorem \ref{thm:tilting-torsion}(2).
\end{Lem}
\begin{proof}
Let $P$ be a progenerator and put $A:=\End_{\H}(P)$. Lemma \ref{lem:cohomological-dimension-one} makes $A$ finite over $R$, and the progenerator gives an $R$-linear Morita equivalence
\[
\Hom_{\H}(P,-)\colon\H\xrightarrow{\ \sim\ }\mod A;
\]
compare \cite[Lemma 3.2.1]{VdB04b}. Torsion decompositions restrict to the Serre subcategory $\Coh_0Y$. By the $\m$-power annihilation criterion in the proof of Lemma \ref{lem:tilting-derived-realization}, the Morita equivalence restricts to
\[
\H\cap D^b_0(\Coh Y)=\F_0[1]*\T_0
\xrightarrow{\ \sim\ }\fl A.
\]
Thus this heart is a length category. A nonzero sheaf in $\T$ has a quotient $\kappa(y)$ for some closed point $y$. Properness gives $y\in C$, so $\kappa(y)\in\T_0$ and $\T_0\ne0$.

For $M\in\Coh Y$, HRS cohomology is concentrated in degrees $0,1$, with $H^1_{\H}(M)=F_M[1]$. Consequently
\[
Q(M):=\Hom_{\H}\bigl(P,H^1_{\H}(M)\bigr)
\]
is an $R$-linear right exact functor to $\mod R$: $H^1_{\H}$ is right exact on $\Coh Y$, and $P$ is projective. Since $P$ is a generator, $Q(M)=0$ if and only if $M\in\T$. Lemma \ref{lem:right-exact-nakayama} now gives
\[
M\in\T\quad\Longleftrightarrow\quad
Q(M/\m M)=0\quad\Longleftrightarrow\quad
M/\m M\in\T_0.
\]
This proves \eqref{eq:recover-torsion-class}, since $\F=\T^{\perp_0}$.
\end{proof}

We next construct the torsion pair associated with a tilting bundle by transporting the standard $t$-structure.

\begin{Prop}\label{prop:bundle-to-torsion}
Let $\E$ be a tilting bundle on $Y$. Then
\[
\T_{\E}:=\E^{\perp_1}=\Fac\E,
\qquad\F_{\E}:=\E^{\perp_0}
\]
form a torsion pair in $\Coh Y$, and $\E$ is a progenerator of $\H_{\E}:=\F_{\E}[1]*\T_{\E}$. Writing $\La:=\End_Y(\E)$, the equivalence $\Phi_{\E}$ of Lemma \ref{lem:tilting-derived-realization} restricts to $\H_{\E}\simeq\mod\La$.
\end{Prop}
\begin{proof}
Lemma \ref{lem:cohomological-dimension-one} gives
\[
\Phi_{\E}(\Coh Y)\subseteq D^{[0,1]}(\mod\La).
\]
Using bounded truncations and orthogonality, we obtain
\[
D^b(\Coh Y)^{\le-1}
\subseteq\Phi_{\E}^{-1}\bigl(D^b(\mod\La)^{\le0}\bigr)
\subseteq D^b(\Coh Y)^{\le0}.
\]
The characterization of intermediate $t$-structures \cite[Lemma 1.1.2]{Pol07} therefore identifies the pulled-back heart with the HRS tilt for the pair $(\E^{\perp_1},\E^{\perp_0})$. Since $\Phi_{\E}(\E)=\La$, the bundle $\E$ is its progenerator.

It remains to identify the torsion class with $\Fac\E$. Since $\E\in\T_{\E}$ and torsion classes are closed under quotients, $\Fac\E\subseteq\T_{\E}$. Conversely, for $T\in\T_{\E}$, take a right $(\add\E)$-approximation $\E^{\oplus a}\to T$ in $\H_{\E}$. It is an epimorphism because $\E$ is a progenerator. Its kernel lies in $\H_{\E}$, so the standard cohomology sequence makes the map an epimorphism also in $\Coh Y$. Thus $T\in\Fac\E$.
\end{proof}

\begin{Cor}\label{lem:tilting-restricts-to-length}
Let $\E$ be a tilting bundle and put $\La:=\End_Y(\E)$. Then $(\T_{\E,0},\F_{\E,0})$ is a torsion pair in $\Coh_0Y$ with $\T_{\E,0}\ne0$, and
\[
\Hom_{D^b(\Coh Y)}(\E,-)\colon
\H_{\E,0}=\F_{\E,0}[1]*\T_{\E,0}\xrightarrow{\ \sim\ }\fl\La.
\]
In particular, $\H_{\E,0}$ is a length category, and
\begin{equation}\label{eq:finite-heart-orthogonality}
\Hom_{D^b(\Coh Y)}(\E,Z[i])=0
\qquad(Z\in\H_{\E,0},\ i\ne0).
\end{equation}
\end{Cor}
\begin{proof}
Apply Proposition \ref{prop:bundle-to-torsion} and Lemma \ref{lem:progenerator-restriction}. The equivalence and the vanishing are the restrictions of $\Phi_{\E}$ to the indicated heart.
\end{proof}

\subsection{Reconstruction from length hearts}

The algebraic reconstruction is supplied by Fushimi's silting--simple-minded correspondence over a complete local base.

\begin{Lem}\label{lem:length-heart-duality}
Let $\K$ be the heart of a bounded $t$-structure on $D^b_0(\Coh Y)$, and suppose that $\K$ is a length category. Then there is a unique basic silting object $P\in\per Y$ such that
\begin{equation}\label{eq:dual-silting-object}
\Hom_{D^b(\Coh Y)}(P,C[i])=0
\qquad(C\in\K,\ i\ne0).
\end{equation}
For this $P$, the functor $\Hom_{D^b(\Coh Y)}(P,-)$ on $\K$ is exact and faithful.
\end{Lem}
\begin{proof}
Via Lemma \ref{lem:tilting-derived-realization}, apply \cite[Theorem 3.6 and Proposition 2.13(1)]{Fus26} to the finite simple-minded collection of simple objects of $\K$; finiteness is Corollary \ref{cor:grothendieck-groups}.
\end{proof}

\begin{Prop}\label{prop:length-heart-to-bundle}
Let $(\T_0,\F_0)$ be a torsion pair in $\Coh_0Y$ with $\T_0\ne0$, and suppose that $\K:=\F_0[1]*\T_0$ is a length category. Then there is a unique basic tilting bundle $\E$ such that
\[
\T_0=\Fac\E\cap\Coh_0Y,
\qquad\F_0=\E^{\perp_0}\cap\Coh_0Y.
\]
\end{Prop}
\begin{proof}
Let $P\in\per Y$ be the basic silting object supplied by Lemma \ref{lem:length-heart-duality}. Every closed point $y\in Y$ lies in $C$, and the simple sheaf $\kappa(y)$ belongs to exactly one of $\T_0$ and $\F_0$. Put $\epsilon_y=0$ in the first case and $\epsilon_y=1$ in the second. Since $\kappa(y)[\epsilon_y]\in\K$, we have
\[
\Hom_{D^b(\Coh Y)}(P,\kappa(y)[i])=0\quad(i\ne\epsilon_y),
\qquad
\Hom_{D^b(\Coh Y)}(P,\kappa(y)[\epsilon_y])\ne0.
\]
Choose a minimal bounded complex $Q_y^\bullet$ of finite free $\O_{Y,y}$-modules representing $P_y$. Minimality gives
\[
\Hom_{D^b(\Coh Y)}(P,\kappa(y)[i])
\cong\Hom_{\O_{Y,y}}(Q_y^{-i},\kappa(y)).
\]
Thus $P_y\cong\O_{Y,y}^{\oplus r_y}[\epsilon_y]$ with $r_y>0$.

The virtual rank of a perfect complex is locally constant, and $Y$ is connected. Hence
\begin{equation}\label{eq:silting-rank-sign}
\rk P=(-1)^{\epsilon_y}r_y
\qquad\text{for every closed point }y.
\end{equation}
A nonzero object of $\T_0$ has a closed-point quotient $\kappa(y)\in\T_0$, so $\rk P>0$. Equation \eqref{eq:silting-rank-sign} forces $\epsilon_y=0$ at every closed point. Every point of $Y$ specializes to a closed point, so the same local form holds everywhere. Thus $P$ is a vector bundle $\E$. Its silting properties already give generation and positive self-extension vanishing, so $\E$ is a basic tilting bundle.

For $T\in\T_0$ and $F\in\F_0$, equation \eqref{eq:dual-silting-object} gives
\[
\Ext_Y^1(\E,T)=0,
\qquad\Hom_Y(\E,F)=0.
\]
Thus $\T_0\subseteq\T_{\E,0}$ and $\F_0\subseteq\F_{\E,0}$. Both pairs are torsion pairs by Corollary \ref{lem:tilting-restricts-to-length}, so these inclusions force equality. Any other basic tilting bundle inducing this pair satisfies \eqref{eq:dual-silting-object} by the same corollary. Uniqueness follows from Lemma \ref{lem:length-heart-duality}.
\end{proof}

\begin{proof}[Proof of Theorem \ref{thm:tilting-torsion}]
Proposition \ref{prop:bundle-to-torsion} gives a map $a\colon(1)\to(2)$. Lemma \ref{lem:progenerator-restriction} gives an injective restriction map $r\colon(2)\to(3)$ and the recovery formula \eqref{eq:recover-torsion-class}. Corollary \ref{lem:tilting-restricts-to-length} and Proposition \ref{prop:length-heart-to-bundle} show that $r\circ a\colon(1)\to(3)$ is bijective. It follows that $a$ and $r$ are both bijective.

For the pair corresponding to $\E$, Proposition \ref{prop:bundle-to-torsion} makes $\E$ a progenerator of its heart. Any two progenerators have the same additive closure, namely the projective objects. Krull--Schmidt therefore identifies $\E$ with its unique basic progenerator. This proves the remaining assertions about the inverse and the projective objects.
\end{proof}

\begin{Rem}
The use of basic bundles is necessary for injectivity on isomorphism classes, since $\Fac\E=\Fac(\E\oplus\E)$. Equivalently, Theorem \ref{thm:tilting-torsion} may be stated for all tilting bundles modulo the relation $\E\sim\E'$ if $\add\E=\add\E'$.

The nonzero conditions on the torsion classes exclude shifted bundles. For example, when $Y=\Spec R$, the pair $(0,\Coh Y)$ has heart $(\Coh Y)[1]$ with progenerator $\O_Y[1]$. Likewise, $\Coh_0Y=\fl R$, and the pair $(0,\Coh_0Y)$ has the length heart $(\fl R)[1]$. These pairs correspond to the shifted object $\O_Y[1]$ rather than a vector bundle in degree zero.
\end{Rem}

\subsection{Canonical generators and regularity}\label{sec:canonical-generators}

We finish this section by recording the form of Van den Bergh's construction used in the birational applications. The results here require only \eqref{eq:geometric-setup}.

\begin{Prop}[Canonical tilting generator]\label{prop:canonical-tilting-general}
There is a unique basic tilting bundle $\P_Y$, up to isomorphism, which is globally generated and satisfies $\O_Y\in\add\P_Y$. Its determinant is $\pi$-ample, and
\begin{equation}\label{eq:canonical-projective-bundles}
\add\P_Y
=\{V\in\vect Y\mid V\text{ is globally generated and }H^1(Y,V^\vee)=0\}.
\end{equation}
A bundle $V$ in this category has $\add V=\add\P_Y$ if and only if $\O_Y\in\add V$ and $\det V$ is $\pi$-ample.
\end{Prop}
\begin{proof}
Van den Bergh constructs a projective generator of ${}^{-1}\!\operatorname{Per}(Y/X)$ and identifies its projective objects and generators by precisely these conditions \cite[Propositions 3.2.5--3.2.7]{VdB04b}. Take its basic representative $\P_Y$. Corollary 3.2.8 of loc.\ cit. makes it tilting; the generator criterion also gives its ample determinant and \eqref{eq:canonical-projective-bundles}.

If $\E$ is another globally generated tilting bundle with $\O_Y\in\add\E$, then pretilting gives $H^1(Y,\E^\vee)=0$. Hence $\E\in\add\P_Y$. Under the equivalence supplied by $\P_Y$, it is a projective module generating the category of perfect complexes, and therefore is a projective generator. Thus $\add\E=\add\P_Y$, proving uniqueness for basic bundles by Krull--Schmidt.
\end{proof}

A tilting equivalence also detects regularity; compare \cite[Corollary 4.2(1)]{IW14} over a field. The required bound on global dimension is the simple-module criterion for semiperfect rings.

\begin{Lem}\label{lem:tilting-regularity}
Let $B$ be a noetherian complete local ring and let $Z$ be a nonempty finite-dimensional noetherian scheme projective over $B$. Suppose that a tilting bundle $\E$ induces an equivalence
\[
D^b(\Coh Z)\simeq D^b(\mod\Lambda),\qquad
\Lambda:=\End_Z(\E),
\]
restricting to $\per Z\simeq\per\Lambda$. Then $\gl\Lambda<\infty$ if and only if $Z$ is regular.
\end{Lem}
\begin{proof}
The algebra $\Lambda$ is finite over $B$, hence left and right noetherian and semiperfect. For representatives $S_1,\ldots,S_r$ of its simple right modules, the standard simple-module criterion gives
\[
\gl\Lambda=\max_{1\le i\le r}\pd_\Lambda S_i.
\]
Consequently $\gl\Lambda<\infty$ is equivalent to $D^b(\mod\Lambda)=\per\Lambda$. Under the given equivalence, this is equivalent to $D^b(\Coh Z)=\per Z$, which in turn is equivalent to regularity of $Z$.
\end{proof}

\section{Numerical stability and the tilting fan}\label{sec:stability-fan}

We use numerical stability on $\Coh_0Y$ in the sense of King. After constructing numerical torsion pairs, we show that bricks on the closed fibre detect walls and TF equivalence. We then describe the semistable category of each nonzero pretilting face and recover its closed cone from the numerical classes. These results identify relative interiors with TF equivalence classes and give the fan theorem. The length hearts of Section \ref{sec:tilting-torsion} identify the full dimensional cones with chambers. We finish by recovering pretilting bundles from their classes in $K_0(\per Y)$ and constructing mutation across facets.

For $P\in\per Y$ and $M\in D^b_0(\Coh Y)$, define
\begin{equation}\label{eq:euler-pairing}
\langle[P],[M]\rangle
:=\sum_{i\in\mathbb Z}(-1)^i\ell_R\Hom_{D^b(\Coh Y)}(P,M[i]).
\end{equation}
The sum is finite, and each summand has finite $R$-length. This gives an Euler pairing on the Grothendieck groups. For $\theta\in V$, we write $\theta(M):=\langle\theta,[M]\rangle$.

\subsection{Completions and numerical classes}

We first supply the completion results needed below. The two extremal completions will also describe the torsion classes on opposite sides of a wall.

\begin{Prop}\label{Bongcomp}
Every basic pretilting bundle $U$ is a summand of a basic tilting bundle. If $U\ne0$, there are tilting completions $\E_{\min}$ and $\E_{\max}$ such that
\begin{equation}\label{eq:extremal-completions}
\Fac\E_{\min}=\Fac U,
\qquad
\Fac\E_{\max}=U^{\perp_1}.
\end{equation}
Thus, for every tilting completion $\E$ of $U$,
\[
\Fac\E_{\min}\subseteq\Fac\E\subseteq\Fac\E_{\max}.
\]
\end{Prop}
\begin{proof}
The assertion for $U=0$ follows from the existence of a tilting bundle. Suppose $U\ne0$; connectedness makes its rank positive. Fix a tilting bundle $G_0$ and an ample line bundle $\L$.

For $N\gg0$, the bundle $(G_0\otimes\L^{-N})^\vee\otimes U$ is globally generated. Put $G:=G_0\otimes\L^{-N}$ and take a left $(\add U)$-approximation $G\to U_0$. It is a monomorphism with locally free cokernel: at each point the global maps $G\to U$ span the full space of fibre maps, and these maps separate vectors in the fibre of $G$. We obtain
\[
0\longrightarrow G\longrightarrow U_0\longrightarrow W\longrightarrow0,
\qquad U_0\in\add U,\quad W\in\vect Y.
\]
The approximation property gives $\Ext_Y^1(W,U)=0$. The same sequence and cohomological dimension one give $\Ext_Y^1(U,W)=0$ and then $\Ext_Y^1(W,W)=0$. Moreover, $G\in\thick(U\oplus W)$. Hence $U\oplus W$ is tilting. Its basic representative $\E_{\min}$ satisfies $\Fac\E_{\min}=\Fac U$, since $W$ is a quotient of $U_0$.

Dually, choose $G':=G_0\otimes\L^N$ with $U^\vee\otimes G'$ globally generated, and take a right $(\add U)$-approximation $U_1\to G'$. It is surjective and has a locally free kernel $W'$:
\[
0\longrightarrow W'\longrightarrow U_1\longrightarrow G'\longrightarrow0.
\]
The approximation property gives $\Ext_Y^1(U,W')=0$. Applying $\Hom_Y(-,U)$ and $\Hom_Y(-,W')$ gives $\Ext_Y^1(W',U)=0=\Ext_Y^1(W',W')$. Thus $U\oplus W'$ is tilting. If $M\in U^{\perp_1}$, this sequence gives $\Ext_Y^1(W',M)=0$, since $\Ext_Y^2(G',M)=0$. Its basic representative $\E_{\max}$ therefore satisfies $\E_{\max}^{\perp_1}=U^{\perp_1}$. Proposition \ref{prop:bundle-to-torsion} proves \eqref{eq:extremal-completions}. The final inclusions follow from $U\in\add\E$.
\end{proof}

\begin{Lem}\label{lem:euler-dual-bases}
Let $\E=\bigoplus_{i=1}^n\E_i$ be a basic tilting bundle. Then $[\E_1],\ldots,[\E_n]$ form a $\mathbb Z$-basis of $K_0(\per Y)$. The length heart
\[
\H_{\E,0}:=(\E^{\perp_0}\cap\Coh_0Y)[1]*(\Fac\E\cap\Coh_0Y)
\]
has simple objects $S_1,\ldots,S_n$, indexed so that
\begin{equation}\label{eq:dual-bases}
\langle[\E_i],[S_j]\rangle=\delta_{ij}c_j,
\qquad c_j:=\ell_R\End_{\H_{\E,0}}(S_j)>0.
\end{equation}
Each $S_j$ is either a sheaf in $\Fac\E\cap\Coh_0Y$ or the shift by $[1]$ of a sheaf in $\E^{\perp_0}\cap\Coh_0Y$.
\end{Lem}
\begin{proof}
Corollary \ref{cor:grothendieck-groups} gives the basis assertion.

Corollary \ref{lem:tilting-restricts-to-length} identifies $\H_{\E,0}$ with $\fl\End_Y(\E)$. Its simples are also simple in the global heart $\H_\E$, since $\fl\End_Y(\E)$ is a Serre subcategory of $\mod\End_Y(\E)$. Pairing the indecomposable projectives in the global heart with their simple tops gives \eqref{eq:dual-bases}. For a simple $S\in\H_{\E,0}$, the standard cohomology triangle is a short exact sequence in this heart,
\[
0\longrightarrow H^{-1}(S)[1]\longrightarrow S\longrightarrow H^0(S)\longrightarrow0.
\]
Simplicity makes one of the two outer terms zero.
\end{proof}

\begin{Cor}\label{cor:number-of-summands}
A basic pretilting bundle $U$ has at most $n$ summands, and it is tilting if and only if $|U|=n$.
\end{Cor}
\begin{proof}
Complete $U$ by Proposition \ref{Bongcomp} and use Corollary \ref{cor:grothendieck-groups}.
\end{proof}

For a basic pretilting bundle $U=\bigoplus_{i=1}^rU_i$, set
\begin{equation}\label{eq:pretilting-cone}
C(U):=\sum_{i=1}^r\mathbb R_{\ge0}[U_i],
\qquad C^\circ(U):=\sum_{i=1}^r\mathbb R_{>0}[U_i].
\end{equation}
Thus $C^\circ(U)$ is the relative interior of the simplicial cone $C(U)$; we put $C(0)=C^\circ(0)=\{0\}$.

\begin{Cor}\label{cor:cone-from-heart}
For a tilting bundle $\E$, writing $\T_{\E,0}:=\Fac\E\cap\Coh_0Y$ and $\F_{\E,0}:=\E^{\perp_0}\cap\Coh_0Y$, we have
\[
C(\E)=\{\theta\in V\mid\theta(T)\ge0\ (T\in\T_{\E,0}),\ \theta(F)\le0\ (F\in\F_{\E,0})\}.
\]
Replacing the inequalities by strict inequalities on nonzero objects gives $C^\circ(\E)$.
\end{Cor}
\begin{proof}
Under the Euler pairing, the right-hand side is the heart cone of $\H_{\E,0}$ \cite[Definition 2.2]{BPPW24}. It suffices to test its simples. In the basis of Lemma \ref{lem:euler-dual-bases}, these tests are exactly nonnegativity of the coefficients. The strict version follows in the same way.
\end{proof}

For a nonzero basic pretilting bundle $U=\bigoplus_iU_i$, put
\[
\Delta(U):=C(U)\cap H
=\Conv\{[U_i]/\rk U_i\}_i.
\]
This is a compact simplex of dimension $|U|-1$, hence of dimension $n-1$ when $U$ is tilting. We will use this rank-one slice both for the type $A$ classification and for the completeness criterion.

\subsection{Numerical torsion pairs on the closed fibre}

The numerical torsion pairs in the module setting of \cite[Proposition 2.11]{Asa21} live in a length category. The category $\Coh_0Y$ need not be a length category: sheaves supported on a curve can have infinite descending chains of subsheaves. We first express the Euler pairing using generic multiplicities and Euler characteristics. For a fixed sheaf, these invariants give the bounds needed to construct torsion decompositions for real parameters. The sum of the generic multiplicities will also bound the lengths of semistable objects and replace module dimension in the brick argument of Theorem \ref{thm:tf-brick-walls}.

Let $Z_1,\ldots,Z_s$ be the one-dimensional irreducible components of $C_{\mathrm{red}}$, with generic points $\eta_i$. For $\xi\in K_0(\per Y)$, put
\[
d_i(\xi):=\deg_k\bigl(\det(\xi)|_{Z_i}\bigr),
\]
and extend $d_i$ linearly to $V$. Degrees are taken over $k$: for a line bundle $L$ on $Z_i$, we use $\deg_k L=\chi_k(Z_i,L)-\chi_k(Z_i,\O_{Z_i})$. For $M\in\Coh_0Y$, set
\[
\chi_R(M):=\ell_R H^0(Y,M)-\ell_R H^1(Y,M)
.
\]
The \emph{generic multiplicities} of $M$ are the lengths
\begin{equation}\label{eq:generic-multiplicity}
\rho_i(M):=\ell_{\O_{Y,\eta_i}}(M_{\eta_i})
\qquad(1\le i\le s).
\end{equation}

\begin{Lem}[Riemann--Roch and numerical coordinates]\label{lem:numerical-devissage}
The lengths \eqref{eq:generic-multiplicity} are finite and define additive functions $\rho_i\colon K_0(\Coh_0Y)\to\mathbb Z$, nonnegative on sheaves. They vanish simultaneously exactly when the support has dimension at most zero. For every $\theta\in V$ and $M\in\Coh_0Y$, we have
\begin{equation}\label{eq:rank-degree-euler}
\theta(M)=(\rk\theta)\chi_R(M)-\sum_{i=1}^s d_i(\theta)\rho_i(M).
\end{equation}
Moreover, the rank-and-degree map
\[
(\rk,d_1,\ldots,d_s)\colon V\longrightarrow\mathbb R^{s+1}
\]
is injective.
\end{Lem}
\begin{proof}
Since $M$ is annihilated by a power of $\m$ and $\eta_i$ is a generic point of $C$, the stalk $M_{\eta_i}$ has zero-dimensional support over $\O_{Y,\eta_i}$, hence finite length. Additivity, nonnegativity and the assertion about support follow from local length.

By additivity and the $\m$-adic filtration, it suffices to prove \eqref{eq:rank-degree-euler} for $M\in\Coh C$. Then $\chi_R(M)=\chi_k(C,M)$ and $\rho_i(M)=\ell_{\O_{C,\eta_i}}(M_{\eta_i})$. For a bundle $\V\in\vect Y$ of rank $r$, Riemann--Roch on $C/k$ gives
\[
\chi_k\bigl(C,(\V|_C)^\vee\otimes M\bigr)
=r\chi_k(C,M)-\sum_{i=1}^s d_i(\V)\rho_i(M).
\]
This form, valid also for nonreduced $C$ over an arbitrary field, follows from the integral-curve formula and coherent d\'evissage; see \cite[Tags 0AYV, 0DJ5, and 01YI]{stacks-project}. It is \eqref{eq:rank-degree-euler} for $\theta=[\V]$, and linearity gives the formula for every $\theta\in V$.

Finally, a class with rank and all degrees zero pairs to zero with every sheaf in $\Coh_0Y$. The dual bases of Lemma \ref{lem:euler-dual-bases} therefore force it to be zero, proving injectivity.
\end{proof}

\begin{Cor}[Rank and determinant]\label{cor:rank-determinant}
Rank and determinant give an isomorphism
\[
(\rk,\det)\colon K_0(\per Y)\xrightarrow{\ \sim\ }\mathbb Z\oplus\Pic(Y),
\]
with inverse
\[
(r,L)\longmapsto(r-1)[\O_Y]+[L].
\]
\end{Cor}
\begin{proof}
A class in the kernel has rank and all degrees zero, so Lemma \ref{lem:numerical-devissage} makes it zero in $V$. The torsion-freeness in Corollary \ref{cor:grothendieck-groups} makes it zero in $K_0(\per Y)$, proving injectivity. The proposed inverse supplies surjectivity.
\end{proof}

We put $\rho_{\mathrm{tot}}(M):=\sum_{i=1}^s\rho_i(M)$. This additive nonnegative integer is positive exactly when $M$ has one-dimensional support. The next lemma uses bounds on the individual $\rho_i$ to maximize numerical values; Proposition \ref{prop:positive-numerical-pairs} uses $\rho_{\mathrm{tot}}$ to control semistable objects.

\begin{Lem}\label{lem:bounded-numerical-tests}
Fix $M\in\Coh_0Y$ and put $h_j:=\ell_R H^j(Y,M)$ for $j=0,1$. For every subsheaf $N\subseteq M$ and quotient $Q$ of $M$,
\[
\begin{gathered}
0\le\rho_i(N),\rho_i(Q)\le\rho_i(M),\\
\chi_R(N)\le h_0,\qquad\chi_R(Q)\ge-h_1.
\end{gathered}
\]
Consequently, for any $\theta\in V_{>0}$, the values $\theta(N)$ on subsheaves of $M$ have a maximum. Moreover, for any compact $B\subseteq V_{>0}$ there is $a>0$ such that, for every $\eta\in B$,
\[
\chi_R(Q)>a\ \Longrightarrow\ \eta(Q)>0,
\qquad
\chi_R(N)<-a\ \Longrightarrow\ \eta(N)<0.
\]
\end{Lem}
\begin{proof}
The bounds on $\rho_i$ follow from additivity and nonnegativity. For $0\to N\to M\to Q\to0$, the maps
\[
H^0(Y,N)\hookrightarrow H^0(Y,M),
\qquad H^1(Y,M)\twoheadrightarrow H^1(Y,Q)
\]
give the bounds on $\chi_R$; the second map is surjective because $H^2(Y,N)=0$ by Lemma \ref{lem:cohomological-dimension-one}.

Fix $\theta\in V_{>0}$. Only finitely many multiplicity vectors $(\rho_1(N),\ldots,\rho_s(N))$ occur. For each vector that occurs, the corresponding values of $\chi_R(N)$ form a nonempty set of integers bounded above by $h_0$, and hence have a largest value. Since $\rk\theta>0$, formula \eqref{eq:rank-degree-euler} shows that this also maximizes $\theta(N)$ among subsheaves with that vector. Taking the largest of these finitely many maxima proves the claim.

On $B$, the rank has a positive lower bound and the $d_i$ are bounded. The sum $\sum_i d_i(\eta)\rho_i(-)$ is therefore bounded in absolute value on all these subobjects and quotients, uniformly for $\eta\in B$. Formula \eqref{eq:rank-degree-euler} proves the last assertions.
\end{proof}

For $\theta\in V$, define four full subcategories of $\Coh_0Y$:
\begin{equation}\label{eq:numerical-classes}
\begin{aligned}
\T_\theta&:=\{M\mid\theta(Q)>0\text{ for every nonzero quotient }Q\text{ of }M\},\\
\overline{\T}_\theta&:=\{M\mid\theta(Q)\ge0\text{ for every quotient }Q\text{ of }M\},\\
\overline{\F}_\theta&:=\{M\mid\theta(N)\le0\text{ for every subsheaf }N\subseteq M\},\\
\F_\theta&:=\{M\mid\theta(N)<0\text{ for every nonzero subsheaf }N\subseteq M\}.
\end{aligned}
\end{equation}
The category of $\theta$-semistable sheaves is
\[
\W_\theta:=\overline{\T}_\theta\cap\overline{\F}_\theta.
\]
Equivalently, $M\in\W_\theta$ if $\theta(M)=0$ and $\theta(Q)\ge0$ for all quotients $Q$ of $M$. We say that two parameters are \emph{TF equivalent} if their classes $\overline{\T}$ and $\overline{\F}$ coincide. On $V_{>0}$ this is also equivalent to equality of both numerical torsion pairs below. This is the TF equivalence used by Asai in the module setting \cite[Section 2.2]{Asa21}.

\begin{Prop}\label{prop:positive-numerical-pairs}
For $\theta\in V_{>0}$, the pairs
\[
(\T_\theta,\overline{\F}_\theta),
\qquad (\overline{\T}_\theta,\F_\theta)
\]
are torsion pairs in $\Coh_0Y$. The category $\W_\theta$ is an abelian length category. Furthermore, $\W_\theta=0$ if and only if
\[
\T_\theta=\overline{\T}_\theta,
\qquad \F_\theta=\overline{\F}_\theta.
\]
\end{Prop}
\begin{proof}
A nonzero zero-dimensional sheaf $Q$ satisfies
\[
\theta(Q)=(\rk\theta)\ell_R H^0(Y,Q)>0
\]
by Lemma \ref{lem:numerical-devissage}. Fix $M$. By Lemma \ref{lem:bounded-numerical-tests}, the values $\theta(N)$ for $N\subseteq M$ have a maximum $a$.

Choose a subsheaf $N$ with $\theta(N)=a$ which is maximal under inclusion; such a choice exists by noetherianity. Every quotient of $N$ has nonnegative value, since its kernel is a subsheaf of $M$. Every nonzero subsheaf of $M/N$ has strictly negative value, by maximality of $a$ and of $N$. This gives the torsion decomposition for $(\overline{\T}_\theta,\F_\theta)$.

Instead choose, among the subsheaves with value $a$, one minimizing $\rho_{\mathrm{tot}}(N)$. A nonzero quotient of $N$ of value zero would have positive $\rho_{\mathrm{tot}}$, so its kernel would contradict this minimality. Thus $N\in\T_\theta$ and $M/N\in\overline{\F}_\theta$. Orthogonality for both decompositions follows by applying the defining inequalities to the image of a nonzero morphism.

Additivity of $\theta$ shows that $\W_\theta$ is closed under kernels, cokernels and extensions. Every nonzero object of $\W_\theta$ has positive integer $\rho_{\mathrm{tot}}$, which is additive in its short exact sequences. Any strict chain of subobjects of $M$ in $\W_\theta$ therefore has length at most $\rho_{\mathrm{tot}}(M)$, proving that $M$ has finite length. Applying the first torsion decomposition to an object of $\overline{\T}_\theta$ gives a quotient in $\W_\theta$; the analogous statement for the other pair proves the asserted criterion for its vanishing.
\end{proof}

\begin{Cor}\label{cor:numerical-openness}
For a fixed $M\in\Coh_0Y$, the loci
\[
\{\theta\in V_{>0}\mid M\in\T_\theta\},
\qquad
\{\theta\in V_{>0}\mid M\in\F_\theta\}
\]
are open. The corresponding loci for $\overline{\T}$ and $\overline{\F}$ are closed in $V_{>0}$.
\end{Cor}
\begin{proof}
By Lemma \ref{lem:bounded-numerical-tests}, only finitely many multiplicity
vectors occur among subsheaves and quotients of $M$. For each vector
$\alpha$ occurring among nonzero quotients, choose $Q_\alpha$ minimizing
$\chi_R$; for each vector $\beta$ occurring among nonzero subsheaves,
choose $N_\beta$ maximizing $\chi_R$. These extrema exist because the
respective Euler characteristics are integers bounded below and above.

Since $\rk\theta>0$, formula \eqref{eq:rank-degree-euler} shows that
$\theta(-)$ increases with $\chi_R(-)$ at fixed multiplicities. Hence
\[
\begin{aligned}
M\in\T_\theta
&\Longleftrightarrow \theta(Q_\alpha)>0\quad\text{for every }\alpha,\\
M\in\F_\theta
&\Longleftrightarrow \theta(N_\beta)<0\quad\text{for every }\beta.
\end{aligned}
\]
These finite systems of strict linear inequalities define open subsets
of $V_{>0}$. Replacing the strict inequalities by weak ones gives the
corresponding loci for $\overline{\T}$ and $\overline{\F}$, which are
therefore closed in $V_{>0}$.
\end{proof}

\subsection{TF equivalence and closed-fibre walls}

A nonzero sheaf $S$ is a \emph{brick} if $\End_Y(S)$ is a division ring. For $\theta\in V_{>0}$, a \emph{$\theta$-stable sheaf} is a simple object of $\W_\theta$; equivalently, $\theta(S)=0$ and $\theta(Q)>0$ for every proper nonzero quotient $Q$ of $S$. For a nonzero $M\in\Coh_0Y$, put
\[
\Theta_M:=\{\theta\in V_{>0}\mid M\in\W_\theta\},
\quad
\mathrm{Wall}:=\bigcup_{0\ne M\in\Coh_0Y}\Theta_M.
\]
We identify $\Coh C$ with the full subcategory of $\Coh_0Y$ annihilated by $\m$ via the closed immersion $C\hookrightarrow Y$.

\begin{Prop}\label{prop:closed-fibre-walls}
Every brick $S\in\Coh Y$ satisfies $\m S=0$, and hence belongs to $\Coh C$. In particular, every stable sheaf lies on the scheme-theoretic closed fibre, and
\begin{equation}\label{eq:closed-fibre-brick-walls}
\mathrm{Wall}
=\bigcup_{\substack{S\in\Coh C\\S\text{ a brick}}}\Theta_S.
\end{equation}
\end{Prop}
\begin{proof}
For $r\in\m$, multiplication by $r$ on the brick $S$ is either zero or an isomorphism. In the latter case, $\End_Y(S)=r\End_Y(S)$. Properness makes $\End_Y(S)$ a nonzero finite $R$-module, so this contradicts Nakayama's lemma. Thus $\m S=0$.

A stable sheaf is a brick by simplicity in $\W_\theta$. By Proposition \ref{prop:positive-numerical-pairs}, every nonzero semistable sheaf has a finite composition series in $\W_\theta$, whose factors are stable sheaves. At least one such factor exists, and the first assertion places it in $\Coh C$. This proves \eqref{eq:closed-fibre-brick-walls}.
\end{proof}

The following segment criterion extends \cite[Theorem 2.16]{Asa21}. Positive rank permits the use of generic multiplicity in place of composition length in the ambient category.

\begin{Thm}\label{thm:tf-brick-walls}
Let $\theta,\eta\in V_{>0}$ be distinct, and put $\theta_t:=(1-t)\theta+t\eta$ for $0\le t\le1$. The following are equivalent:
\begin{enumerate}
\item $\theta$ and $\eta$ are TF equivalent;
\item the category $\W_{\theta_t}$ is independent of $t\in[0,1]$;
\item for every brick $S\in\Coh C$, the intersection $[\theta,\eta]\cap\Theta_S$ is not a singleton.
\end{enumerate}
Condition (3) is unchanged if bricks in all of $\Coh Y$ are used.
\end{Thm}
\begin{proof}
Suppose (1) holds. The defining inequalities give
\[
\overline{\T}_\theta\subseteq\overline{\T}_{\theta_t},
\qquad
\F_\theta\subseteq\F_{\theta_t}.
\]
Both pairs are torsion pairs by Proposition \ref{prop:positive-numerical-pairs}, so these inclusions force equality. Applying the same argument to $(\T,\overline{\F})$ shows that every $\theta_t$ is TF equivalent to $\theta$. Thus (2) holds. Since $\theta\ne\eta$, condition (2) implies (3).

For $(3)\Rightarrow(1)$, suppose that (1) fails. Torsion decompositions give, after exchanging $\theta$ and $\eta$ if necessary, a nonzero intersection among
\begin{equation}\label{eq:tf-crossing-intersections}
\overline{\T}_\theta\cap\F_\eta,
\qquad
\T_\theta\cap\overline{\F}_\eta.
\end{equation}
Since $\rk\eta>0$, nonzero sheaves in either intersection have no zero-dimensional subsheaf. Thus $\rho_{\mathrm{tot}}$ is positive on them and strictly decreases on proper quotients. Choose a nonzero $M$ of minimal $\rho_{\mathrm{tot}}$ in one such intersection. The rest is exactly Asai's argument \cite[proof of Theorem 2.16]{Asa21}, with $\rho_{\mathrm{tot}}$ replacing module dimension; the endomorphism rings here have finite $R$-length. It shows that $M$ is a brick and $[\theta,\eta]\cap\Theta_M$ is a singleton. Proposition \ref{prop:closed-fibre-walls} places $M$ on $C$, contradicting (3), and also gives the final assertion.
\end{proof}

\subsection{Pretilting cones and semistable reduction}

The following is a geometric analogue of Yurikusa's semistability theorem for two-term presilting complexes over finite-dimensional algebras \cite[Theorem 1.4 and Proposition 3.3]{Yur18}. It identifies the numerical torsion classes and the semistable category on each pretilting cone.

\begin{Prop}\label{prop:pretilting-numerical-classes}
Let $U$ be a basic pretilting bundle and $\theta\in C^\circ(U)$. Then
\begin{equation}\label{eq:pretilting-numerical-classes}
\T_\theta=\Fac U\cap\Coh_0Y,
\qquad \overline{\T}_\theta=U^{\perp_1}\cap\Coh_0Y,
\qquad \overline{\F}_\theta=U^{\perp_0}\cap\Coh_0Y.
\end{equation}
In particular,
\begin{equation}\label{eq:semistable-orthogonal}
\W_\theta=\{M\in\Coh_0Y\mid\Hom_Y(U,M)=0=\Ext_Y^1(U,M)\}.
\end{equation}
\end{Prop}
\begin{proof}
The case $U=0$ is immediate. Write $U=\bigoplus_iU_i$ and $\theta=\sum_i a_i[U_i]$ with $a_i>0$. As in \cite[Proposition 3.3]{Yur18}, use the Euler pairing, with $R$-lengths in place of vector-space dimensions:
\[
\theta(M)=\sum_i a_i\bigl(
\ell_R\Hom_Y(U_i,M)-\ell_R\Ext_Y^1(U_i,M)\bigr)
\qquad(M\in\Coh_0Y).
\]
On $\Coh_0Y$, $\theta$ is nonnegative on $U^{\perp_1}$, nonpositive on $U^{\perp_0}$, and strictly positive on nonzero sheaves in $\Fac U$; the last assertion uses pretilting and cohomological dimension one.

We now use torsion decompositions as in \cite[proof of Theorem 1.4]{Yur18}. Propositions \ref{Bongcomp} and \ref{prop:bundle-to-torsion} give torsion pairs
\[
(\Fac U,U^{\perp_0}),
\qquad (U^{\perp_1},\E_{\max}^{\perp_0})
\]
in $\Coh Y$, which restrict to $\Coh_0Y$. The signs above and the first decomposition give $\T_\theta=\Fac U\cap\Coh_0Y$ and $\overline{\F}_\theta=U^{\perp_0}\cap\Coh_0Y$.

For the second decomposition, every nonzero $F\in\E_{\max}^{\perp_0}\cap\Coh_0Y$ has $\Hom_Y(U,F)=0$ because $U\in\add\E_{\max}$. Also $\Ext_Y^1(U,F)\ne0$, since otherwise $F$ belongs to both classes of that torsion pair. Hence $\theta(F)<0$. Within $\Coh_0Y$, all quotients of a sheaf in $U^{\perp_1}$ have nonnegative value, while any sheaf outside it has a nonzero quotient of negative value. This proves $\overline{\T}_\theta=U^{\perp_1}\cap\Coh_0Y$. Intersecting with $\overline{\F}_\theta$ gives \eqref{eq:semistable-orthogonal}.
\end{proof}

The extremal completions identify the orthogonal category in $\Coh Y$ and its semistable part on the closed fibre. This extends the idempotent-quotient description of \cite[Proposition 4.1]{Asa21} to our complete local base.

\begin{Thm}\label{thm:semistable-reduction}
Let $U\ne0$ be basic pretilting with $r:=|U|$, and let $\E_{\min}$ and $\E_{\max}$ be its extremal completions from Proposition \ref{Bongcomp}. For $\sigma\in\{\min,\max\}$, put
\[
\Lambda_\sigma:=\End_Y(\E_\sigma),
\qquad
B_\sigma:=\Lambda_\sigma/\Lambda_\sigma e_U\Lambda_\sigma,
\]
where $e_U$ is the idempotent corresponding to the summand $U$. Put
\[
\widetilde{\W}_U:=U^{\perp_0}\cap U^{\perp_1}\subseteq\Coh Y,
\qquad
\W_U:=\widetilde{\W}_U\cap\Coh_0Y.
\]
There are exact $R$-linear equivalences
\begin{equation}\label{eq:semistable-reduction-equivalences}
\begin{aligned}
\Hom_Y(\E_{\max},-)&\colon\widetilde{\W}_U\xrightarrow[\simeq]{}\mod B_{\max},\\
\Ext_Y^1(\E_{\min},-)&\colon\widetilde{\W}_U\xrightarrow[\simeq]{}\mod B_{\min}.
\end{aligned}
\end{equation}
They restrict to $\W_U\simeq\fl B_\sigma$ for each $\sigma$. For $\theta\in C^\circ(U)$, one has $\W_U=\W_\theta$; this is a length category with exactly $n-r$ isomorphism classes of simple objects. If $\E_{\max}=\bigoplus_{i=1}^n\E_i$ is indexed so that $U=\bigoplus_{i=1}^r\E_i$, these simples are $S_{r+1},\ldots,S_n$ from Lemma \ref{lem:euler-dual-bases}; in particular, they are sheaves. Thus
\[
\dim C(U)+\operatorname{rank}_{\mathbb Z}K_0(\W_U)=n.
\]
\end{Thm}
\begin{proof}
For either completion, consider
\[
\B_\sigma:=\{Z\in\H_{\E_\sigma}\mid
\Hom_{D^b(\Coh Y)}(U,Z)=0\}.
\]
Since $U$ is a projective object of $\H_{\E_\sigma}$, Morita equivalence identifies $\B_\sigma$ with the Serre subcategory $\mod B_\sigma\subseteq\mod\Lambda_\sigma$.

Take $Z\in\B_{\max}$ and its short exact sequence in the HRS heart
\[
0\longrightarrow F[1]\longrightarrow Z\longrightarrow T\longrightarrow0.
\]
Since $\B_{\max}$ is Serre, $F[1]\in\B_{\max}$, so $\Ext_Y^1(U,F)=0$. Hence $F\in U^{\perp_1}=\Fac\E_{\max}$. But $F\in\E_{\max}^{\perp_0}$, so $F=0$. Thus $Z=T$ is a sheaf in $\widetilde{\W}_U$. Conversely, every $M\in\widetilde{\W}_U$ lies in $\Fac\E_{\max}$ and is annihilated by $\Hom_Y(U,-)$, giving $\B_{\max}=\widetilde{\W}_U$.

For $Z\in\B_{\min}$, the same sequence has $T\in\B_{\min}$. Now $T\in\Fac\E_{\min}=\Fac U$ and $\Hom_Y(U,T)=0$, forcing $T=0$. Since $\E_{\min}^{\perp_0}=U^{\perp_0}$, we get $Z=F[1]$ with $F\in\widetilde{\W}_U$. The converse is immediate, so $\B_{\min}=\widetilde{\W}_U[1]$. Morita equivalence gives \eqref{eq:semistable-reduction-equivalences}.

These $R$-linear equivalences restrict to objects annihilated by a power of $\m$, giving $\W_U\simeq\fl B_\sigma$. Equation \eqref{eq:semistable-orthogonal} identifies $\W_U$ with $\W_\theta$. Its simples are exactly those not indexed by summands of $U$, so Lemma \ref{lem:euler-dual-bases} gives their description and the rank formula.
\end{proof}

The same dual simples recover the closed cone intrinsically from its weak numerical classes, as in the closure criterion of \cite[Lemma 2.15]{Asa21}.

\begin{Thm}\label{thm:pretilting-cone-inequalities}
For every basic pretilting bundle $U$,
\begin{equation}\label{eq:pretilting-cone-inequalities}
C(U)=\left\{\eta\in V\ \middle|\
\begin{aligned}
\eta(M)&\ge0&&\text{for all }M\in U^{\perp_1}\cap\Coh_0Y,\\
\eta(N)&\le0&&\text{for all }N\in U^{\perp_0}\cap\Coh_0Y
\end{aligned}\right\}.
\end{equation}
Equivalently, for $\theta\in C^\circ(U)$,
\begin{equation}\label{eq:pretilting-cone-tf-closure}
C(U)=\{\eta\in V\mid
\overline{\T}_\theta\subseteq\overline{\T}_\eta,\
\overline{\F}_\theta\subseteq\overline{\F}_\eta\}.
\end{equation}
\end{Thm}
\begin{proof}
If $\eta=\sum_i a_i[U_i]$ with $a_i\ge0$, the Euler pairing gives the asserted signs. Conversely, suppose $\eta$ satisfies the inequalities. For $U=0$, they make $\eta$ vanish on all of $\Coh_0Y$, so $\eta=0$ by Lemma \ref{lem:euler-dual-bases}.

For $U\ne0$, use the indexing of $\E_{\max}$ and its dual simples in Theorem \ref{thm:semistable-reduction}. Since $\Fac\E_{\max}=U^{\perp_1}$ and $\E_{\max}^{\perp_0}\subseteq U^{\perp_0}$, Corollary \ref{cor:cone-from-heart} gives
\[
\eta=\sum_{i=1}^n a_i[\E_i]\in C(\E_{\max}),\qquad a_i\ge0.
\]
For $j>r$, we have $S_j\in\W_U$, so both inequalities apply and $c_ja_j=\eta(S_j)=0$. Thus $\eta=\sum_{i=1}^r a_i[U_i]\in C(U)$.
Finally, Proposition \ref{prop:pretilting-numerical-classes} identifies the two test categories with $\overline{\T}_\theta$ and $\overline{\F}_\theta$. They are quotient closed and subobject closed, respectively, so their sign conditions are equivalent to the inclusions in \eqref{eq:pretilting-cone-tf-closure}.
\end{proof}

\begin{Cor}\label{cor:semistable-reduction-map}
For $U\ne0$ basic pretilting, restriction of the Euler pairing gives a surjective linear map
\[
\pi_U\colon V\longrightarrow
\Hom_{\mathbb R}(K_0(\W_U)\otimes_{\mathbb Z}\mathbb R,\mathbb R),
\qquad \eta\longmapsto\eta|_{\W_U},
\]
with $\ker\pi_U=\operatorname{span}_{\mathbb R}C(U)$.
\end{Cor}
\begin{proof}
The classes of $S_{r+1},\ldots,S_n$ form a basis of $K_0(\W_U)$. Their pairings with the completion classes are given by \eqref{eq:dual-bases}. Thus $[\E_{r+1}],\ldots,[\E_n]$ map to a basis of the target, and the kernel is spanned by $[U_1],\ldots,[U_r]$.
\end{proof}

We now combine the closed-cone description with the simple-object count in the semistable category to identify the entire TF class of each relative interior. This gives the geometric counterpart of \cite[Theorem 3.11]{Asa21}.

\begin{Thm}\label{thm:tf-cones}
For every basic pretilting bundle $U$, the set $C^\circ(U)$ is an entire TF equivalence class in $V$.
\end{Thm}
\begin{proof}
Proposition \ref{prop:pretilting-numerical-classes} proves constancy on $C^\circ(U)$. Fix $\theta\in C^\circ(U)$ and let $\eta$ be TF equivalent to it. Theorem \ref{thm:pretilting-cone-inequalities} gives $\eta\in C(U)$, which already proves the assertion for $U=0$.

Suppose $U\ne0$. A closed-point sheaf belongs to $\W_0$ but not to $\W_\theta$, since $\theta$ has positive rank. Equality of semistable categories therefore excludes $\eta=0$. There is a nonzero direct summand $U'$ of $U$ with $\eta\in C^\circ(U')$. By Theorem \ref{thm:semistable-reduction}, the common category $\W_\theta=\W_\eta$ has both $n-|U|$ and $n-|U'|$ simple objects. Thus $|U'|=|U|$, so $U'=U$ and $\eta\in C^\circ(U)$.
\end{proof}

The following intersection formula is a geometric counterpart of \cite[Proposition 3.9(3)]{Asa21}; we prove it using minimal tilting completions and Theorem \ref{thm:tilting-torsion}.

\begin{Thm}\label{thm:tilting-fan}
For $U,U'\in\ptilt Y$, let $W$ be the basic bundle with $\add W=\add U\cap\add U'$. Then
\begin{equation}\label{eq:pretilting-cone-intersection}
C(U)\cap C(U')=C(W).
\end{equation}
Consequently, $\Sigma(Y):=\{C(U)\mid U\in\ptilt Y\}$ is a rational simplicial fan, and each cone is a face of a full dimensional tilting cone.
\end{Thm}
\begin{proof}
Let $0\ne\theta\in C(U)\cap C(U')$, and choose summands $A$ of $U$ and $B$ of $U'$ with $\theta\in C^\circ(A)\cap C^\circ(B)$. Propositions \ref{Bongcomp} and \ref{prop:pretilting-numerical-classes} give basic minimal tilting completions $\E_A,\E_B$ with
\[
\Fac\E_A\cap\Coh_0Y=\T_\theta=\Fac\E_B\cap\Coh_0Y.
\]
Theorem \ref{thm:tilting-torsion} therefore gives $\E_A\cong\E_B=:\E$. In the summand basis of $\E$ from Lemma \ref{lem:euler-dual-bases}, the summands of $A$ and $B$ are both exactly those with positive coefficients in $\theta$, so $A\cong B$. Hence $\theta\in C(W)$; the reverse inclusion and the zero case are immediate.

The intersection formula, Proposition \ref{Bongcomp} and Lemma \ref{lem:euler-dual-bases} prove the fan assertions, since summands give faces and every pretilting bundle has a tilting completion.
\end{proof}

\begin{Rem}
Under the Euler pairing, Corollary \ref{cor:cone-from-heart} identifies $C(\E)$ with the heart cone of the HRS heart $\H_{\E,0}$ for every $\E\in\tilt Y$. Since every pretilting cone is a face of a tilting cone, $\Sigma(Y)$ is a subfan of the heart fan of $\Coh_0Y$ constructed by Broomhead--Pauksztello--Ploog--Woolf \cite[Theorem-Definition 3.7]{BPPW24}.
\end{Rem}

\subsection{Chambers}

The following chamber correspondence parallels \cite[Theorem 3.17]{Asa21}. Its proof uses the length-heart characterization from Section \ref{sec:tilting-torsion}.

The wall set defined above is $\mathrm{Wall}=\{\theta\in V_{>0}\mid\W_\theta\ne0\}$. Its closure is taken in $V_{>0}$. A chamber is a connected component of $V_{>0}\setminus\overline{\mathrm{Wall}}$.

\begin{Thm}\label{thm:tilting-chambers}
The assignment $\E\mapsto C^\circ(\E)$ gives a bijection from basic tilting bundles to chambers. More precisely,
\[
\coprod_{\E\in\tilt Y}C^\circ(\E)
=V_{>0}\setminus\overline{\mathrm{Wall}}.
\]
For rational parameters, one also has
\[
\theta\in V_{>0}\cap(K_0(\per Y)\otimes\mathbb Q),\quad\W_\theta=0
\quad\Longleftrightarrow\quad
\theta\in C^\circ(\E)\text{ for some }\E\in\tilt Y.
\]
\end{Thm}
\begin{proof}
By Proposition \ref{prop:pretilting-numerical-classes}, every open tilting cone is disjoint from $\mathrm{Wall}$; being open, it also avoids $\overline{\mathrm{Wall}}$.

For rational $\theta$ of positive rank with $\W_\theta=0$, we follow \cite[proof of Proposition 3.19]{Asa21}. Proposition \ref{prop:positive-numerical-pairs} gives the torsion pair $(\T_\theta,\F_\theta)$. After multiplication by a positive integer, $\theta$ takes positive integer values on every nonzero object of the HRS heart $\F_\theta[1]*\T_\theta$, so this heart is length. Since $\T_\theta$ contains every closed-point sheaf, Theorem \ref{thm:tilting-torsion} reconstructs a basic tilting bundle $\E$ inducing this pair, and Corollary \ref{cor:cone-from-heart} gives $\theta\in C^\circ(\E)$.

For real $\theta\in V_{>0}\setminus\overline{\mathrm{Wall}}$, the argument of \cite[proof of Theorem 3.17]{Asa21} applies: choose a rational $\eta$ in a convex neighbourhood of $\theta$ contained in this complement. Theorem \ref{thm:tf-brick-walls} makes the two parameters TF equivalent, so the rational case and Theorem \ref{thm:tf-cones} put $\theta$ in the same open tilting cone as $\eta$. This proves the displayed equality.

By Theorem \ref{thm:tilting-fan}, the open cones in this decomposition are pairwise disjoint and connected, hence are precisely the chambers.
\end{proof}

\subsection{Classes, restriction and wall crossing}

We recover pretilting bundles from their classes, identify their fans by restriction to the full closed fibre, and construct exchange sequences across facets.

\begin{Cor}\label{cor:class-determines-pretilting}
A pretilting bundle is determined up to isomorphism by its rank and determinant, equivalently by its class in $K_0(\per Y)$. For basic pretilting bundles, equality of their cones also implies isomorphism.
\end{Cor}
\begin{proof}
By Corollary \ref{cor:rank-determinant}, it suffices to take $P,Q$ with the same class $\theta$. If $\theta=0$, both bundles vanish by rank. Otherwise, their basic representatives $A,B$ satisfy $\theta\in C^\circ(A)\cap C^\circ(B)$. By \eqref{eq:pretilting-cone-intersection}, the intersection is the cone of their common summands. As a face containing a relative interior point of each cone, it corresponds to all summands of both $A$ and $B$. Thus $A\cong B$. In a tilting completion of $A$, Lemma \ref{lem:euler-dual-bases} then identifies the multiplicities in $P$ and $Q$, proving $P\cong Q$.

For basic bundles, the same intersection formula shows that equality of cones forces all their summands to be common, proving the second assertion.
\end{proof}

\begin{Cor}[Numerical triviality]\label{cor:numerically-trivial-pretilting}
Let $F$ be a pretilting bundle of rank $r$. If $d_i(F)=0$ for every one-dimensional component of the closed fibre, then
\[
F\cong\O_Y^{\oplus r}.
\]
\end{Cor}
\begin{proof}
Lemma \ref{lem:numerical-devissage} and the torsion-freeness in Corollary \ref{cor:grothendieck-groups} give $[F]=r[\O_Y]$ in $K_0(\per Y)$. Since $H^1(Y,\O_Y)=0$, the bundle $\O_Y^{\oplus r}$ is pretilting. Corollary \ref{cor:class-determines-pretilting} identifies the two bundles.
\end{proof}

\begin{Prop}[Closed-fibre invariance]\label{prop:closed-fibre-fan}
The closed fibre $C\to\Spec k$ satisfies \eqref{eq:geometric-setup}. Restriction gives a bijection on isomorphism classes of pretilting bundles, preserving direct-sum decompositions and compatibility. In particular, it induces bijections
\[
\ptilt Y\xrightarrow{\ \sim\ }\ptilt C,
\qquad
\tilt Y\xrightarrow{\ \sim\ }\tilt C,
\]
and an isomorphism
\[
K_0(\per Y)\xrightarrow{\ \sim\ }K_0(\per C)
\]
identifying their pretilting fans, ranks and degrees.
\end{Prop}
\begin{proof}
Put $I:=\m\O_Y$. As $I$ is a quotient of $\O_Y^{\oplus N}$, right exactness of $H^1$ gives $H^1(Y,I)=H^1(C,\O_C)=0$. Moreover $\Gamma(Y,I)=\m$: it contains $\m$ and cannot contain a unit on the nonempty closed fibre. Thus $H^0(C,\O_C)=k$.

Let $i\colon C\hookrightarrow Y$ be the closed immersion and write $i^*\colon K_0(\per Y)\to K_0(\per C)$. By Corollary \ref{cor:rank-determinant}, $i^*$ is the direct sum of the identity on $\mathbb Z$ and Picard restriction. Lemma \ref{lem:picard-restriction} therefore makes it an isomorphism. It preserves ranks and degrees.

For a sheaf on $C$, its subobjects, quotients and Euler pairing agree in the two settings. Proposition \ref{prop:closed-fibre-walls} thus identifies the wall sets under $i^*$. Theorem \ref{thm:tilting-chambers}, whose reconstruction step uses Theorem \ref{thm:tilting-torsion}, identifies the chambers with open tilting cones on both $Y$ and $C$. Taking closures therefore identifies their maximal cones. Every pretilting cone is a face of a tilting cone, and every such face is pretilting, by Theorem \ref{thm:tilting-fan}. Hence
\[
i^*(\Sigma(Y))=\Sigma(C).
\]

It remains to identify this correspondence with restriction of bundles. Corollary \ref{lem:closed-fibre-detects-ext} gives, for bundles $F,G$ on $Y$,
\[
\Ext_Y^1(F,G)=0
\quad\Longleftrightarrow\quad
\Ext_C^1(F|_C,G|_C)=0.
\]
In particular, restriction preserves pretilting bundles. Let $\E=\bigoplus_j\E_j$ be basic tilting on $Y$, and let $T=\bigoplus_jT_j$ be basic tilting on $C$ with $C(T)=i^*(C(\E))$. The summand classes of each tilting bundle form a lattice basis, so they are the primitive generators of the rays of its cone. After reindexing, $[\E_j|_C]=[T_j]$, and Corollary \ref{cor:class-determines-pretilting} gives $\E_j|_C\cong T_j$. Thus restriction realizes the bijection on basic tilting bundles and their summands. Completing basic pretilting bundles by Proposition \ref{Bongcomp} gives the bijection on all faces, and arbitrary multiplicities give the assertion for all pretilting bundles. The displayed Ext criterion preserves compatibility.
\end{proof}

\begin{Rem}
Proposition \ref{prop:closed-fibre-fan} is a geometric analogue of Kimura's reduction theorem for two-term silting complexes under $\Lambda\mapsto\Lambda\otimes_Rk=\Lambda/\m\Lambda$, for a module-finite algebra over a complete local noetherian ring \cite[Theorem 4.3]{Kim24}. Compare also Fushimi's invariance of silting objects under derived base change to the residue field \cite[Theorem 4.9]{Fus26}.
\end{Rem}

We now describe the two tilting completions of an almost complete pretilting bundle and the exchange sequence between them.

\begin{Thm}\label{thm:facet-mutation}
Let $U\ne0$ be basic pretilting with $|U|=n-1$, and let $\E_-:=\E_{\min}$ and $\E_+:=\E_{\max}$ be its extremal completions from Proposition \ref{Bongcomp}.
\begin{enumerate}
\item These are the only two basic tilting completions of $U$. Their cones meet along $C(U)$ and lie on opposite sides of its span.
\item Write $\E_-=U\oplus V_-$ and $\E_+=U\oplus V_+$ with $V_-,V_+$ indecomposable. There is a short exact sequence
\[
0\longrightarrow V_+\longrightarrow U_0\longrightarrow V_-\longrightarrow0,
\qquad U_0\in\add U.
\]
The first and second maps are respectively minimal left and right $(\add U)$-approximations. In particular, $\E_-$ and $\E_+$ are related by an irreducible mutation.
\end{enumerate}
\end{Thm}
\begin{proof}
The category $\per Y$ is Krull--Schmidt, and its finite $R$-module Hom spaces ensure the approximation hypotheses of \cite[Section 2.4]{AI}.

Write $\E_-=U\oplus V_-$ and take a minimal right $(\add U)$-approximation $b\colon U_0\to V_-$. Proposition \ref{Bongcomp} gives $V_-\in\Fac\E_-=\Fac U$, so $b$ is surjective and has a locally free kernel:
\[
0\longrightarrow K\xrightarrow{a}U_0\xrightarrow{b}V_-\longrightarrow0.
\]
Its triangle defines the right silting mutation of $\E_-$ at $V_-$. By the mutation and inverse-mutation results of Aihara--Iyama \cite[Theorem 2.31, Proposition 2.33 and Theorem 2.35]{AI}, $K$ is indecomposable, $U\oplus K$ is basic silting and differs from $\E_-$, and $a$ is a minimal left $(\add U)$-approximation. Since $K$ is a vector bundle, $U\oplus K$ is a tilting bundle.

Every tilting completion has facet $C(U)$. The fan property in Theorem \ref{thm:tilting-fan} allows at most one completion on each side of its span, so $\E_-$ and $U\oplus K$ are the only completions and their cones lie on opposite sides. The extremal property in Proposition \ref{Bongcomp} identifies $U\oplus K$ with $\E_+$, hence $K\cong V_+$. The displayed sequence gives (2), and the fan intersection formula gives (1).
\end{proof}

For $\E=U\oplus V\in\tilt Y$ with $U\ne0$ and $V$ indecomposable, write $\mu_V(\E)$ for the unique other tilting completion of $U$.

For $\theta\in C^\circ(U)$, Theorem \ref{thm:semistable-reduction} with $r=n-1$ shows that $\W_\theta$ has exactly one simple object and gives its idempotent-quotient description for either completion.

\begin{Rem}
The hypothesis $U\ne0$ in Theorem \ref{thm:facet-mutation} is necessary. If $Y=\Spec R$, then $n=1$, the zero bundle has $n-1$ summands, and its only basic tilting completion is $\O_Y$. The assertions about two adjacent completions concern nonzero facets in the positive-rank region.
\end{Rem}

\section{Completeness, finiteness and mutation connectivity}\label{sec:completeness}

Retain the general assumptions \eqref{eq:geometric-setup}. We characterize completeness of the tilting fan by finiteness conditions and then deduce mutation connectivity. We call the family $\{\Delta(\E)\mid\E\in\tilt Y\}$ locally finite if every compact subset of $H$ meets only finitely many of its members. Completeness in $V_{>0}$ means that their union is $H$.

\begin{Thm}\label{thm:complete-finite}
Under \eqref{eq:geometric-setup}, the following are equivalent:
\begin{enumerate}
\item $H=\bigcup_{\E\in\tilt Y}\Delta(\E)$;
\item the family $\{\Delta(\E)\mid\E\in\tilt Y\}$ is locally finite;
\item there are finitely many basic tilting bundles modulo tensoring by line bundles;
\item the ranks of indecomposable pretilting bundles are uniformly bounded.
\end{enumerate}
\end{Thm}

The proof follows Asai's use of open sets associated with presilting objects \cite[Lemma 4.3 and proof of Proposition 4.9]{Asa21}. Here closed-fibre duality gives the corresponding open sets directly on $H$. Completeness provides a finite cover of a compact Picard fundamental set, which bounds ranks; lattice points then give finiteness modulo line bundles.

\subsection{Degrees and the Picard lattice}

\begin{Lem}\label{lem:picard-lattice}
The map
\[
\mu\colon H\longrightarrow\mathbb R^s,\qquad
\mu(\theta):=(d_1(\theta),\ldots,d_s(\theta))
\]
is an affine isomorphism, taking $[\O_Y]$ to zero. The image
\[
\La_{\mathrm{Pic}}:=\{(d_1(L),\ldots,d_s(L))\mid L\in\Pic(Y)\}\subseteq\mathbb R^s
\]
is a full lattice. Tensoring by $L$ acts on $\mathbb R^s$ by translation by its degree vector. In particular, $H/\La_{\mathrm{Pic}}$ is compact.
\end{Lem}
\begin{proof}
Corollary \ref{cor:rank-determinant} identifies $H$, with origin $[\O_Y]$, with $\Pic(Y)\otimes_{\mathbb Z}\mathbb R$. Lemma \ref{lem:numerical-devissage} makes $\mu$ injective. Moreover,
\[
d_i(\theta\otimes L)=d_i(\theta)+(\rk\theta)d_i(L),
\]
so the action on $H$ is translation.

For each $Z_i$, choose a regular closed point $p_i$ outside all other components of $C_{\mathrm{red}}$; such points exist in the dense regular locus of the integral curve $Z_i$. This point is a Cartier divisor on $C_{\mathrm{red}}$, and the degree vector of $\O_{C_{\mathrm{red}}}(p_i)$ is $[\kappa(p_i):k]e_i$. By Lemma \ref{lem:picard-restriction}, this line bundle lifts to $Y$. These vectors span $\mathbb R^s$, proving surjectivity of $\mu$. They also lie in $\La_{\mathrm{Pic}}\subseteq\mathbb Z^s$, so this subgroup is a full lattice. Its index in $\mathbb Z^s$ need not be one.
\end{proof}

\subsection{Open sets detecting direct summands}

For a bundle $V$ on $Y$, write $V_C:=V|_C$. Let $\omega_C^\bullet$ be the relative dualizing complex of $C/k$ and put $\omega_C:=H^{-1}(\omega_C^\bullet)$. Duality for proper schemes of dimension at most one gives
\begin{equation}\label{eq:closed-fibre-duality}
\Hom_C(M,V_C\otimes\omega_C)
\cong\Hom_k\bigl(\Ext_C^1(V_C,M),k\bigr)
\qquad(M\in\Coh C);
\end{equation}
see \cite[Tag 0BS2]{stacks-project}.

The numerical torsion pairs restrict to $\Coh C$, since this subcategory is closed under subobjects and quotients. Proposition \ref{prop:pretilting-numerical-classes} gives, for $\theta\in C^\circ(W)\cap H$,
\begin{equation}\label{eq:closed-fibre-numerical-classes}
\T_\theta\cap\Coh C=\Fac W_C,
\qquad
\overline{\T}_\theta\cap\Coh C=W_C^{\perp_1},
\end{equation}
where the right-hand categories are taken in $\Coh C$.

\begin{Lem}\label{lem:open-summand-detection}
For an indecomposable pretilting bundle $U$, the set
\[
\mathcal N_U:=\{\theta\in H\mid
U_C\in\T_\theta,\quad U_C\otimes\omega_C\in\F_\theta\}
\]
is open. For every basic pretilting bundle $W$ and $\theta\in C^\circ(W)\cap H$,
\begin{equation}\label{eq:open-summand-detection}
\theta\in\mathcal N_U
\quad\Longleftrightarrow\quad U\in\add W.
\end{equation}
\end{Lem}
\begin{proof}
Openness follows from Corollary \ref{cor:numerical-openness}. Fix
$\theta\in C^\circ(W)\cap H$. If $U\in\add W$, both defining conditions
of $\mathcal N_U$ follow from \eqref{eq:closed-fibre-numerical-classes},
\eqref{eq:closed-fibre-duality} and torsion-pair orthogonality.

Conversely, suppose $\theta\in\mathcal N_U$. Proposition \ref{Bongcomp}
gives a tilting bundle $\E$ with $\Fac\E=\Fac W$. Since
$U_C\in\Fac W_C$, Lemma \ref{lem:progenerator-restriction} applied to
$\E$ gives $U\in\Fac W$. Thus a right $(\add W)$-approximation yields
an exact sequence of vector bundles
\[
0\longrightarrow K\longrightarrow W_0\longrightarrow U\longrightarrow0,
\qquad W_0\in\add W.
\]
The approximation property and pretilting give $\Ext_Y^1(W,K)=0$.
Hence $K_C\in\overline{\T}_\theta$ by Corollary
\ref{lem:closed-fibre-detects-ext} and
\eqref{eq:closed-fibre-numerical-classes}. Orthogonality with
$U_C\otimes\omega_C\in\F_\theta$ and
\eqref{eq:closed-fibre-duality} give $\Ext_C^1(U_C,K_C)=0$.
Corollary \ref{lem:closed-fibre-detects-ext} now gives
$\Ext_Y^1(U,K)=0$, so the sequence splits and $U\in\add W$.
\end{proof}

\begin{proof}[Proof of Theorem \ref{thm:complete-finite}]
Fix a compact fundamental set $F\subseteq H$ for $\La_{\mathrm{Pic}}$ by Lemma \ref{lem:picard-lattice}.

$(1)\Rightarrow(4)$: Completeness places every point of $H$ in $C^\circ(W)\cap H$ for some nonzero basic pretilting $W$. Lemma \ref{lem:open-summand-detection}, applied to a summand of $W$, shows that the $\mathcal N_U$ cover $H$. Choose a finite subcover $\mathcal N_{U_1},\ldots,\mathcal N_{U_t}$ of $F$. For any indecomposable pretilting bundle $V$, choose $L\in\Pic(Y)$ with $[V\otimes L]/\rk V\in F$. Formula \eqref{eq:open-summand-detection} gives $U_j\in\add(V\otimes L)$ for some $j$, hence $V\otimes L\cong U_j$. Thus $\rk V\le\max_j\rk U_j$.

$(4)\Rightarrow(3)$: Let $b\in\mathbb Z_{>0}$ bound these ranks. A basic tilting bundle has $n$ summands, so its rank is at most $nb$. After a line-bundle twist, its class lies in
\[
K_0(\per Y)\cap\bigcup_{r=1}^{nb}rF.
\]
This is a finite set, since $K_0(\per Y)$ is a lattice in $V$ and the union is compact. Corollary \ref{cor:class-determines-pretilting} now gives (3).

$(3)\Rightarrow(2)$: All tilting simplices are Picard-lattice translates of finitely many compact simplices. Such a family is locally finite.

$(2)\Rightarrow(1)$: Put $S:=\bigcup_\E\Delta(\E)$, which is closed by local finiteness. If $n=1$, then $H$ is a point and the assertion is immediate. Assume $n\ge2$, and let $Z$ be the union of faces of codimension at least two in the tilting simplices. Locally, $Z$ is a finite union of such faces, so it is closed and has empty interior.

We show that $H\setminus Z$ is polygonally connected. Given $x,y\in H\setminus Z$, choose a bounded convex open neighbourhood $B\subseteq H$ containing them. Only finitely many of the faces defining $Z$ meet $\overline B$. For each such face $F$, the affine spans of $F\cup\{x\}$ and $F\cup\{y\}$ are proper affine subspaces of $H$, since $F$ has codimension at least two. Choose $z\in B$ outside their union. Then both $[x,z]$ and $[z,y]$ avoid these faces and stay in $B$, so they give a polygonal path in $H\setminus Z$.

Every point of $S\setminus Z$ lies in a simplex interior or in the relative interior of a facet. By Theorems \ref{thm:tilting-fan} and \ref{thm:facet-mutation}, each facet has exactly two adjacent simplices, on opposite sides, whose union contains a neighbourhood of every point of its relative interior. Thus $S\setminus Z$ is open in $H\setminus Z$. It is also nonempty and closed there, hence equals $H\setminus Z$ by connectedness. Since this complement is dense and $S$ is closed, $S=H$.
\end{proof}

\begin{Cor}[Mutation connectivity]\label{cor:complete-mutation-connectivity}
If the equivalent conditions of Theorem \ref{thm:complete-finite} hold, then any two basic tilting bundles on $Y$ are related by a finite sequence of left or right irreducible mutations.
\end{Cor}
\begin{proof}
If $n=1$, every tilting simplex is the point $H$, so Corollary \ref{cor:class-determines-pretilting} gives a unique basic tilting bundle. Assume $n\ge2$, and use $Z$ as in the proof of Theorem \ref{thm:complete-finite}. Given $\E,\E'\in\tilt Y$, choose $x,y$ in the interiors of their simplices and a bounded convex open neighbourhood $B$ containing them.

In the two-segment construction of that proof, choose the intermediate point $z$ also outside the facet hyperplanes of the finitely many tilting simplices meeting $\overline B$. This only adds finitely many proper affine subspaces to avoid. The path $\gamma=[x,z]\cup[z,y]$ stays in $H\setminus Z$, and neither segment lies in a facet hyperplane. Each segment meets each such hyperplane at most once, so the path crosses finitely many facets, transversely and in their relative interiors, one at a time.

By completeness, the successive simplices along this path give a finite gallery from $\Delta(\E)$ to $\Delta(\E')$. Theorem \ref{thm:tilting-fan} identifies each common facet with a nonzero pretilting bundle having $n-1$ summands. Theorem \ref{thm:facet-mutation} therefore identifies every step with an irreducible mutation.
\end{proof}

For complex single-curve cDV contractions, Shimpi \cite[Theorem 3.2]{Shi25} classifies all basic tilting complexes up to shifts and mutation.

\section{Tilting fans determined by split ADE highest-root cycles}\label{sec:ade-fibres}

Theorem \ref{thm:ade-cycle-fan} identifies the tilting fan with the affine ADE alcove fan and its parabolic restrictions, and gives the Cayley graph description for full cycles. Mutations correspond to facet adjacency by Theorems \ref{thm:tilting-fan} and \ref{thm:facet-mutation}; Corollary \ref{cor:ade-mutation-symmetries} describes the geometric normalizer action. Section \ref{sec:full-ade-fan} determines the full-cycle fan using root tests for bricks, affine simple-root bases from tilting hearts and completeness. Section \ref{sec:ade-stein} then passes to Stein contractions by numerical descent, and closed-fibre invariance transfers the classification to $Y$.

\subsection{Intrinsic cycles and the alcove classification}

Let $\Delta$ be a connected finite ADE diagram with vertices $1,\ldots,m$, Cartan matrix $A$, and highest-root vector $a=(a_i)$. On $Q_\Delta:=\mathbb Z^\Delta$ put
\[
q(b):=\tfrac12 b^{\mathsf T}Ab,\qquad (b,c):=b^{\mathsf T}Ac.
\]
The nonnegative vectors with $q=1$ are the positive roots, and every positive root is at most $a$.

\begin{Def}\label{def:ade-cycle-scheme}
A \emph{split ADE highest-root cycle} of type $\Delta$ is a proper one-dimensional Cohen--Macaulay $k$-scheme $Z$ with the following properties.
\begin{enumerate}
\item The reduction $Z_{\mathrm{red}}$ is a nodal tree of curves $D_i\cong\mathbb P^1_k$ with dual graph $\Delta$ and $k$-rational nodes. At the generic point $\eta_i$ of $D_i$, we require $\ell_{\O_{Z,\eta_i}}(\O_{Z,\eta_i})=a_i$.
\item Put $\I_i=\ker(\O_Z\to\O_{D_i})$. For each $i$ with $a_i>1$, require
\begin{equation}\label{eq:ade-conormal}
\I_i/\I_i^2\cong\omega_{D_i}^{-1}.
\end{equation}
\end{enumerate}
Here \emph{split} refers only to the components and nodes of $Z_{\mathrm{red}}$ in (1). No retraction $Z\to Z_{\mathrm{red}}$ is required.
\end{Def}

The conormal condition specifies the canonical normal direction on each nonreduced component: $(\I_i/\I_i^2)^\vee\cong\omega_{D_i}$. We first establish the local normal forms and the successive layers of subcycles. Corollary \ref{cor:ade-cycle-euler} then gives their Euler characteristics and the dualizing sheaves of primary thickenings, while Lemma \ref{lem:ade-cycle-calculations} gives the cohomology and line-bundle properties of $Z$.

Here $\I_i$ is the ideal of $D_i$ in the entire $Z$, so its conormal also controls the gluing at the nodes. In type $A$, every $a_i=1$, and Cohen--Macaulayness makes $Z$ reduced. Thus split ADE highest-root cycles of type $A_m$ are exactly reduced split nodal chains of $m$ projective lines. Since $Z$ is proper of dimension one over $k$, it is projective \cite[Tag 0A26]{stacks-project}.

\begin{Lem}[Local normal forms]\label{lem:ade-local-structure}
Definition \ref{def:ade-cycle-scheme} implies the following normal forms.
\begin{enumerate}
\item At a closed point $p\in D_i$ away from the nodes, with $B=\O_{D_i,p}$, there is a $k$-algebra isomorphism
\[
\O_{Z,p}\cong B[t]/(t^{a_i}).
\]
\item At the rational node $p\in D_i\cap D_j$,
\begin{equation}\label{eq:ade-node}
\widehat{\O}_{Z,p}\cong k[[x,y]]/(x^{a_i}y^{a_j}).
\end{equation}
\end{enumerate}
Consequently $Z$ is Gorenstein and each $\I_i$ is locally principal along $D_i$.
\end{Lem}
\begin{proof}
For $a_i>1$, Nakayama makes $\I_i$ locally principal, with a possibly zero-divisor generator $t$.

For (1), put $A=\O_{Z,p}$ and $B=\O_{D_i,p}$. If $a_i=1$, purity gives $A=B$. Otherwise the nilradical is $(t)$, and purity gives $t^{a_i}=0$. Lifting an affine coordinate and inverting the denominators gives a section $B\to A$, hence a surjection $B[t]/(t^{a_i})\to A$. It is generically an isomorphism by the prescribed length; its kernel is zero since the source is $B$-free.

For (2), if $a_i=a_j=1$, purity gives the reduced node. Otherwise choose $x\in\I_i$, $y\in\I_j$ lifting the branch parameters, using a generator whenever the corresponding multiplicity exceeds one. They generate the maximal ideal, giving a surjection
\[
R:=k[[X,Y]]\twoheadrightarrow\widehat{\O}_{Z,p},
\qquad (X,Y)\longmapsto(x,y).
\]
The Cohen--Macaulay quotient has an unmixed height-one kernel, hence a principal ideal. Its radical is $(XY)$, so the kernel is $(X^rY^s)$. Flatness of completion preserves the branch lengths, giving $r=a_i$, $s=a_j$ and \eqref{eq:ade-node}.

Both local forms are hypersurfaces and make the component ideals principal. These properties descend from completion, proving the final assertions.
\end{proof}

For $0\le b\le a$, define the subcycle $Z_b\subseteq Z$ and the primary thickenings $B_{i,n}$ by
\[
\O_{Z_b}:=\O_Z\Big/\prod_i\I_i^{b_i},
\qquad B_{i,n}:=Z_{ne_i}=V(\I_i^n)\quad(1\le n\le a_i).
\]
Thus $Z_0=\varnothing$ and $Z_a=Z$. Write $p_{ij}:=D_i\cap D_j$ for adjacent vertices.

\Needspace{11\baselineskip}
\begin{Lem}[Subcycles and successive layers]\label{lem:ade-subcycle-filtration}
Let $0\le b\le a$.
\begin{enumerate}
\item If $b\ne0$, then $Z_b$ is Gorenstein. Its generic length along $D_i$ is $b_i$ for each $i$ with $b_i>0$.
\item If $b+e_j\le a$, restriction gives an exact sequence
\begin{equation}\label{eq:ade-cycle-step}
0\longrightarrow
\omega_{D_j}^{-b_j}\Bigl(-\sum_{i\sim j}b_i p_{ij}\Bigr)
\longrightarrow\O_{Z_{b+e_j}}
\longrightarrow\O_{Z_b}\longrightarrow0.
\end{equation}
Its kernel is a line bundle on $D_j$ of degree $(b,e_j)$.
\end{enumerate}
\end{Lem}
\begin{proof}
For (1), the local forms of Lemma \ref{lem:ade-local-structure} identify $Z_b$ away from the nodes with $B[t]/(t^{b_i})$ when $b_i>0$, and at a node with completed local ring $k[[x,y]]/(x^{b_i}y^{b_j})$ when $b_i+b_j>0$. These are hypersurfaces, hence Gorenstein, and have the stated generic lengths.

For (2), the same local forms show that the kernel in \eqref{eq:ade-cycle-step} is a line bundle on $D_j$. The $b_j$-th conormal power contributes $\omega_{D_j}^{-b_j}$, and the neighbouring branches require vanishing of order $b_i$ at $p_{ij}$. When $b_j>0$, one has $a_j>1$, so the conormal identification \eqref{eq:ade-conormal} applies. The resulting degree is $2b_j-\sum_{i\sim j}b_i=(b,e_j)$. In particular, the successive layers on $B_{i,n}$ are $\omega_{D_i}^{-r}$ for $0\le r<n$.
\end{proof}

\begin{Cor}[Euler characteristics and dualizing sheaves]\label{cor:ade-cycle-euler}
Let $0<b\le a$.
\begin{enumerate}
\item For a line bundle $L$ on $Z_b$, of degree $d_i$ on each component $D_i$ in its support, one has
\begin{equation}\label{eq:ade-cycle-euler}
\chi_k(Z_b,L)=q(b)+\sum_{b_i>0}b_i d_i.
\end{equation}
\item For $1\le n\le a_i$ and a line bundle $L$ on $B_{i,n}$ with $\deg(L|_{D_i})=d_i$, one has
\[
\chi_k(B_{i,n},L)=n(n+d_i),
\qquad
\omega_{B_{i,n}}|_{D_i}\cong\omega_{D_i}^{\otimes n}.
\]
\end{enumerate}
\end{Cor}
\begin{proof}
For (1), apply Lemma \ref{lem:ade-subcycle-filtration}(2) along a filtration from $Z_0$ to $Z_b$, tensoring each step with the restriction of $L$. Since $\chi_k(\O_{\mathbb P^1}(d))=d+1$ and $q(c+e_j)-q(c)=1+(c,e_j)$, induction gives \eqref{eq:ade-cycle-euler}.

For (2), take $b=ne_i$ in (1), so $q(b)=n^2$. Duality gives $\chi_k(\omega_{B_{i,n}})=-\chi_k(\O_{B_{i,n}})=-n^2$, and the same formula makes the degree of $\omega_{B_{i,n}}|_{D_i}$ equal to $-2n$. Since $D_i\cong\mathbb P^1_k$, this proves the asserted isomorphism.
\end{proof}

For the degree and nefness criteria, compare \cite[Lemmas 3.4.3--3.4.5]{VdB04b}. The proof below uses the intrinsic cycle structure over the given field $k$.

\begin{Lem}[Line bundles on the cycle]\label{lem:ade-cycle-calculations}
One has
\[
H^0(Z,\O_Z)=k,\quad H^1(Z,\O_Z)=0,\quad
\Pic Z\cong\mathbb Z^\Delta,
\]
where the last isomorphism is given by multidegree. Nef line bundles are globally generated and have zero first cohomology.
\end{Lem}
\begin{proof}
For the split nodal tree $T:=Z_{\mathrm{red}}$, the assertions follow from the normalization sequences and induction on leaves: the components are $\mathbb P^1_k$ and the nodes are $k$-rational.

Put $\mathbf1:=(1,\ldots,1)$, so $Z_{\mathbf1}=T$ and $q(\mathbf1)=1$. For a positive root $\mathbf1\le b<a$, the inequality $(b,a-b)=-q(a-b)<0$ gives a $j$ with $b_j<a_j$ and $(b,e_j)<0$. Positive definiteness and integrality give $q(b+e_j)=2+(b,e_j)=1$. Iteration gives a root filtration from $T$ to $Z$.

Along this filtration, \eqref{eq:ade-cycle-step} has kernel $\O_{D_j}(-1)$, which is square-zero since $b_j\ge1$. Its acyclicity preserves $H^0(\O)$ and $H^1(\O)$, while the square-zero unit sequence preserves $\Pic$ \cite[Tag 0C6R]{stacks-project}.

For nef $L$ of multidegree $d$, the twisted kernels are $\O_{D_j}(d_j-1)$, whose $H^1$ vanishes. Hence
\[
H^1(Z,L)=0,\qquad
H^0(Z,L)\twoheadrightarrow H^0(T,L|_T).
\]
Global generation follows from that on $T$ and Nakayama's lemma.
\end{proof}

For $J\subseteq\Delta$, set $I:=\Delta\setminus J$. Let $L_J$ be the unique line bundle on $Z$ of degree zero on $D_j$ for $j\in J$ and degree one on $D_i$ for $i\in I$. It is globally generated by Lemma \ref{lem:ade-cycle-calculations}. Define
\[
g_J:Z\longrightarrow C_J
\]
to be the Stein contraction of the morphism given by its complete linear system. Exactly the components indexed by $J$ are contracted. The target is projective over $k$. The cases $J=\varnothing$ and $J=\Delta$ give $C_J=Z$ and $C_J=\Spec k$, respectively.

On $\mathbb R^I$ define the arrangement
\begin{equation}\label{eq:restricted-affine-ade-arrangement}
\mathscr A_{\Delta,J}:=
\left\{
\sum_{i\in I}b_i z_i=n
\ \middle|\
\beta=\sum_i b_i\alpha_i\in\Phi_\Delta^+,\quad
n\in\mathbb Z,\quad (b_i)_{i\in I}\ne0
\right\}.
\end{equation}
Duplicate hyperplanes are identified, and the restricted coefficient vectors retain their actual multiplicities. Write $\widetilde\Delta$ for the affine Dynkin diagram and $\operatorname{Cham}(\widetilde\Delta,J)$ for the chambers of its restricted Tits cone. Rank-one slicing identifies them with the chambers of $\mathscr A_{\Delta,J}$, the $J$-alcoves \cite[Theorem 1.12, Definitions 2.5--2.6 and Proposition 2.7]{IWTC}. The level coordinates are identified explicitly in Section \ref{sec:ade-stein}.

\Needspace{14\baselineskip}
\begin{Thm}[Tilting fans and affine alcoves]\label{thm:ade-cycle-fan}
Assume \eqref{eq:geometric-setup} and $C\cong C_J$ over $k$, where $Z$ is a split ADE highest-root cycle of type $\Delta$ and $J\subseteq\Delta$.
\begin{enumerate}
\item In the rank-and-degree coordinates of Lemma \ref{lem:ade-numerical-coordinates}, the assignment $U\mapsto\Delta(U)$ is a bijection from $\ptilt Y\setminus\{0\}$ onto the nonempty faces of the closed $J$-alcoves of $\mathscr A_{\Delta,J}$. On tilting bundles this gives $\tilt Y\cong\operatorname{Cham}(\widetilde\Delta,J)$.
\item If $J=\varnothing$, the underlying undirected mutation graph is isomorphic to $\operatorname{Cay}(W_{\mathrm{aff}}(\Delta),S_{\mathrm{aff}})$, where $S_{\mathrm{aff}}$ is the set of affine simple reflections.
\end{enumerate}
\end{Thm}

For general $J$, identify $\mathbb R^I$ with the subspace $\{z_j=0\mid j\in J\}\subseteq\mathbb R^\Delta$. The arrangement $\mathscr A_{\Delta,J}$ is obtained by intersecting the hyperplanes of $\mathscr A_{\Delta,\varnothing}$ with this subspace, discarding empty and whole-subspace intersections and identifying duplicates. This parabolic restriction need not be a Coxeter arrangement. Its alcoves carry the parabolic labels and normalizer symmetries described in Section \ref{sec:ade-mutation-symmetries}.

Corollary \ref{cor:duval-cdv-tilting} applies this classification to partial crepant resolutions of Du Val singularities and characteristic-zero cDV singularities.

\subsection{The full ADE fan}\label{sec:full-ade-fan}

To determine the full-cycle fan, we first prove root tests for bricks and then apply them to the simple objects of tilting hearts. The root tests use only the local normal forms, the conormal degrees and duality on the curve. The field $k$ remains arbitrary.

For a pure sheaf $P$, write $P^\dagger:=\mathcal Hom_Z(P,\omega_Z)$.

\begin{Lem}[Finite-length modifications]\label{lem:ade-finite-length-invariance}
Let $P,Q$ be pure one-dimensional sheaves. The index
\[
\chi_k\mathcal Hom_Z(P,Q)-\chi_k(P^\dagger\otimes Q)
\]
is unchanged when either argument is replaced by a subsheaf with finite-length quotient. If $\rho(P)+\rho(Q)\le a$, the same holds for
\[
\chi_k\mathcal Hom_Z(P,Q)-\chi_k\mathcal Ext_Z^1(P,Q).
\]
\end{Lem}
\begin{proof}
Let $M,N$ be the completed stalks of $P,Q$ at a closed point $p$. In the local form $A_0=R/(f)$ of Lemma \ref{lem:ade-local-structure}, $R$ is regular of dimension two. Equal-rank free resolutions of length one \cite[Tag 00O7]{stacks-project} give, naturally in both variables,
\[
\Ext_R^1(M,N)
\cong\Ext_R^1(M,R)\otimes_RN
\cong\Hom_{A_0}(M,A_0)\otimes_{A_0}N.
\]
Fixing a local generator of $\omega_Z$ identifies the maps on $\Ext_R^1$ with those on the completed stalks of $P^\dagger\otimes Q$.

For a finite-length modification, let $\alpha,\beta$ be the induced maps on $\Hom_R$ and $\Ext_R^1$, in their functorial directions, and put $\operatorname{ind}u:=\ell\operatorname{coker}u-\ell\ker u$. The long exact sequences identify $\operatorname{ind}\alpha-\operatorname{ind}\beta$, in the first and second variables respectively, with
\[
\begin{aligned}
\ell\Ext_R^1(K,N)-\ell\Ext_R^2(K,N)&=0,\\
\ell\Hom_R(M,K)-\ell\Ext_R^1(M,K)&=0,
\end{aligned}
\]
where $K$ is the quotient. Both vanishings follow from the equal-rank resolutions and $\Hom_R(K,R)=\Ext_R^1(K,R)=0$. Thus $\alpha,\beta$ have equal indices. Naturality, preservation of length under completion and the weights $[\kappa(p):k]$ give the first assertion.

Under $\rho(P)+\rho(Q)\le a$, every $R$-extension has a pure middle term of branch lengths at most the orders of $f$. Hence $f$ annihilates it generically and therefore everywhere by purity. The natural map $\Ext_{A_0}^1(M,N)\to\Ext_R^1(M,N)$ is thus an isomorphism, so the same argument proves the second assertion.
\end{proof}

\begin{Lem}[Cartan estimates]\label{lem:ade-cartan-estimates}
For pure one-dimensional sheaves $P,Q$, let $\hom(P,Q)$ and $\operatorname{ext}^1(P,Q)$ denote the $k$-dimensions of $\Hom_Z(P,Q)$ and $\Ext_Z^1(P,Q)$, respectively. Then
\begin{equation}\label{eq:ade-bounded-global}
\begin{aligned}
(\rho(P),\rho(Q))&\le\hom(P,Q)+\hom(Q,P),\\
(\rho(P),\rho(Q))&\ge\hom(P,Q)-\operatorname{ext}^1(P,Q)
\quad\text{if }\rho(P)+\rho(Q)\le a.
\end{aligned}
\end{equation}
\end{Lem}
\begin{proof}
Put $H=\mathcal Hom_Z(P,Q)$, $T=P^\dagger\otimes Q$ and $E=\mathcal Ext_Z^1(P,Q)$. We first compute
\begin{equation}\label{eq:ade-cartan-indices}
\begin{aligned}
\chi_k(H)-\chi_k(T)&=(\rho(P),\rho(Q)),\\
\chi_k(H)-\chi_k(E)&=(\rho(P),\rho(Q))
\quad\text{if }\rho(P)+\rho(Q)\le a.
\end{aligned}
\end{equation}
At the generic points, each pure sheaf is a sum of cyclic modules, hence is generically isomorphic to a sum of $\O_{B_{i,n}}$. A common coherent lattice gives finite-length modifications. Lemma \ref{lem:ade-finite-length-invariance} and biadditivity therefore reduce the calculation to $P=\O_{B_{i,n}}$, $Q=\O_{B_{j,m}}$; the bound, when required, holds for every pair of summands.

Finite duality gives $P^\dagger=\omega_{B_{i,n}}$. If $i=j$, put $s=\min(n,m)$. The local forms identify $H=\I_i^{m-s}/\I_i^m$ and $T=\omega_{B_{i,n}}|_{B_{i,s}}$. These are line bundles of degrees $2(m-s)$ and $-2n$ by Lemma \ref{lem:ade-subcycle-filtration}(2) and Corollary \ref{cor:ade-cycle-euler}(2). The latter gives $\chi_k(H)-\chi_k(T)=2nm$. If $i\ne j$, then $H=0$, while $T$ has length $nm$ at an adjacent node and vanishes otherwise. This proves the first identity.

For the second, the bound gives $f=ghc$ in the completed equations $P=R/(g)$, $Q=R/(h)$. Since $\I_i^n=(g)\cong R/(hc)$, the completed stalks of $L:=\I_i^n\otimes Q$ are $R/(h)$; faithful flatness makes $L$ a line bundle on $B_{j,m}$. As $ne_i+e_j\le a$, Lemma \ref{lem:ade-subcycle-filtration}(2) applied to
\[
L|_{D_j}\cong\I_i^n/(\I_i^n\I_j)
=\ker\bigl(\O_{Z_{ne_i+e_j}}\longrightarrow\O_{B_{i,n}}\bigr)
\]
gives $\deg(L|_{D_j})=n(e_i,e_j)$. Applying $\mathcal Hom_Z(-,Q)$ to the ideal sequence yields
\[
0\longrightarrow H\longrightarrow Q\longrightarrow L^\vee
\longrightarrow E\longrightarrow0,
\]
where the dual is taken on $B_{j,m}$. Corollary \ref{cor:ade-cycle-euler}(2) gives
\[
\begin{aligned}
\chi_k(H)-\chi_k(E)&=\chi_k(Q)-\chi_k(L^\vee)\\
&=m^2-m\bigl(m-n(e_i,e_j)\bigr)=nm(e_i,e_j).
\end{aligned}
\]

For general $P,Q$, Cohen--Macaulay biduality and curve duality \cite[Tags 0B5A and 0BS2]{stacks-project} give $P\cong P^{\dagger\dagger}$ and $h^1(T)=\hom(Q,P)$. Together with $\operatorname{ext}^1(P,Q)=h^1(H)+h^0(E)$ from local-to-global Ext, \eqref{eq:ade-cartan-indices} gives both bounds by discarding nonnegative cohomology dimensions.
\end{proof}

\begin{Prop}[Three root tests for bricks]\label{prop:ade-root-tests}
Let $F,G$ be nonzero pure one-dimensional bricks on $Z$.
\begin{enumerate}
\item $\End_Z(F)=k$ and $q(\rho(F))=1$.
\item If $\Hom_Z(F,G)=\Hom_Z(G,F)=0$, then $(\rho(F),\rho(G))\le0$.
\item If $\Ext_Z^1(F,G)=0$, then $(\rho(F),\rho(G))\ge0$.
\end{enumerate}
\end{Prop}
\begin{proof}
Put $\beta=\rho(F)$ and $e=\dim_k\End_Z(F)$. Choose line bundles $L_i$ of multidegree $e_i$ by Lemma \ref{lem:ade-cycle-calculations}. The division algebra $\End_Z(F)$ acts on $H^j(F)$ and $H^j(F\otimes L_i)$, so $e$ divides their Euler characteristics. Riemann--Roch (Lemma \ref{lem:numerical-devissage}) gives $e\mid\beta_i$. Since $\beta/e$ is a nonzero integral vector, positive definiteness and Lemma \ref{lem:ade-cartan-estimates} give
\[
e^2\le q(\beta)\le e.
\]
Thus $e=q(\beta)=1$, proving (1).

The upper bound of Lemma \ref{lem:ade-cartan-estimates} gives (2). For (3), put $\gamma=\rho(G)$. If $(\beta,\gamma)<0$, then (1) gives
\[
0<q(\beta+\gamma)=2+(\beta,\gamma)\le1.
\]
Hence $\beta+\gamma$ is a positive root and $\beta+\gamma\le a$. The lower bound of Lemma \ref{lem:ade-cartan-estimates} would then give $(\beta,\gamma)\ge\hom(F,G)\ge0$, a contradiction.
\end{proof}

We now identify the affine root lattice attached to $Z$. By projectivity and Lemma \ref{lem:ade-cycle-calculations}, $Z\to\Spec k$ satisfies the standing setup. Put $\mathsf N:=K_0(\Coh Z)$ and choose
\[
\delta=[\kappa(p)],\qquad \gamma_i=[\O_{D_i}(-1)],
\]
where $p$ is a $k$-rational point. Let $\L_i$ have multidegree $e_i$. Corollary \ref{cor:rank-determinant} shows that $[\O_Z]$ and $[\L_i]-[\O_Z]$ form an integral basis of $K_0(\per Z)$. Their Euler matrix against $\delta,\gamma_i$ is $\operatorname{diag}(1,-I_m)$. Since the two Grothendieck groups are free of the same rank, $\delta,\gamma_i$ form an integral basis as well, and the Euler pairing gives
\[
\mathsf N\xrightarrow[\sim]{(\chi_k,\rho)}\mathbb Z\oplus Q_\Delta,
\qquad K_0(\per Z)\cong\mathsf N^\vee,
\qquad \langle(r,d),(h,\beta)\rangle=rh-d\cdot\beta.
\]
Define the Cartan form on $\mathsf N$ by $(x,y):=\rho(x)^{\mathsf T}A\rho(y)$ and $q(x):=\tfrac12(x,x)$. The basis
\[
\alpha_0=[\O_Z]=\delta+\sum_i a_i\gamma_i,\qquad
\alpha_i=-\gamma_i,\qquad
\delta=\alpha_0+\sum_i a_i\alpha_i
\]
has the affine Cartan matrix. Thus $\mathsf N$ is the affine root lattice, with radical $\mathbb Z\delta$ and real roots $n\delta+\beta$ for $\beta\in\Phi_\Delta$. Under Euler duality, evaluation at $\delta$ is rank.

\begin{Rem}
If $i\colon Z\hookrightarrow S$ realizes $Z$ in a regular $k$-surface with $D_i\cdot D_j=-A_{ij}$, the surface Euler formula gives
\[
([F],[G])=\chi_k\RHom_S(i_*F,i_*G)
\qquad(F,G\in\Coh Z).
\]
For the minimal resolution of a complex Kleinian singularity, this is the Euler form of the supported $2$-Calabi--Yau category \cite[Section 2.1]{Bri09}. The intrinsic construction above requires no ambient surface.
\end{Rem}

\Needspace{11\baselineskip}
\begin{Lem}[An elementary root-basis criterion]\label{lem:ade-obtuse-basis}
Suppose $s_0,\ldots,s_m$ is an integral basis of $\mathsf N$ such that
\[
q(s_i)=1,\qquad (s_i,s_j)\le0\ (i\ne j),\qquad
\delta=\sum_i r_i s_i\quad(r_i>0).
\]
Then every $b\in\mathsf N$ with $q(b)=1$ has coefficients of one sign in this basis. Consequently,
\[
\{\ell\mid\ell(\delta)=1,\ \ell(s_i)>0\text{ for all }i\}
\]
is one alcove of the affine ADE arrangement.
\end{Lem}
\begin{proof}
If $b=u-v$ has nonzero positive and negative parts in the stated basis, neither part is a multiple of $\delta$, since their coefficient supports are disjoint and the coefficients of $\delta$ are all positive. Thus $q(u),q(v)\ge1$, while $(u,v)\le0$. Hence
\[
1=q(b)=q(u)+q(v)-(u,v)\ge2,
\]
a contradiction. Therefore no real-root hyperplane meets the displayed open simplex. Its bounding hyperplanes are themselves real-root hyperplanes, so the simplex is an entire alcove.
\end{proof}

\begin{Prop}[Full-cycle fan]\label{prop:full-ade-fan}
For $Z\to\Spec k$, the rank-one tilting simplices are exactly the closed alcoves of $\mathscr A_{\Delta,\varnothing}$, and the pretilting cones are the cones over all their faces.
\end{Prop}
\begin{proof}
\emph{The simple-root basis.}
Let $\E=\bigoplus_{i=0}^m\E_i$ be basic tilting, with heart simples $S_i=F_i[\epsilon_i]$, $\epsilon_i\in\{0,1\}$, and put $s_i=[S_i]$. Since every zero-dimensional sheaf belongs to $\Fac\E$, torsion-pair orthogonality and simplicity make each $F_i$ pure unless $S_i$ itself is zero-dimensional. The latter is impossible: $\Hom(\E_j,Q)\ne0$ for every $j$ and every nonzero zero-dimensional $Q$, contradicting simple--projective orthogonality since $m+1\ge2$. Proposition \ref{prop:ade-root-tests}(1) therefore gives $\End(S_i)=k$ and $q(s_i)=1$.

By Lemma \ref{lem:euler-dual-bases}, the bases $[\E_i]$ and $s_i$ are Euler dual, so
\begin{equation}\label{eq:ade-positive-point-class}
\delta=\sum_i r_i s_i,\qquad r_i:=\rk\E_i>0.
\end{equation}

For distinct simples with equal shifts, Hom orthogonality and Proposition \ref{prop:ade-root-tests}(2) give $(s_i,s_j)\le0$. For opposite shifts, order them as $S_i=F_i$, $S_j=F_j[1]$. Orthogonality says $\Ext_Z^1(F_i,F_j)=0$, so Proposition \ref{prop:ade-root-tests}(3), with the sign from the shift, again gives $(s_i,s_j)\le0$. Lemma \ref{lem:ade-obtuse-basis} now makes the $s_i$ an affine simple-root basis.

\emph{Alcoves and exhaustivity.}
Corollary \ref{cor:cone-from-heart} identifies $C^\circ(\E)$ with positivity on the $s_i$, so Lemma \ref{lem:ade-obtuse-basis} makes $C^\circ(\E)\cap H$ one entire alcove. The affine ADE alcoves form a locally finite family: on a compact set, only finitely many integral levels of the finitely many finite roots can occur. The tilting simplices are therefore locally finite, so Theorem \ref{thm:complete-finite} gives $H=\bigcup_\E\Delta(\E)$. Hence every alcove occurs. Completion by Proposition \ref{Bongcomp} gives all faces.
\end{proof}

\subsection{Stein contractions and Tits cone intersections}\label{sec:ade-stein}

\begin{Lem}[Numerical descent of rigid bundles]\label{lem:ade-rigid-descent}
Let $g:Z\to C$ be a projective Stein contraction, with contracted set $J$. Then $\mathbb Rg_*\O_Z=\O_C$, and pullback gives a bijection
\[
\ptilt C\ \longleftrightarrow\
\{F\in\ptilt Z\mid d_j(F)=0\text{ for every }j\in J\}.
\]
It preserves indecomposable summands, compatibility, ranks and degrees on surviving components.
\end{Lem}
\begin{proof}
The Stein condition gives $g_*\O_Z=\O_C$. The map is an isomorphism away from finitely many points, so $R^1g_*\O_Z$ has finite support. Leray and $H^1(Z,\O_Z)=0$ make it zero. Projection formula therefore identifies Hom and Ext between bundles downstairs with those between their pullbacks.

Conversely, let $F$ be pretilting of rank $r$ and of zero $J$-degree. The same finite-support argument gives $R^1g_*\mathcal End(F)=0$. For each $c\in C$, flat base change to $A=\widehat{\O}_{C,c}$ makes $F_A$ pretilting on $Z_A:=Z\times_C\Spec A$. This morphism satisfies the standing setup, and every one-dimensional component of its closed fibre is contracted by $g$. Corollary \ref{cor:numerically-trivial-pretilting} thus gives
\[
F_A\cong\O_{Z_A}^{\oplus r}.
\]
Flat base change and faithful flatness of completion show that $g_*F$ is locally free and that $g^*g_*F\to F$ is an isomorphism. Full faithfulness proves the remaining assertions.
\end{proof}

\begin{Rem}\label{rem:ade-stein-uniqueness}
The descent proof uses only $H^1(Z,\O_Z)=0$ and numerical triviality of pretilting bundles, so it applies to projective Stein contractions of arbitrary connected projective curves with this vanishing.

In particular, $C_J$ depends only on $Z$ and $J$: any such contraction descends $L_J$ to a line bundle $M$ of positive degree on every surviving component, hence ample. Projection formula identifies its section ring with that of $L_J$, giving
\[
C\cong\operatorname{Proj}\bigoplus_{n\ge0}H^0(Z,L_J^{\otimes n})\cong C_J.
\]
\end{Rem}

\begin{Lem}[Rank and degree coordinates]\label{lem:ade-numerical-coordinates}
Assume \eqref{eq:geometric-setup} and $C\cong C_J$. Multidegree and rank-and-degree identify
\[
\Pic(Y)\cong\mathbb Z^I,\qquad
K_0(\per Y)\cong\mathbb Z\oplus\mathbb Z^I.
\]
Thus $H\cong\mathbb R^I$ with coordinates $z_i=d_i$, and $n=|I|+1$.
\end{Lem}
\begin{proof}
Every line bundle on $Z$ is pretilting by Lemma \ref{lem:ade-cycle-calculations}. Lemma \ref{lem:ade-rigid-descent} therefore identifies $\Pic(C_J)$ with the zero-$J$-degree subgroup $\mathbb Z^I\subseteq\Pic(Z)$. Picard restriction (Lemma \ref{lem:picard-restriction}) and Corollary \ref{cor:rank-determinant} give the assertions.
\end{proof}

\begin{proof}[Proof of Theorem \ref{thm:ade-cycle-fan}]
\emph{(1)} Closed-fibre invariance reduces the classification to $C_J$, which satisfies the standing setup by Lemma \ref{lem:ade-rigid-descent}. Put $V_J:=\{d_j=0\mid j\in J\}$ and $H_J:=V_J\cap\{\rk=1\}$. Numerical descent gives
\begin{equation}\label{eq:ade-fan-restriction}
g_J^*\Sigma(C_J)=\{\sigma\in\Sigma(Z)\mid\sigma\subseteq V_J\}.
\end{equation}
Since $H_J$ is a flat of the full root arrangement, the sign vector of a relatively open face of $\mathscr A_{\Delta,J}$ extends uniquely to the full arrangement: walls containing $H_J$ have sign zero, and disjoint walls have fixed nonzero signs. The zero equations force the ambient face into $H_J$, so it coincides with the restricted face; conversely, restricting an ambient face contained in $H_J$ gives a restricted face. Taking closures, Proposition \ref{prop:full-ade-fan} identifies these faces with rank-one slices of pretilting cones. Corollary \ref{cor:class-determines-pretilting} gives injectivity, and tilting completion identifies the maximal faces with tilting bundles.

\emph{(2)} For $J=\varnothing$, let $A_0$ be the fundamental alcove and label bundles by $\Delta(\E_w)=\overline{wA_0}$. The affine Coxeter labelling \cite[Theorems 1.12 and 1.20, Proposition 2.7]{IWTC} identifies alcove adjacency with $wA_0\longleftrightarrow wsA_0$ for $s\in S_{\mathrm{aff}}$. Theorems \ref{thm:tilting-fan} and \ref{thm:facet-mutation} therefore make $w\mapsto\E_w$ the asserted graph isomorphism.
\end{proof}

For comparison with Iyama--Wemyss \cite[Definitions 2.5--2.6 and Proposition 2.7]{IWTC}, put
\[
\mathsf N_J:=\mathsf N/\sum_{j\in J}\mathbb Z\alpha_j,\qquad J_{\mathrm{aff}}:=J.
\]
We view $J_{\mathrm{aff}}$ as a subset of the affine vertices, with vertex $0$ retained. The dual space is
\[
\Theta_{J_{\mathrm{aff}}}=\{\ell\mid\ell(\alpha_j)=0\ (j\in J)\}=V_J.
\]
Its rank-one slice is $\operatorname{Level}(J_{\mathrm{aff}})$. In coordinates,
\[
\vartheta_0=r-\sum_i a_id_i,\qquad \vartheta_i=d_i,\qquad
\vartheta(\delta)=r.
\]
The walls are the nonempty proper zero sets of the images of affine real roots in $\mathsf N_J$. For $J\ne\varnothing$, the quadratic form does not descend, since $q(\alpha_j)=1$. If a restricted finite root is $c\gamma$, its walls are $\gamma\cdot z=n/c$; primitive normalization must also rescale the constants.

\begin{Ex}
For $\Delta=D_4$, retain only the central component and contract the three outer components. The restricted positive-root coefficients are $1$ and $2$, so the rank-one walls are
\[
z\in\mathbb Z\quad\text{or}\quad 2z\in\mathbb Z,
\qquad\text{equivalently}\qquad z\in\tfrac12\mathbb Z.
\]
This is the $D_4$ level arrangement of Iyama--Wemyss \cite[Example 2.9]{IWTC}. The primitive classes at $z=0,\frac12,1$ are $(1,0),(2,1),(1,1)$, respectively. Each adjacent tilting bundle consequently has summand ranks $1$ and $2$, and there are two tilting bundles modulo the Picard translation $z\mapsto z+1$.
\end{Ex}

\begin{Cor}[Du Val and cDV partial resolutions]\label{cor:duval-cdv-tilting}
Let $\pi\colon Y\to\Spec R$ be a projective crepant partial resolution, where $Y$ is normal and $R$ is a complete local $k$-algebra with algebraically closed residue field $k$. Assume that $R$ is Du Val, or that $\operatorname{char}k=0$ and $R$ is cDV. Let $\Delta$ be the Dynkin diagram of $R$ in the surface case and of a general Du Val section in the cDV case, and let $J$ index the curves contracted from the minimal surface resolution to the corresponding partial resolution.

Then Theorem \ref{thm:ade-cycle-fan} applies with these data. In particular, in rank-and-degree coordinates there is a bijection
\[
\tilt Y\xrightarrow{\ \sim\ }\operatorname{Cham}(\widetilde\Delta,J),
\qquad \E\longmapsto\operatorname{Int}\Delta(\E).
\]
\end{Cor}
\begin{proof}
These morphisms satisfy \eqref{eq:geometric-setup}; for cDV singularities see \cite[Sections 1.2 and 2.2]{Wem23survey}.

First suppose $R$ is Du Val, and write $\widetilde S\xrightarrow{h}Y\to\Spec R$ for its minimal resolution. By Artin's description of the closed fibre \cite{Art66}, it is the highest-root cycle $Z=\sum_i a_iD_i$. This is a Cohen--Macaulay divisor on the regular surface $\widetilde S$, and for $a_i>1$ one has $\I_i/\I_i^2\cong\O_{D_i}(-D_i)\cong\O_{D_i}(2)$. Hence $Z$ satisfies Definition \ref{def:ade-cycle-scheme}.

Put $C=Y\times_Rk$ and $g=h|_Z\colon Z\to C$. Since $Rh_*\O_{\widetilde S}=\O_Y$ and $R^1h_*$ is right exact, the global generation of $\m\O_{\widetilde S}$ gives $R^1h_*(\m\O_{\widetilde S})=0$. Pushing forward the restriction sequence therefore gives a surjection $\O_C\twoheadrightarrow g_*\O_Z$. Its kernel has finite support and no global sections, since $H^0(C,\O_C)=H^0(Z,\O_Z)=k$. Thus $g_*\O_Z=\O_C$, and Remark \ref{rem:ade-stein-uniqueness} identifies $C$ with $C_J$.

For cDV $R$ in characteristic zero, Reid's general elephant gives a general $t\in\m$ such that $R/(t)$ is Du Val and $Y_t:=Y\times_R R/(t)$ is a crepant partial resolution \cite[Theorem 1.14]{Rei83}. Its scheme-theoretic closed fibre is again $C$, so the surface case identifies $C$ with $C_J$. Theorem \ref{thm:ade-cycle-fan} now applies in both cases, and Corollary \ref{cor:class-determines-pretilting} gives the stated bijection.
\end{proof}

Over $\mathbb C$, compare Shimpi's classification of intermediate hearts for flopping contractions over isolated cDV singularities and for partial Kleinian resolutions \cite[Theorem 4.6 and the Introduction]{Shi26}.

\subsection{Numerical exchange and normalizer symmetries}\label{sec:ade-mutation-symmetries}

The fan of Theorem \ref{thm:ade-cycle-fan} is complete, so Theorem \ref{thm:complete-finite} and Corollary \ref{cor:complete-mutation-connectivity} give finiteness modulo line bundles and mutation connectivity. We next describe the exchange coefficients and the normalizer symmetries.

\begin{Cor}[Numerical exchange]\label{cor:ade-numerical-exchange}
Let $U\oplus \E_i$ and $U\oplus \E_i'$, with $U=\bigoplus_{k\ne i}\E_k$, correspond to adjacent $J$-alcoves. Write their adjacent-ray relation as
\[
[\E_i]+[\E_i']=\sum_{k\ne i}b_k[\E_k].
\]
Then the middle term of the exchange sequence in Theorem \ref{thm:facet-mutation} is $\bigoplus_{k\ne i}\E_k^{\oplus b_k}$. These coefficients are thus determined by the adjacent integral cones.
\end{Cor}
\begin{proof}
Theorem \ref{thm:facet-mutation} supplies the exchange sequence. Compare its class relation in the integral basis given by $U\oplus\E_i$.
\end{proof}

Write $W:=W_{\mathrm{aff}}(\Delta)$, and let $W_J$ be its standard parabolic subgroup generated by the vertices in $J$. By Iyama--Wemyss \cite[Theorem 1.12 and Proposition 2.7]{IWTC}, a $J$-alcove has a parabolic label $(x,K)$, where $K\subseteq\widetilde\Delta$, $x$ is minimal in $xW_K$, and $xW_Kx^{-1}=W_J$. We call $K$ its parabolic type.

\begin{Cor}[Normalizer symmetries]\label{cor:ade-mutation-symmetries}
Under the assumptions of Theorem \ref{thm:ade-cycle-fan}, the action on the restricted arrangement induces a free left action of
\[
G_J:=N_W(W_J)/W_J
\]
on $\tilt Y$ preserving mutation adjacency. Its orbits are precisely the subsets of a fixed parabolic type.
\end{Cor}
\begin{proof}
The normalizer acts on the restricted arrangement through $G_J$, preserving the rank-one level since $W$ fixes $\delta$. Its action on chambers has the stated properties by \cite[Lemmas 3.23--3.24]{GHLY26}. Transport this action through Theorem \ref{thm:ade-cycle-fan}, using Theorems \ref{thm:tilting-fan} and \ref{thm:facet-mutation} to identify mutation adjacency.
\end{proof}

\begin{Rem}[Wall types and geometric symmetries]
For $J=\varnothing$, use the alcove labelling of the proof of Theorem \ref{thm:ade-cycle-fan}(2). If $F_s$ is the wall of $A_0$ fixed by $s\in S_{\mathrm{aff}}$, the wall of type $s$ at $wA_0$ is $wF_s$, with reflection $wsw^{-1}$. Crossing it gives the mutation $\mu_s(\E_w):=\E_{ws}$ and hence the right action
\[
\E_w\cdot v:=\E_{wv}.
\]
The geometric left action of Corollary \ref{cor:ade-mutation-symmetries} is $g\cdot\E_w=\E_{gw}$ and preserves adjacency; the right action need not. These actions concern isomorphism classes of bundles, with no assertion of a faithful action by derived autoequivalences.

For general $J$, mutation may change parabolic type; the full wall-crossing structure is described by the Brink--Howlett groupoid \cite[Theorem 4.7]{GHLY26}.
\end{Rem}

\subsection{Explicit bundles in type \texorpdfstring{$A$}{A}}\label{sec:type-a}

Assume that $C$ is a split ADE highest-root cycle of type $A_m$ (Definition \ref{def:ade-cycle-scheme}), and write
\begin{equation}\label{eq:split-chain}
C=C_1\cup\cdots\cup C_m,\qquad C_i\cong\mathbb P^1_k,\qquad m\ge1,
\end{equation}
where consecutive components meet at one rational ordinary node and there are no other intersections. The following gives the explicit type $A$ form of Theorem \ref{thm:ade-cycle-fan}.

Under \eqref{eq:split-chain}, Picard and rank-and-degree coordinates are $\mathbb Z^m$ and $\mathbb Z^{m+1}$. On $H$ use the cumulative coordinates
\begin{equation}\label{eq:type-a-coordinates}
x_0:=0,\qquad x_j:=\sum_{i=1}^j d_i\quad(1\le j\le m).
\end{equation}
For $a\in\mathbb Z^m$, let $L_a$ denote the unique line bundle with cumulative degree vector $a$. Thus $L_0=\O_Y$ and $L_a\otimes L_b=L_{a+b}$; set $a_0=0$ when taking coordinate differences.

\begin{Thm}\label{thm:type-a-classification}
Assume \eqref{eq:split-chain}. Every indecomposable pretilting bundle is a line bundle. The simplices $\Delta(\E)$ are exactly the closures of the alcoves of the arrangement
\begin{equation}\label{eq:affine-a-arrangement}
x_j-x_i\in\mathbb Z\qquad(0\le i<j\le m),\qquad x_0=0.
\end{equation}
Explicitly, the basic tilting bundles are, without repetition,
\begin{equation}\label{eq:explicit-type-a-bundles}
\E_{b,\sigma}=\bigoplus_{j=0}^m
L_{b+e_{\sigma(1)}+\cdots+e_{\sigma(j)}},
\qquad b\in\mathbb Z^m,\quad\sigma\in\mathfrak S_m.
\end{equation}
Here $e_1,\ldots,e_m$ are the standard basis vectors in cumulative coordinates, and the sum for $j=0$ is empty. There are $m!$ such bundles modulo tensoring by line bundles.
\end{Thm}
\begin{proof}
Apply Theorem \ref{thm:ade-cycle-fan} with $\Delta=A_m$ and $J=\varnothing$. The closed alcoves are
\[
b+\{1\ge x_{\sigma(1)}\ge\cdots\ge x_{\sigma(m)}\ge0\}.
\]
Every primitive ray through an alcove vertex has rank one. An indecomposable pretilting bundle occurs in a tilting bundle, so its class is primitive and it is a line bundle, uniquely determined by its degree. The vertices give \eqref{eq:explicit-type-a-bundles}, and tensoring translates $b$, leaving precisely $m!$ possibilities.
\end{proof}

We retain the following criterion for the later flag-model arguments:
\begin{equation}\label{eq:line-compatibility}
\Ext_Y^1(L_a,L_b)=0=\Ext_Y^1(L_b,L_a)
\quad\Longleftrightarrow\quad
|(a_j-a_i)-(b_j-b_i)|\le1
\quad(0\le i<j\le m).
\end{equation}
By completion, compatibility means that $a,b$ lie in a common alcove. The inequalities say that $b-a$ has entries all in $\{0,1\}$ or all in $\{0,-1\}$; ordering the differing coordinates first gives such an alcove. Conversely, every root coordinate varies by at most one on a closed alcove.

\subsection{Minimal resolutions of complete toric surface singularities}\label{sec:toric-surfaces}

Let $k$ be an arbitrary field, let $N$ be a lattice of rank two with dual $M$, and let $\sigma\subset N_{\mathbb R}$ be a full dimensional strongly convex rational cone which is not regular. Set
\begin{equation}\label{eq:complete-toric-surface}
S:=\sigma^\vee\cap M,\qquad A:=k[S],\qquad
\mathfrak q:=(\chi^u\mid u\in S\setminus\{0\}),\qquad
R:=\widehat{A_{\mathfrak q}}.
\end{equation}
The toric variety is split and $A/\mathfrak q=k$. Let $\pi_{\mathrm{tor}}\colon\widetilde X_\sigma\to X_\sigma:=\Spec A$ be its minimal toric resolution, and let
\[
\pi\colon Y:=\widetilde X_\sigma\times_A R\longrightarrow\Spec R
\]
be its base change, with $m$ exceptional curves.

\begin{Prop}\label{prop:toric-split-fibre}
The morphism $\pi$ is the projective minimal resolution of $\Spec R$. It satisfies
\[
\mathbb R\pi_*\O_Y=\O_{\Spec R},
\]
and its scheme-theoretic closed fibre is a reduced split nodal chain of $m$ projective lines. Thus both \eqref{eq:geometric-setup} and \eqref{eq:split-chain} hold.
\end{Prop}
\begin{proof}
Hirzebruch--Jung theory gives a split nodal exceptional chain $E=\sum_iE_i$ with $E_i^2=-b_i\le-2$ \cite[Theorems 10.2.3, 10.2.5 and Corollary 10.4.9]{CLS11}. The subdivision is a sequence of star subdivisions, hence projective \cite[Proposition 11.1.6(c)]{CLS11}.

For the valency $t_i\le2$, one has $(-E)\cdot E_i=b_i-t_i\ge0$. These are all the complete invariant curves, so $\O_{\widetilde X_\sigma}(-E)$ is globally generated \cite[Theorem 6.3.12]{CLS11}. Since every nonzero $u\in S$ is positive on each interior ray,
\[
\Gamma(\widetilde X_\sigma,\O(-E))
=\bigoplus_{u\in S\setminus\{0\}}k\chi^u
=\mathfrak q.
\]
global generation yields
\begin{equation}\label{eq:toric-fibre-ideals}
\mathfrak q\O_{\widetilde X_\sigma}=\O_{\widetilde X_\sigma}(-E).
\end{equation}
Thus the scheme-theoretic closed fibre is exactly the reduced divisor $E$.

Toric vanishing gives $\mathbb R(\pi_{\mathrm{tor}})_*\O_{\widetilde X_\sigma}=\O_{X_\sigma}$ \cite[Theorem 9.2.5]{CLS11}. Both arguments are combinatorial and hold over any field. Flat base change gives the vanishing, closed fibre and projectivity for $\pi$. Excellence makes the completion map regular, so $Y$ is regular \cite[Tags 07QS and 07R6]{stacks-project}; the self-intersections remain $-b_i\le-2$, proving minimality.
\end{proof}

Theorem \ref{thm:type-a-classification} now applies: the tilting fan and bundle classification depend only on $m$, independently of the Hirzebruch--Jung coefficients $b_i$.

\begin{Ex}\label{ex:toric-one-fifth}
Let $k$ be any field, and take
\[
N=\mathbb Z^2+\mathbb Z\frac{(1,2)}5,
\qquad\sigma=\mathbb R_{\ge0}^2.
\]
Then $M=\{(u,v)\in\mathbb Z^2\mid u+2v\equiv0\pmod5\}$, and the semigroup ring is
\[
A=k[x^5,x^3y,xy^2,y^5],
\qquad R=\widehat{A_{\mathfrak q}},
\quad\mathfrak q=(x^5,x^3y,xy^2,y^5).
\]
Hirzebruch--Jung theory gives two exceptional curves with self-intersections $-2$ and $-3$.

This singularity is not Gorenstein. Indeed, the toric canonical-divisor criterion requires $(1,1)\in M$, whereas $1+2\not\equiv0\pmod5$; see \cite[Proposition 8.2.12]{CLS11}. In the cumulative degree coordinates \eqref{eq:type-a-coordinates}, the two basic tilting bundles modulo line bundles are
\[
\O_Y\oplus L_{(1,0)}\oplus L_{(1,1)},
\qquad
\O_Y\oplus L_{(0,1)}\oplus L_{(1,1)}.
\]
\end{Ex}

\part{An Auslander--McKay embedding over a normal Gorenstein base}\label{part:auslander-mckay}

Throughout Part \ref{part:auslander-mckay}, let $(R,\m,k)$ be a noetherian complete local normal Gorenstein domain and put $X:=\Spec R$. We consider projective birational morphisms $\pi\colon Y\to X$ from normal integral schemes satisfying
\begin{equation}\label{eq:crepant-model-hypotheses}
\dim(Y\times_Rk)\le1,\qquad
\pi^!\O_X\cong\O_Y.
\end{equation}
Let $\mathsf{Crep}_1(R)$ denote their category, with $X$-morphisms. Lemma \ref{lem:crepancy-and-duality} identifies the last condition with Cohen--Macaulayness of $Y$ and triviality of its canonical sheaf. In particular, these models are Gorenstein.

Section \ref{sec:general-reflexive-pushforward} first proves the structure-sheaf vanishing in Proposition \ref{prop:crepant-vanishing}. Thus Part \ref{part:tilting bundles} applies, and Proposition \ref{prop:canonical-tilting-general} supplies the canonical generator $\P_Y$. We then relate pretilting bundles to modifying modules, establish compatibility with modifying mutation, and prove a three-dimensional realization criterion from completeness of the tilting fan. Section \ref{sec:general-auslander-mckay} reconstructs $Y$ from $\pi_*\P_Y$ and identifies its contractions with canonical summands, giving the Auslander--McKay embedding. Its final subsection identifies the null category of a contraction with modules over an idempotent quotient and detects its exceptional locus. No type $A$ assumption on the closed fibre is imposed.

\section{Crepant models and modifying modules}\label{sec:general-reflexive-pushforward}

\subsection{Crepancy and vanishing}\label{sec:crepant-vanishing}

Put $N:=\dim R$. Since $R$ is Gorenstein, take $\O_X[N]$ as the normalized dualizing complex of $X$. For every normal integral projective birational morphism $\pi\colon Y\to X$, define its \emph{canonical sheaf} by
\begin{equation}\label{eq:canonical-sheaf-normal-model}
\omega_Y:=\mathcal H^0(\pi^!\O_X).
\end{equation}
Proper duality makes $\pi^!\O_X[N]$ a dualizing complex, and normality makes $\omega_Y$ a rank-one reflexive sheaf \cite[Tags 0AWI and 0AWN]{stacks-project}.

In this setting, we call $\pi$ \emph{crepant} if $\omega_Y\cong\O_Y$.

\begin{Lem}[Crepancy and duality]\label{lem:crepancy-and-duality}
For a normal integral projective birational morphism $\pi\colon Y\to X$,
\begin{equation}\label{eq:crepancy-cm-duality}
\pi^!\O_X\cong\O_Y
\quad\Longleftrightarrow\quad
Y\text{ is Cohen--Macaulay and }\omega_Y\cong\O_Y.
\end{equation}
Under these equivalent conditions, $Y$ is Gorenstein.
\end{Lem}
\begin{proof}
The dualizing complex $\pi^!\O_X[N]$ is concentrated in degree $-N$ precisely when $Y$ is Cohen--Macaulay \cite[Tag 0AWT]{stacks-project}; the degree is fixed by its generic restriction. In that case $\pi^!\O_X\cong\omega_Y$, proving the equivalence. Under these conditions $\O_Y[N]$ is dualizing, hence $Y$ is Gorenstein \cite[Tag 0BFQ]{stacks-project}.
\end{proof}

\begin{Rem}[Canonical divisors]
Fix a generator $\eta$ of $\omega_X=\O_X$. Its orders along prime divisors under the birational identification define canonical Weil divisors $K_X=0$ and $K_Y$, with $\omega_Y\cong\O_Y(K_Y)$. If $\omega_Y\cong\O_Y$, a generator has generic value $h\eta$ for some $h\in(\Frac R)^\times$. Proper birationality makes $\pi$ an isomorphism over every height-one point of $X$ \cite[Tag 0BFP]{stacks-project}, so $h\in R_{\mathfrak p}^\times$ for every height-one prime $\mathfrak p$. Normality gives $h\in R^\times$, and rescaling shows that $\eta$ itself generates $\omega_Y$. Thus our definition of crepancy agrees with $K_{Y/X}=K_Y-\pi^*K_X=0$.
\end{Rem}

We now prove the birational vanishing input by ruling out poles along the closed fibre, then apply duality.

\begin{Lem}[Birational restriction and splitting]\label{lem:birational-split-unit}
Let $(A,\mathfrak n)$ be a noetherian normal local domain and let $g\colon Z\to\Spec A$ be proper and birational, with $Z$ normal integral. Write $F:=g^{-1}(\mathfrak n)$ for the closed fibre. Then:
\begin{enumerate}
\item The restriction map
\begin{equation}\label{eq:birational-h1-injectivity}
H^1(Z,\O_Z)\longrightarrow H^1(Z\setminus F,\O_Z)
\end{equation}
is injective.
\item Assume in addition that $\dim F\le1$. If the unit $A\to\mathbb Rg_*\O_Z$ admits a retraction in $D^b(\mod A)$, then it is an isomorphism.
\end{enumerate}
\end{Lem}
\begin{proof}
(1) If $\dim A\le1$, then $g$ is an isomorphism by \cite[Tag 0BFP]{stacks-project}, so the assertion is immediate. Assume $\dim A\ge2$, and put $W:=Z\setminus F$, $C:=F_{\mathrm{red}}$ and $K:=\Frac A$. Normality gives $g_*\O_Z=A$, so Stein factorization makes $C$ connected \cite[Tag 03H0]{stacks-project}. Since $C$ is proper and reduced, $H^0(C,\O_C)$ is a finite-dimensional reduced connected $A/\mathfrak n$-algebra, hence a field. We will rule out poles by constructing a nonzero element of this field with a zero on $C$.

On a finite affine cover $\{U_\alpha\}$, a class in the kernel of \eqref{eq:birational-h1-injectivity} has a \v{C}ech representative $a_\alpha-a_\beta$, where $a_\alpha\in K$ is regular on $U_\alpha\cap W$ and the differences are regular on overlaps. Suppose poles occur along prime divisors $E_1,\ldots,E_r\subseteq F$. Each $E_j$, being of codimension one in $Z$, is an irreducible component of $C$. Write $\nu_j$ for its valuation and $d_j:=-\nu_j(a_\alpha)>0$ for the pole order, which is well defined since the differences are regular.

Choose any $f\in\mathfrak n\setminus\{0\}$, set $e_j:=\nu_j(f)>0$, and choose $i$ with $d_i/e_i$ maximal. Put $d:=d_i$, $e:=e_i$ and $b_\alpha:=f^d a_\alpha^e$. Then
\[
\nu_j(b_\alpha)=de_j-ed_j\ge0
\]
at every pole divisor meeting $U_\alpha$. There are no other poles, so each $b_\alpha$ is regular by normality. On overlaps, the difference-of-powers identity gives
\[
\nu_j(b_\alpha-b_\beta)
\ge de_j-(e-1)d_j\ge d_j>0.
\]
Thus the restrictions agree on every $E_j$. On every other irreducible component of $C$, each $b_\alpha$ vanishes: at its generic point $a_\alpha$ is regular and $f$ vanishes. Since $C$ is reduced, the restrictions glue to $b\in H^0(C,\O_C)$. Moreover, $\nu_i(b_\alpha)=0$ gives $b\ne0$, so $b$ is invertible.

To obtain a zero, choose a height-one prime $\mathfrak p$ minimal over $(f)$, so $\mathfrak p\ne\mathfrak n$. Since $g$ is proper and birational, it is surjective. Choose $y\in g^{-1}(\mathfrak p)$ and let $H:=\overline{\{y\}}_{\mathrm{red}}\subset Z$. Properness makes $g(H)$ closed; as it contains $\mathfrak p$, it also contains $\mathfrak n$. Hence $H\cap C\ne\varnothing$.

At the generic point $y\in W$ of $H$, the functions $a_\alpha$ are regular and $f$ vanishes. Thus each regular function $b_\alpha$ vanishes on $H\cap U_\alpha$, so $b$ vanishes on $H\cap C$, contradicting its invertibility. No poles occur, so normality makes every $a_\alpha$ regular and the original cocycle is a coboundary.

(2) The bound $\dim F\le1$ and upper semicontinuity make every fibre at most one-dimensional. Normality gives $g_*\O_Z=A$, and \cite[Tag 02V7]{stacks-project} gives $R^ig_*\O_Z=0$ for $i>1$.

Suppose $M:=R^1g_*\O_Z\ne0$. Choose $\mathfrak p$ minimal in its support, and put $B:=A_{\mathfrak p}$ and $N:=M_{\mathfrak p}$. Then $N\ne0$ has finite length, and the base change $g_B\colon Z_B:=Z\times_A B\to\Spec B$ satisfies the same hypotheses. Flat base change and the retraction give
\begin{equation}\label{eq:split-birational-unit}
\mathbb Rg_{B*}\O_{Z_B}\cong B\oplus N[-1],
\end{equation}
with the unit as the first summand. Let $U_B:=\Spec B\setminus\{\mathfrak pB\}$. The restriction map
\[
N=H^1(Z_B,\O_{Z_B})\longrightarrow H^1(g_B^{-1}(U_B),\O_{Z_B})
\]
is zero: the $B$-summand contributes no $H^1$ on the affine base, while $N|_{U_B}=0$. But (1), applied to $g_B$, makes this map injective. This contradiction proves $M=0$, so the unit is an isomorphism.
\end{proof}

\begin{Prop}[Crepant vanishing]\label{prop:crepant-vanishing}
For every model satisfying \eqref{eq:crepant-model-hypotheses}, the natural unit is an isomorphism
\begin{equation}\label{eq:crepant-structure-vanishing}
\O_X\xrightarrow{\ \sim\ }\mathbb R\pi_*\O_Y.
\end{equation}
\end{Prop}
\begin{proof}
By Lemma \ref{lem:crepancy-and-duality} and the comparison of generators above, the canonical generic identification extends to an isomorphism $\O_Y\cong\omega_Y\cong\pi^!\O_X$. Compose its derived pushforward with the proper-duality counit to obtain
\[
t\colon\mathbb R\pi_*\O_Y\longrightarrow R.
\]
The composite with the unit is the identity at the generic point, hence on $R$. Thus $t$ retracts the unit, and Lemma \ref{lem:birational-split-unit}(2) applies.
\end{proof}

This establishes \eqref{eq:geometric-setup} for every object of $\mathsf{Crep}_1(R)$ before using the tilting theory of Part \ref{part:tilting bundles}.

\subsection{Reflexive pushforward}

For a finite module $M$, write $M^*:=\Hom_R(M,R)$. A module $M\in\refl R$ is \emph{modifying} if $\End_R(M)\in\CM R$. It is \emph{maximal modifying} if it is modifying and
\[
M\oplus N\text{ is modifying},\quad N\in\refl R
\quad\Longrightarrow\quad N\in\add M.
\]
The endomorphism algebra of a nonzero modifying module is called a \emph{modifying algebra}; that of a maximal modifying module is called a \emph{maximal modifying algebra} (MMA). A generator means a module with $R\in\add M$. We write $\modif R$ for the set of isomorphism classes of basic modifying modules, including $0$, and $\MM R\subseteq\modif R$ for the maximal modifying ones. Completeness gives Krull--Schmidt decompositions. We use basic representatives throughout; arbitrary modifying modules are obtained by allowing multiplicities of their indecomposable summands.

Let $\mathsf{ModGen}(R)$ be the poset category of isomorphism classes of basic modifying generators, with
\[
\Hom_{\mathsf{ModGen}(R)}(M,N)=
\begin{cases}
\{\ast\},&\add M\subseteq\add N,\\
\varnothing,&\text{otherwise}.
\end{cases}
\]
These notions follow \cite[Section 1]{IW}. For a vector bundle $V$ on a model $\pi\colon Y\to X$, put
\[
\Psi_Y(V):=(\pi_*V)^{**},
\]
and write $\Psi$ when the model is fixed.

\begin{Lem}\label{lem:reflexive-hom-comparison}
For vector bundles $V,W$ on $Y$, reflexive pushforward induces a natural isomorphism
\[
\Hom_Y(V,W)^{**}\cong\Hom_R(\Psi_Y(V),\Psi_Y(W)).
\]
\end{Lem}
\begin{proof}
A proper birational morphism to a normal scheme is an isomorphism over an open subset containing every codimension-one point \cite[Tag 0BFP]{stacks-project}. The two modules in the display are reflexive and agree there. Reflexive extension gives the assertion, compatibly with composition.
\end{proof}

The next lemma uses the duality and reflexive-extension arguments of \cite[Lemma 3.2.9 and Proposition 3.2.10]{VdB04b} and \cite[Lemma 4.8]{IW14}.

\begin{Lem}\label{lem:pretilting-endomorphisms}
Let $U$ be a non-zero pretilting bundle on $Y$.
\begin{enumerate}
\item $\Psi_Y$ gives an equivalence
\[
\add U\xrightarrow{\ \sim\ }\add\Psi_Y(U),
\]
preserving generic ranks.
\item $\Psi_Y(U)$ is modifying.
\item If $\O_Y\in\add U$, then $\pi_*U$ is already reflexive and maximal Cohen--Macaulay, and is a modifying generator.
\end{enumerate}
\end{Lem}
\begin{proof}
Put $B:=\End_Y(U)$. Since $U$ is pretilting, $\mathbb R\pi_*\mathcal End_Y(U)$ is concentrated in degree zero and equals $B$. Proper duality, the hypothesis $\pi^!\O_X\cong\O_Y$ and the trace pairing give
\[
\RHom_R(B,R)
\cong\mathbb R\pi_*\mathcal End_Y(U)^\vee
\cong B.
\]
In particular, $B^*\cong B$ and $\Ext_R^{>0}(B,R)=0$.

(1) A dual module over a normal domain is reflexive, so $B$ is reflexive. For $V,W\in\add U$, the module $\Hom_Y(V,W)$ is a direct summand of a finite direct sum of copies of $B$, hence is reflexive as well. Lemma \ref{lem:reflexive-hom-comparison} therefore identifies the natural map
\[
\Hom_Y(V,W)\xrightarrow{\ \sim\ }\Hom_R(\Psi_Y(V),\Psi_Y(W)).
\]
These identifications respect composition, proving full faithfulness. Splitting idempotents proves essential surjectivity onto $\add\Psi_Y(U)$, and birationality gives preservation of generic ranks.

(2) By (1), $B\cong\End_R(\Psi_Y(U))$. Since $R$ is Gorenstein, the vanishing $\Ext_R^{>0}(B,R)=0$ implies that $B$ is maximal Cohen--Macaulay by local duality. Thus $\Psi_Y(U)$ is modifying.

(3) If $\O_Y\in\add U$, then $\pi_*U=\Hom_Y(\O_Y,U)$ is a direct summand of a finite direct sum of copies of $B$. It is therefore reflexive and maximal Cohen--Macaulay. Hence $\pi_*U\cong\Psi_Y(U)$ is modifying by (2), and $\pi_*\O_Y=R$ belongs to its additive closure.
\end{proof}

\subsection{Reflexive twists and module mutations}

For a rank-one reflexive module $L$, put
\[
\mathsf T_L(N):=(L\otimes_RN)^{**}\qquad(N\in\refl R).
\]
Restriction to codimension one, followed by reflexive extension, shows that $\mathsf T_L$ is an $R$-linear equivalence of $\refl R$ with inverse $\mathsf T_{L^*}$, and that
\begin{equation}\label{eq:reflexive-twist-hom}
\Hom_R(\mathsf T_L(N),\mathsf T_L(N'))\cong\Hom_R(N,N').
\end{equation}
These isomorphisms respect composition. In particular, reflexive twisting preserves modifying modules, maximal modifying modules and basicness.

We use the Iyama--Wemyss mutation construction \cite[Section 6.1]{IW}, with the subscript denoting the summand being replaced, as in \cite[Section 4.1]{IW18}. For a basic modifying module $M=X\oplus N$ with $N\ne0$, choose a minimal right $(\add N)$-approximation $r_X$ and a minimal right $(\add N^*)$-approximation $r_{X^*}$:
\[
r_X\colon N_X\longrightarrow X,\qquad
r_{X^*}\colon N_{X^*}\longrightarrow X^*,
\quad N_X\in\add N,\quad N_{X^*}\in\add N^*.
\]
The right and left mutations are the basic representatives of
\begin{equation}\label{eq:iw-mutation-convention}
\mu_X^+(M):=N\oplus\Ker r_X,\qquad
\mu_X^-(M):=N\oplus(\Ker r_{X^*})^*.
\end{equation}
These definitions use kernels of the approximations, which need not be surjective.

\begin{Lem}\label{lem:twist-preserves-mutation}
For a rank-one reflexive module $L$, the equivalence $\mathsf T_L$ of \eqref{eq:reflexive-twist-hom} satisfies
\[
\mathsf T_L(\mu_X^\pm(M))
\cong\mu_{\mathsf T_L(X)}^\pm(\mathsf T_L(M)).
\]
\end{Lem}
\begin{proof}
An additive equivalence preserves minimal approximations. The kernel of a morphism between reflexive modules is reflexive: its image is torsionfree, and the depth lemma applied to the kernel--image sequence gives the $S_2$ condition for the kernel. Thus the ordinary module kernel is also the categorical kernel in $\refl R$, and $\mathsf T_L$ preserves it. This proves the assertion for right mutation. For left mutation, use
\[
(\mathsf T_L(A))^*\cong\mathsf T_{L^*}(A^*)
\]
and the same assertion for $\mathsf T_{L^*}$ on the dual modules. These identifications follow by restriction to codimension one and reflexive extension.
\end{proof}

\subsection{Compatibility with modifying mutation}\label{sec:oriented-modifying-mutation}

The following comparison will identify the mutation of a tilting bundle with modifying mutation after reflexive pushforward.

\begin{Lem}[Oriented pushforward comparison]\label{lem:oriented-pushforward-mutation}
Suppose that $U\ne0$ and $E_\pm=U\oplus V_\pm$ are basic pretilting bundles. Let
\[
V_+\xrightarrow{i}U_0\xrightarrow{r}V_-,
\qquad U_0\in\add U,
\]
be a complex which is short exact at the generic point, where $i,r$ are respectively minimal left and right $(\add U)$-approximations. Put $N:=\Psi_Y(U)$, $X_\pm:=\Psi_Y(V_\pm)$ and $M_\pm:=N\oplus X_\pm$. Then
\[
\mu_{X_-}^+(M_-)\cong M_+,\qquad
\mu_{X_+}^-(M_+)\cong M_-.
\]
Under the approximation hypotheses, generic exactness is equivalent to $ri=0$ and $\rk U_0=\rk V_++\rk V_-$. In particular, the lemma applies to the exchange sequence of Theorem \ref{thm:facet-mutation}.
\end{Lem}
\begin{proof}
Lemma \ref{lem:pretilting-endomorphisms}(1) preserves the approximations on $\add E_\pm$. At a height-one prime of $R$, the morphism $\pi$ is an isomorphism and $U$ becomes a nonzero free module over a discrete valuation ring. Thus $r$ is split surjective and $i$ is split injective there. Their zero composite and the rank equality make the complex short exact in codimension one.

The kernel of $\Psi_Y(r)$ is reflexive by the kernel argument in Lemma \ref{lem:twist-preserves-mutation}. It agrees with $X_+$ in codimension one, so reflexive extension identifies it with $X_+$. This gives the right mutation formula. Similarly, $\Psi_Y(i)^*$ is a minimal right $(\add N^*)$-approximation with reflexive kernel $X_-^*$, proving the left mutation formula.
\end{proof}

\subsection{Complete fans and maximal modifying modules}\label{sec:complete-fan-realization}

In dimension three, completeness of the tilting fan combines with Iyama--Wemyss's reflexive tilting correspondence to realize modifying modules on a fixed model.

\begin{Prop}[Realization from a complete fan]\label{prop:complete-fan-realization}
Assume $\dim R=3$ and put $n:=\operatorname{rank}_{\mathbb Z}K_0(\per Y)$. Suppose that $R$ admits an MM module with $n$ nonisomorphic indecomposable summands and that the tilting fan of $Y$ is complete in the positive-rank region. Then reflexive pushforward gives surjections
\begin{equation}\label{eq:complete-fan-pushforward}
\Psi_Y\colon\ptilt Y\twoheadrightarrow\modif R,
\qquad
\Psi_Y\colon\tilt Y\twoheadrightarrow\MM R.
\end{equation}
It preserves the number of indecomposable summands and satisfies $\Psi_Y^{-1}(\MM R)=\tilt Y$. If $\pi$ is small, both maps are bijections.
\end{Prop}
\begin{proof}
Put $M:=\pi_*\P_Y$. Lemma \ref{lem:pretilting-endomorphisms} makes $M$ modifying with $|M|=n$. For an MM module $N_0$ as in the hypothesis, $\Hom_R(N_0,M)$ is reflexive partial tilting over $\End_R(N_0)$ \cite[Theorem 4.17]{IW}. Its $n$ summands already exhaust a Bongartz completion, so it is tilting and $M$ is MM \cite[Lemma 4.13 and Theorem 4.17]{IW}.

Set $\Lambda:=\End_R(M)\cong\End_Y(\P_Y)$, $F:=\Hom_R(M,-)$ and $\Phi:=\RHom_Y(\P_Y,-)$. For any nonzero basic modifying module $N$, $F(N)$ is reflexive partial tilting of projective dimension at most one \cite[Theorem 4.17]{IW}. Thus $Q:=\Phi^{-1}F(N)$ is perfect of positive rank $\rk_R N$. Completeness puts $[Q]$ in some $C(E)$, with $E=\bigoplus_jE_j\in\tilt Y$. Since the $[E_j]$ form an integral basis, there are integers $b_j\ge0$ such that
\begin{equation}\label{eq:complete-fan-realization-class}
U:=\bigoplus_jE_j^{\oplus b_j},\qquad [U]=[Q].
\end{equation}

Since $N':=\Psi_Y(U)$ is modifying by Lemma \ref{lem:pretilting-endomorphisms}, the complex $Q':=\Phi^{-1}F(N')$ is perfect. Pushforward followed by $\pi_*U\to N'$ defines the natural $\Lambda$-linear comparison $\alpha_U\colon\Hom_Y(\P_Y,U)\to F(N')$.

For $\mathfrak p\ne\m$, reflexivity makes $F(N)_{\mathfrak p}$ Cohen--Macaulay over $R_{\mathfrak p}$, whose dimension is at most two. Auslander--Buchsbaum for the symmetric order $\Lambda_{\mathfrak p}$ makes this module projective \cite[Lemmas 2.10 and 2.16]{IW}. After base change to $A:=\widehat{R_{\mathfrak p}}$, one therefore has $Q_A\in\add((\P_Y)_A)$. The base-changed morphism satisfies \eqref{eq:geometric-setup}, so Corollary \ref{cor:class-determines-pretilting} gives $U_A\cong Q_A\in\add((\P_Y)_A)$.
Flat base change \cite[Tag 02KH]{stacks-project} and the commutation of finite-module duals with flat tensor products identify $(\alpha_U)_A$ with the comparison on $\add((\P_Y)_A)$. This is an isomorphism by the base change of Lemma \ref{lem:pretilting-endomorphisms}(1). The same membership gives $\Ext_Y^{>0}(\P_Y,U)\otimes_R A=0$. Faithful flatness of $R_{\mathfrak p}\to A$ therefore makes $\alpha_U$ an isomorphism away from $\m$ and shows that all positive cohomology modules of $\Phi(U)$ have finite length.

Since $\Phi(U)$ has no negative cohomology, the natural maps
\[
\Phi(U)\longleftarrow H^0\Phi(U)
=\Hom_Y(\P_Y,U)\xrightarrow{\alpha_U}F(N')
\]
have cones in $D^b_{\fl}(\mod\Lambda)$. By Lemma \ref{lem:tilting-derived-realization}, their inverse images have cones supported on $C$, so both maps become isomorphisms on $Y\setminus C$. Thus
\[
Q'|_{Y\setminus C}\cong U|_{Y\setminus C}.
\]
Normality and $\operatorname{codim}_Y C\ge2$ identify their ranks and determinants globally. Corollary \ref{cor:rank-determinant} gives $[Q']=[U]=[Q]$.

The two-term presilting complexes $F(N')$ and $F(N)$ have the same class in $K_0(\proj\Lambda)$. Their reductions modulo $\m$ are presilting by \cite[Proposition 4.2(a)]{Kim24} and are isomorphic by $g$-vector injectivity over arbitrary fields \cite[Theorem 6.5(a)]{DIJ19}, in the presilting formulation of \cite[Proposition 3.9(1)]{Asa21}. The same reduction criterion makes their direct sum presilting. In a common Bongartz completion the summand classes are linearly independent \cite[Theorem 2.27]{AI}, so $F(N')\cong F(N)$. Reflexive equivalence gives $N'\cong N$, and Lemma \ref{lem:pretilting-endomorphisms}(1) makes $U$ basic. Also $\Psi_Y(0)=0$, proving the first surjection.

By Bongartz completion, basic modifying modules are MM precisely when they have $n$ summands. Together with Lemma \ref{lem:pretilting-endomorphisms} and Corollary \ref{cor:number-of-summands}, this proves the remaining summand-count assertions and the second surjection. If $\pi$ is small, $\Pic(Y)\to\Cl(R)$ is injective. Equal reflexive pushforwards then give equal ranks and determinants, so Corollary \ref{cor:class-determines-pretilting} proves injectivity on $\ptilt Y$.
\end{proof}

\begin{Cor}[Realization on cDV minimal models]\label{cor:cdv-modifying-realization}
Let $R$ be a complete local cDV $k$-algebra with algebraically closed residue field of characteristic zero, and let $\pi\colon Y\to\Spec R$ be a projective $\mathbb Q$-factorial terminalization. With $(\Delta,J)$ as in Corollary \ref{cor:duval-cdv-tilting}, reflexive pushforward gives
\begin{equation}\label{eq:cdv-chamber-realization}
\tilt Y\cong\operatorname{Cham}(\widetilde\Delta,J)
\overset{\Psi_Y}{\twoheadrightarrow}\MM R.
\end{equation}
The map $\Psi_Y\colon\ptilt Y\twoheadrightarrow\modif R$ is also surjective. If $R$ is isolated, both maps are bijections.
\end{Cor}
\begin{proof}
The model is crepant and Gorenstein, and its fibres have dimension at most one by the cDV geometry used in Corollary \ref{cor:duval-cdv-tilting}. Lemma \ref{lem:crepancy-and-duality} therefore puts it in $\mathsf{Crep}_1(R)$. Its terminal singularities are isolated hypersurface singularities \cite[Remark 4.17]{IW14}, so $\pi_*\P_Y$ is MM by \cite[Corollary 4.19]{IW14}. Corollary \ref{cor:duval-cdv-tilting} gives the chamber identification and completeness, so Proposition \ref{prop:complete-fan-realization} applies. If $R$ is isolated, it is terminal \cite[Theorem 1.1]{Rei83}; crepancy then excludes exceptional divisors, so $\pi$ is small.
\end{proof}

Over $\mathbb C$, taking $M=\pi_*\P_Y$ in Iyama--Wemyss's construction gives an index injection $\MM R\hookrightarrow\operatorname{Cham}(\widetilde\Delta,J)$, bijective when $R$ is isolated \cite[Theorem 0.13 and Corollary 9.8]{IWTC}. The map \eqref{eq:cdv-chamber-realization} instead sends each tilting chamber to its reflexive pushforward. In the non-isolated case it can have nontrivial fibres, so these constructions need not be inverse; Section \ref{sec:pushforward-fibres} describes the fibres explicitly for $uv=f$.

\section{The Auslander--McKay embedding and contractions}\label{sec:general-auslander-mckay}

Module blowups reconstruct the models and realize every summand of a canonical generator containing the structure sheaf as a contraction. We first prove these assertions, then deduce the curve--summand correspondence. The final subsection describes the null categories and exceptional loci of these contractions.

\subsection{Module blowups and intermediate models}

Let $M$ be a nonzero finite torsionfree $R$-module of rank $r$, and put $K:=\Frac R$. A generic trivialization defines the determinant image
\begin{equation}\label{eq:module-determinant-ideal}
\mathfrak d(M):=\im\bigl(\textstyle\bigwedge_R^rM
\longrightarrow\bigwedge_K^r(M\otimes_RK)\cong K\bigr).
\end{equation}
It is a fractional ideal, defined up to multiplication by an element of $K^\times$. After scaling it to an ordinary ideal, set
\[
\operatorname{Bl}_M X:=\operatorname{Bl}_{\mathfrak d(M)}X.
\]
The resulting scheme is independent of this choice. For an integral birational $X$-scheme $t\colon W\to X$, the \emph{strict transform} of $M$ is $t^*M$ modulo its torsion subsheaf. Blowing up $M$ is universal among such schemes on which this strict transform is locally free. On $q\colon\operatorname{Bl}_M X\to X$, the strict transform $V$ is a vector bundle with $\det V\cong\O_{\operatorname{Bl}_M X}(1)$. This is Villamayor's module blowup \cite{Vil06}; see \cite[Theorem 2.10]{BCZ12} for its determinant description and universal property, and \cite[Section 2.2]{GI18} for its relation to flattening blowing up.

\begin{Lem}[Intermediate models]\label{lem:intermediate-contraction-vanishing}
Let $\pi\colon Y\to X$ belong to $\mathsf{Crep}_1(R)$, and let $q\colon Z\to X$ be projective and birational, with $Z$ normal integral. If $\pi=q\circ h$, then $Z\in\mathsf{Crep}_1(R)$ and $\mathbb Rh_*\O_Y\cong\O_Z$. Moreover, the reduced closed locus
\[
Z_h:=\{z\in Z\mid\dim h^{-1}(z)=1\}
\]
is finite over $X$.
\end{Lem}
\begin{proof}
The graph makes $h$ projective, and birationality makes it surjective. Thus both $h$ and $q$ have fibres of dimension at most one, and normality gives $h_*\O_Y=\O_Z$. A positive-dimensional fibre of $Z_h\to X$ would give a fibre of $\pi$ of dimension at least two. Hence $Z_h\to X$ is proper and quasi-finite, so finite.

Put $F:=R^1h_*\O_Y$, whose scheme-theoretic support is finite over $X$. Applying $\mathbb Rq_*$ to the truncation triangle for $\mathbb Rh_*\O_Y$ gives
\[
0=R^1\pi_*\O_Y\longrightarrow q_*F\longrightarrow R^2q_*\O_Z=0,
\]
by Proposition \ref{prop:crepant-vanishing} for $\pi$ and the fibre bound for $q$. Faithfulness of finite pushforward on the support gives $F=0$, hence $\mathbb Rh_*\O_Y\cong\O_Z$. Proper duality now yields
\[
q^!\O_X
\cong\mathbb R\mathcal Hom_Z(\mathbb Rh_*\O_Y,q^!\O_X)
\cong\mathbb Rh_*\pi^!\O_X
\cong\O_Z,
\]
so $Z\in\mathsf{Crep}_1(R)$.
\end{proof}

\begin{Prop}[Module blowups and reconstruction]\label{cor:module-reconstruction}
Let $\pi\colon Y\to X$ belong to $\mathsf{Crep}_1(R)$, and let $E$ be a positive-rank globally generated vector bundle. Put $M:=\pi_*E$ and $Z:=\operatorname{Bl}_M X$, with structure morphism $q$ and strict transform $V$ of $M$.
\begin{enumerate}
\item
One has $Z\in\mathsf{Crep}_1(R)$, and $\pi$ factors uniquely as $q\circ h$, with $h^*V\cong E$.
\item
If $\det E$ is $\pi$-ample, then $h$ is an isomorphism. In particular,
\begin{equation}\label{eq:module-reconstruction}
Y\cong_X\operatorname{Bl}_{\pi_*E}X.
\end{equation}
\end{enumerate}
\end{Prop}
\begin{proof}
(1) Proposition \ref{prop:crepant-vanishing} and Lemma \ref{lem:cohomological-dimension-one} give $\Gamma(Y,\O_Y)=R$, $H^1(Y,\O_Y)=0$ and right exactness of $H^1(Y,-)$.

For a globally generated bundle $W$ of rank $s>0$, choose generators of $\Gamma(Y,W)$ and their evaluation sequence
\[
0\longrightarrow\mathcal K\longrightarrow\mathcal F:=\O_Y^{\oplus a}
\longrightarrow W\longrightarrow0.
\]
Then $H^1(Y,\mathcal K)=0$. The kernel of $\bigwedge^s\mathcal F\twoheadrightarrow\det W$ is a quotient of $\mathcal K\otimes\bigwedge^{s-1}\mathcal F$ \cite[Tag 00DO]{stacks-project}, so its $H^1$ also vanishes. Taking sections therefore gives a surjection
\[
\bigwedge_R^s\Gamma(Y,W)\twoheadrightarrow\Gamma(Y,\det W).
\]
Put $L:=\det E$ and $I:=\mathfrak d(M)$, using compatible generic trivializations. Applying this to $W=E^{\oplus n}$ gives $\Gamma(Y,L^{\otimes n})=\mathfrak d(M^{\oplus n})=I^n$, hence
\begin{equation}\label{eq:determinant-section-ring}
\bigoplus_{n\ge0}\Gamma(Y,L^{\otimes n})
=\bigoplus_{n\ge0}I^n.
\end{equation}
The left-hand ring is the ring of global functions on the total space of $L^\vee$. It is normal because $Y$ is normal \cite[Tag 0358]{stacks-project}; thus $Z=\operatorname{Proj}_R\bigoplus_{n\ge0}I^n$ is normal.

Evaluation identifies $E$ with the strict transform of $M$ on $Y$, so universality gives $h\colon Y\to Z$. The induced surjection $h^*V\twoheadrightarrow E$ is an isomorphism between bundles of equal rank. Lemma \ref{lem:intermediate-contraction-vanishing} puts $Z$ in $\mathsf{Crep}_1(R)$.

(2) We have $h^*\O_Z(1)\cong\det E$. If $\det E$ is $\pi$-ample, its restriction to each $h$-fibre is both ample and trivial, so every $h$-fibre is zero-dimensional. Thus $h$ is finite and birational with normal target, hence an isomorphism.
\end{proof}

Proposition \ref{prop:canonical-tilting-general} and Lemma \ref{lem:pretilting-endomorphisms}(1), (3) give a basic modifying generator
\begin{equation}\label{eq:canonical-modifying-generator}
M_Y:=\pi_*\P_Y=\Psi_Y(\P_Y),
\qquad Y\cong_X\operatorname{Bl}_{M_Y}X,
\end{equation}
where reconstruction follows by applying Proposition \ref{cor:module-reconstruction}(2) to $\P_Y$. Compare \cite[Corollary 4.4]{Wem18} in the three-dimensional crepant setting.

\subsection{Canonical summands and contractions}

For compatibility of canonical generators with contractions, compare \cite[Proposition 2.14]{Wem18}, where the dual generator is used.

\begin{Thm}[Auslander--McKay embedding and contractions]\label{thm:general-auslander-mckay-embedding}
Let $\pi\colon Y\to X$ belong to $\mathsf{Crep}_1(R)$.
\begin{enumerate}
\item \emph{Auslander--McKay embedding.}
The canonical generators define a fully faithful functor
\begin{equation}\label{eq:general-auslander-mckay-embedding}
\mathsf{Crep}_1(R)^{\op}\lhook\joinrel\longrightarrow\mathsf{ModGen}(R),
\qquad Y\longmapsto M_Y=\pi_*\P_Y.
\end{equation}
Explicitly, for any two such models $Y,Z$,
\begin{equation}\label{eq:general-model-morphisms}
\Hom_X(Y,Z)\ne\varnothing
\quad\Longleftrightarrow\quad
\add M_Z\subseteq\add M_Y,
\end{equation}
and this Hom set is a singleton whenever it is nonempty.

\item \emph{Realization of summands.}
Let $U$ be a direct summand of $\P_Y$ containing $\O_Y$, and put $N:=\pi_*U$. Then $Z:=\operatorname{Bl}_N X$ belongs to $\mathsf{Crep}_1(R)$, and $\pi$ factors through $h\colon Y\to Z$, with
\[
h^*\P_Z\cong U,\qquad M_Z\cong N.
\]
\end{enumerate}
\end{Thm}
\begin{proof}
For any $h\colon Y\to Z$ in $\mathsf{Crep}_1(R)$, put $U:=h^*\P_Z$. This bundle is globally generated, and Lemma \ref{lem:intermediate-contraction-vanishing} and the projection formula give
\[
H^1(Y,U^\vee)=H^1(Z,\P_Z^\vee)=0,
\qquad \pi_*U=M_Z.
\]
Thus $U\in\add\P_Y$ by \eqref{eq:canonical-projective-bundles}. Fully faithful pullback makes $U$ basic, so it is a direct summand containing $\O_Y$.

Conversely, let $U$ be such a summand and put $N:=\pi_*U$. Proposition \ref{cor:module-reconstruction}(1) constructs $h\colon Y\to Z:=\operatorname{Bl}_N X$ in $\mathsf{Crep}_1(R)$, with $h^*V\cong U$ for the strict transform $V$ of $N$. The projection formula gives $h_*U\cong V$, so $V$ contains $\O_Z$; it is basic because $h^*$ is fully faithful. Moreover, $V$ is globally generated, $\det V$ is ample over $X$, and
\[
H^1(Z,V^\vee)=H^1(Y,U^\vee)=0.
\]
Proposition \ref{prop:canonical-tilting-general} gives $V\cong\P_Z$ and $M_Z\cong N$, proving (2).

Lemma \ref{lem:pretilting-endomorphisms}(1) identifies the summands of $\P_Y$ containing $\O_Y$ with those of $M_Y$ containing $R$. Hence $\add M_Z\subseteq\add M_Y$ gives an $X$-morphism $Y\to\operatorname{Bl}_{M_Z}X\cong_X Z$ by the construction above and \eqref{eq:canonical-modifying-generator}. The converse follows from the first paragraph. Two $X$-morphisms $Y\to Z$ agree generically and therefore everywhere, since $Y$ is integral and $Z$ is separated over $X$. This proves (1).
\end{proof}

Consequently, for fixed $Y$, contractions $h\colon Y\to Z$ in $\mathsf{Crep}_1(R)$, up to $X$-isomorphism of their targets, correspond bijectively to direct summands $N$ of $M_Y$ containing $R$, taken up to isomorphism, via
\begin{equation}\label{eq:canonical-summand-contractions}
(h\colon Y\to Z)\longleftrightarrow N=M_Z.
\end{equation}

Van den Bergh's curve--summand correspondence \cite[Sections 3.4--3.5, especially Proposition 3.5.4 and Theorem 3.5.5]{VdB04b} now follows from the contractions just constructed.

\begin{Lem}[Curves and canonical summands]\label{lem:canonical-curve-summands}
Let $\pi\colon Y\to X$ belong to $\mathsf{Crep}_1(R)$, and let $C_1,\ldots,C_s$ be the one-dimensional irreducible components of its closed fibre. The canonical decompositions
\[
\P_Y=\O_Y\oplus\bigoplus_{i=1}^s P_i,\qquad
M_Y=R\oplus\bigoplus_{i=1}^s M_i,\qquad M_i:=\pi_*P_i,
\]
can be indexed uniquely by these curves so that
\begin{equation}\label{eq:canonical-curve-degrees}
d_j(P_i)=0\quad(j\ne i),\qquad d_i(P_i)>0.
\end{equation}
Thus $C_i\longleftrightarrow M_i$ is a bijection between the closed-fibre curves and the nonfree indecomposable summands of $M_Y$.
\end{Lem}
\begin{proof}
Lemma \ref{lem:euler-dual-bases}, Corollary \ref{cor:rank-determinant} and Lemma \ref{lem:picard-lattice} give $|\P_Y|=s+1$. For each nontrivial summand $P_i$, Theorem \ref{thm:general-auslander-mckay-embedding}(2) gives $h_i\colon Y\to Z_i$ with
\[
h_i^*\P_{Z_i}\cong\O_Y\oplus P_i.
\]
The same count shows that the closed fibre of $Z_i$ has exactly one curve $D_i$. Its generic point is outside the positive-dimensional fibre locus of $h_i$, which is finite over $X$. Normality therefore makes $h_i$ an isomorphism there, so exactly one $C_j$ dominates $D_i$. Thus $\det P_i\cong h_i^*\det\P_{Z_i}$ has positive degree on that $C_j$ and degree zero on all other curves.

Since $\det\P_Y$ is $\pi$-ample, every $C_j$ occurs this way. Equal counts give the unique indexing in \eqref{eq:canonical-curve-degrees}. Lemma \ref{lem:pretilting-endomorphisms}(1), (3) transfers the correspondence to the modules $M_i$.
\end{proof}

Consequently, the contraction associated with $U=\O_Y\oplus\bigoplus_{i\in I}P_i$ contracts precisely the curves $C_j$ with $j\notin I$: its pulled-back ample determinant has degree zero exactly on these curves. Thus every subset of the closed-fibre curves has a unique contraction in $\mathsf{Crep}_1(R)$, up to $X$-isomorphism of the target.

Part \ref{part:MM} proves essential surjectivity of the embedding for $uv=f$ by classifying all modifying generators. We next describe the null category of each contraction.

\subsection{Contraction algebras and null categories}\label{sec:general-contraction-algebras}

The Auslander--McKay embedding identifies contractions with additive inclusions of canonical generators. The corresponding idempotent quotient describes the sheaves annihilated by the contraction. Its support detects where the contraction fails to be an isomorphism, as in Donovan--Wemyss \cite[Theorems 1.2 and 1.4]{DW19contractions}. The categorical description follows from Theorem \ref{thm:semistable-reduction}; the support statement is proved by localization.

Let $h\colon Y\to Z$ be an $X$-morphism in $\mathsf{Crep}_1(R)$, and write $\pi=q\circ h$. Put
\[
P:=\P_Y,\qquad U:=h^*\P_Z,\qquad \Lambda:=\End_Y(P).
\]
Theorem \ref{thm:general-auslander-mckay-embedding} identifies $U$ with a summand of $P$ containing $\O_Y$, and every such summand arises from a unique contraction. Choose its summand idempotent $e_U\in\Lambda$, and define
\begin{equation}\label{eq:general-contraction-algebra}
A_h:=\Lambda/\Lambda e_U\Lambda,\qquad
\mathcal N_h:=\{F\in\Coh Y\mid\mathbb Rh_*F=0\}.
\end{equation}

\begin{Thm}[Contractions and idempotent quotients]\label{thm:contraction-null-category}
With the notation above, one has the following.
\begin{enumerate}
\item The functor $\Ext_Y^1(P,-)$ gives an exact equivalence
\begin{equation}\label{eq:contraction-null-equivalence}
\mathcal N_h\xrightarrow{\ \sim\ }\mod A_h.
\end{equation}
It restricts to $\mathcal N_h\cap\Coh_0Y\simeq\fl A_h$. This finite-length category is $\W_U$ from Theorem \ref{thm:semistable-reduction}, and has $|\P_Y|-|\P_Z|$ simple objects.
\item For the reduced closed locus
\[
Z_h:=\{z\in Z\mid\dim h^{-1}(z)=1\}
\]
one has an equality of closed sets
\begin{equation}\label{eq:contraction-algebra-support}
\operatorname{Supp}_R A_h=q(Z_h).
\end{equation}
In particular, $A_h=0$ if and only if $h$ is an isomorphism.
\end{enumerate}
\end{Thm}
\begin{proof}
(1) Both $U$ and $P$ are globally generated and contain $\O_Y$ in their additive closures. Hence
\[
\Fac U=\Fac\O_Y=\Fac P.
\]
Thus $P\cong\E_{\min}(U)$ by \eqref{eq:extremal-completions} and the uniqueness in Theorem \ref{thm:tilting-torsion}.
Derived adjunction and generation by $\P_Z$ give, for $F\in\Coh Y$,
\[
\mathbb Rh_*F=0
\quad\Longleftrightarrow\quad
\RHom_Y(U,F)=0
\quad\Longleftrightarrow\quad
F\in U^{\perp_0}\cap U^{\perp_1}.
\]
The last equivalence uses cohomological dimension one. Hence $\mathcal N_h=\widetilde{\W}_U$, and Theorem \ref{thm:semistable-reduction} gives all assertions in (1), with $B_{\min}=A_h$ and $|U|=|\P_Z|$.

(2) Localize at $\mathfrak p\in\Spec R$. Tilting equivalences and the inclusion $U\in\add P$ localize. If $h_{\mathfrak p}$ is an isomorphism, then $U_{\mathfrak p}$ is a tilting generator in $\add P_{\mathfrak p}$, so its idempotent is full and $(A_h)_{\mathfrak p}=0$.

Conversely, if $(A_h)_{\mathfrak p}=0$, fullness of the idempotent gives $P_{\mathfrak p}\in\add U_{\mathfrak p}$. The idempotent defining this summand descends to $Z_{\mathfrak p}$, since $h_{\mathfrak p*}\O_{Y_{\mathfrak p}}=\O_{Z_{\mathfrak p}}$ makes pullback fully faithful on vector bundles. Hence $P_{\mathfrak p}\cong h_{\mathfrak p}^*V$ for a vector bundle $V$ on $Z_{\mathfrak p}$. The line bundle $\det P_{\mathfrak p}$ is ample over $\Spec R_{\mathfrak p}$ and trivial on every $h_{\mathfrak p}$-fibre. Thus $h_{\mathfrak p}$ has no positive-dimensional fibres. It is finite and birational with normal target, hence is an isomorphism.

We have proved $(A_h)_{\mathfrak p}=0$ exactly when $h_{\mathfrak p}$ is an isomorphism. The non-isomorphism locus in $Z$ is $Z_h$, which is finite over $X$ by Lemma \ref{lem:intermediate-contraction-vanishing}. Its image is therefore closed, and the local criterion proves \eqref{eq:contraction-algebra-support}.
\end{proof}

\begin{Cor}[Detecting exceptional divisors]\label{cor:contraction-algebra-dimension}
If $h$ is not an isomorphism, then
\begin{equation}\label{eq:contraction-algebra-dimension}
\dim\operatorname{Exc}(h)=\dim_R A_h+1.
\end{equation}
Consequently,
\[
\begin{aligned}
h\text{ has an exceptional divisor}
&\quad\Longleftrightarrow\quad \dim_R A_h=\dim R-2,\\
h\text{ is small}
&\quad\Longleftrightarrow\quad \dim_R A_h\le\dim R-3.
\end{aligned}
\]
Here small means that the exceptional set has codimension at least two, including the isomorphism case. If $\dim R=3$, then $h$ is small if and only if $\ell_R(A_h)<\infty$.
\end{Cor}
\begin{proof}
The map $Z_h\to X$ is finite, and the fibres of $h^{-1}(Z_h)\to Z_h$ have dimension one. Thus Theorem \ref{thm:contraction-null-category}(2) gives
\[
\dim\operatorname{Exc}(h)=\dim Z_h+1=\dim_R A_h+1.
\]
Since $\dim Y=\dim R$, the codimension statements follow. A finite $R$-module has support of dimension at most zero exactly when it has finite $R$-length, giving the last assertion.
\end{proof}

The category in \eqref{eq:contraction-null-equivalence} is the abelian category of coherent sheaves annihilated by $\mathbb Rh_*$. The equivalence does not assert that the full triangulated kernel of $\mathbb Rh_*$ is the bounded derived category of the ordinary quotient $A_h$. Part \ref{part:MM} computes these quotients and their wall categories explicitly for flag models.

\part{Modifying modules and flag models for \texorpdfstring{$uv=f$}{uv=f}}\label{part:MM}

Throughout Part \ref{part:MM}, let $(S,\n,k)$ be a regular complete local ring with $d:=\dim S\ge1$ and take $0\neq f\in\n$. Put $R:=S[[u,v]]/(f-uv)$ and $\m:=(\n,u,v)R$. Then $(R,\m,k)$ is a $(d+1)$-dimensional Gorenstein complete local ring, and Corollary \ref{cor:uvf-normality} shows that it is a normal domain. The notions and functors of Part \ref{part:auslander-mckay} will be used throughout. The flag models constructed below satisfy \eqref{eq:crepant-model-hypotheses}, while their explicit closed fibres also bring the type $A$ results of Part \ref{part:tilting bundles} into play.

Fix a prime factorization, allowing associate factors,
\[
f=f_1\cdots f_n=\varepsilon\prod_{\alpha=1}^{\ell}p_\alpha^{a_\alpha},
\qquad \varepsilon\in S^\times,\quad a_\alpha\ge1,
\]
where the $p_\alpha$ are pairwise nonassociate. Write $a=(a_1,\ldots,a_\ell)$, so $n=\sum_\alpha a_\alpha$, and write $s\sim t$ when $s,t\in S$ are associate.

Section \ref{sec:modifying-modules} classifies modifying modules and computes their mutations, proving connectivity of the MM mutation graph. Section \ref{sec:flag-models} constructs the flag models and canonical tilting generators, proves the Auslander--McKay correspondence in relative dimension one, and determines when an NCCR exists. Section \ref{sec:reconstruction-pushforward} uses mutation compatibility to realize all modifying modules on one maximal-flag model. It then describes the fibres, quotient mutation graph and wall categories of this realization.

\section{Classification and mutation of modifying modules}\label{sec:modifying-modules}

The definitions of modifying modules, maximal modifying algebras and reflexive twists are those of Part \ref{part:auslander-mckay}. We now classify these modules using the prime factors of $f$.

\subsection{Singular locus}

\begin{Lem}\label{uvreg}
The open subscheme $D(u)\cup D(v)\subseteq\Spec R$ is regular.
\end{Lem}
\begin{proof}
By symmetry it suffices to consider $P\in D(u)\subseteq\Spec R$. Let $Q$ be its inverse image in $A:=S[[u,v]]$. If $f-uv\in Q^2A_Q$, applying the derivation $\partial/\partial v$ gives
\[
-u=\frac{\partial}{\partial v}(f-uv)\in QA_Q,
\]
contrary to $u\notin Q$. Thus $R_P\cong A_Q/(f-uv)$ is regular.
\end{proof}

As an immediate corollary, we obtain the following.

\begin{Cor}\label{cor:uvf-normality}
$R$ is a normal domain.
\end{Cor}
\begin{proof}
It is enough to show that $R$ satisfies the $(R_1)$ condition. Take $P\in\Spec R$ with $\height P\le1$. Since $\dim R/(u,v)=\dim S/(f)=d-1=\dim R-2$ holds and $R$ is a catenary ring, we have $P\in D(u)\cup D(v)$. Thus by Lemma \ref{uvreg}, $R_P$ is regular.
\end{proof}

Using Lemma \ref{uvreg}, we can determine the singular locus of $R$.

\begin{Prop}
Consider the surjection $\varpi\colon R\to R/(u,v)\cong S/(f)$. This induces a bijection
\[\varpi^{-1}\colon\Sing S/(f)\xrightarrow[\cong]{}\Sing R.\]
\end{Prop}
\begin{proof}
By Lemma \ref{uvreg}, we have $\Sing R\subseteq V(u,v)\cong\Spec S/(f)$. Let $P\in V(u,v)\subseteq\Spec R$ and take $\p\in V(f)\subseteq\Spec S$ such that $P=\varpi^{-1}(\p/(f))=\p R+(u,v)R$ holds. We prove that $R_P$ is singular if and only if so is $S_\p/(f)$.

Remark that the $(u,v)$-adic completion of $R_P$ is isomorphic to $S_\p[[u,v]]/(f-uv)$. Thus we have
\begin{equation*}
\begin{split}
R_P\text{ is singular}&\Leftrightarrow S_\p[[u,v]]/(f-uv)\text{ is singular}\\
&\Leftrightarrow f-uv\in(\p S_\p[[u,v]]+(u,v))^2\subseteq S_\p[[u,v]]\\
&\Leftrightarrow f\in(\p S_\p)^2\\
&\Leftrightarrow S_\p/(f)\text{ is singular}\qedhere
\end{split}
\end{equation*}
\end{proof}

\subsection{Divisor class groups}\label{sec:divisor-class-groups}

Use the prime factorization fixed at the start of Part \ref{part:MM}. For a divisor $z$ of $f$, write $I(z):=(u,z)$. Put
\[
f_I:=\prod_{j\in I}f_j,\quad T_I:=I(f_I),\quad f_\emptyset:=1,
\qquad \mathfrak P_\alpha:=I(p_\alpha).
\]
We write ideal classes additively: $[L]+[N]=[(L\otimes_RN)^{**}]$ and $[R]=0$. The height-one primes containing $u$ are precisely the $\mathfrak P_\alpha$. Their normalized valuations satisfy
\begin{equation}\label{eq:factor-ideal-valuations}
\nu_\alpha(u)=a_\alpha,\qquad \nu_\alpha(p_\beta)=\delta_{\alpha\beta},
\end{equation}
since $v$ and $f/p_\alpha^{a_\alpha}$ are units in $R_{\mathfrak P_\alpha}$.

\begin{Thm}\label{thm:divisor-class-group}
There is an isomorphism
\begin{equation}\label{eq:divisor-class-group}
\mathbb Z^{\ell}/\mathbb Z(a_1,\ldots,a_{\ell})
\xrightarrow{\ \sim\ }\Cl(R),\qquad \overline{e_\alpha}\longmapsto[\mathfrak P_\alpha].
\end{equation}
For $I\subseteq\{1,\ldots,n\}$, the module $T_I$ is maximal Cohen--Macaulay and rank-one reflexive, with
\begin{equation}\label{eq:factor-ideal-classes}
[T_I]=\sum_\alpha m_\alpha(I)[\mathfrak P_\alpha],\qquad
m_\alpha(I):=\#\{j\in I\mid(f_j)=(p_\alpha)\}.
\end{equation}
Thus $[T_I]+[T_J]=[T_{I\sqcup J}]$ for disjoint $I,J$, and $[T_{\{1,\ldots,n\}}]=0$.
\end{Thm}

This extends \cite[Theorem 5.19]{IW18}. We separate generation from relations: a finite-quotient argument gives the generators, and an elementary colon calculation in $R$ determines all relations.

\begin{Lem}\label{lem:class-group-generation}
We have $\Pic(R[1/u])=0$, and $[\mathfrak P_1],\ldots,[\mathfrak P_{\ell}]$ generate $\Cl(R)$.
\end{Lem}
\begin{proof}
Put $A:=R[1/u]$ and extend $L\in\Pic(A)$ to a rank-one reflexive $R$-module $M$. We use Gabber's finite-quotient argument \cite[proof of Lemma 6.2]{Poi22}. The module $M$ is free at the two generic points of $R/\n R=k[[u,v]]/(uv)$ by Lemma \ref{uvreg}. Trivializing its dual over the corresponding semilocal ring and clearing denominators identifies $M$ with an ideal $I\subseteq R$ which is the unit ideal at both points.

Thus $D:=R/I$ has finite-dimensional $D/\n D$. The $\n$-adic and $\m$-adic topologies on this complete quotient agree. Complete Nakayama makes monomial lifts of a $k$-basis of $D/\n D$ into $S$-generators of $D$. Consequently,
\[
B:=\bigl(S[u,v]/(uv-f)\bigr)_{(\n,u,v)}\longrightarrow D
\]
is surjective. Since $\widehat B=R$ and $D$ is $\m$-adically complete, its kernel $J$ satisfies, compatibly with $R\to D$,
\[
R/JR\cong\widehat{B/J}\cong D=R/I,
\qquad\text{hence }JR=I.
\]
Faithfully flat descent along $B_u\to A$ makes $J_u$ invertible, since $J_u\otimes_{B_u}A=I_u\cong L$. As $B_u$ is a localization of the UFD $S[u,u^{-1}]$, it follows that $J_u\cong B_u$, and hence $L\cong A$.

Finally, $A$ is regular by Lemma \ref{uvreg}, so $\Cl(A)=\Pic(A)=0$. The localization sequence for divisor class groups now gives the generation assertion.
\end{proof}

\begin{proof}[Proof of Theorem \ref{thm:divisor-class-group}]
Every element of $R$ has a unique expansion
\[
x=\sum_{j\ge0}\alpha_ju^j+\sum_{r\ge1}\beta_rv^r,
\qquad\alpha_j,\beta_r\in S.
\]
For existence, choose a representative $\sum_{i,j\ge0}h_{ij}u^iv^j$. Reducing mixed monomials by $uv=f$ gives
\[
\alpha_j=\sum_{q\ge0}h_{j+q,q}f^q,
\qquad
\beta_r=\sum_{q\ge0}h_{q,q+r}f^q.
\]
These sums converge in the complete ring $S$ because $f\in\n$, yielding the claimed expansion. For uniqueness, reduce a zero expansion in $R/(f)\cong(S/(f))[[u,v]]/(uv)$. The pure $u$-terms and positive $v$-terms are independent there, so every coefficient is divisible by $f$. Cancel the non-zero divisor $f$ in $R$ and repeat: every coefficient lies in $\bigcap_{q\ge0}f^qS=0$.

For $0\ne s\in S$ and $N\ge0$, choose $c_jS=(f^{N-j}):_Ss$ for $0\le j\le N$, with $c_N=1$, using that $S$ is a UFD. Coefficient comparison gives $sx\in u^NR$ exactly when $\alpha_j\in c_jS$ for $j<N$ and $\beta_r\in c_0S$. Grouping the $v$-series with $c_0$ and the $u$-tail with $u^N$ yields
\begin{equation}\label{eq:power-colon-generators}
(u^N):_Rs=\sum_{j=0}^N c_ju^jR.
\end{equation}
If this ideal is principal, Nakayama makes one $u^jc_j$ a generator. It contains $u^N\notin\n R$, as $R/\n R\cong k[[u,v]]/(uv)$, so $c_j$ is a unit and the ideal equals $(u^j)$.

Lemma \ref{lem:class-group-generation} gives the surjection $\Phi\colon\mathbb Z^\ell\to\Cl(R)$, $\Phi(e_\alpha)=[\mathfrak P_\alpha]$, and \eqref{eq:factor-ideal-valuations} gives $\Phi(a)=0$. Shift $b\in\ker\Phi$ by a multiple of $a$ to make $b\ge0$, and set $s:=\prod_\alpha p_\alpha^{Na_\alpha-b_\alpha}$ with $Na_\alpha\ge b_\alpha$. The divisorial ideal
\[
J:=(u^N):_Rs=R\cap(u^N/s)R
\]
has valuation $b_\alpha$ at $\mathfrak P_\alpha$ and zero elsewhere, so $[J]=\Phi(b)=0$. It is principal, hence $J=(u^j)$ by the preceding calculation. Valuations give $b=ja$; undoing the shift proves $\ker\Phi=\mathbb Za$ and \eqref{eq:divisor-class-group}.

For $I\ne\emptyset$, the quotient $R/T_I\cong(S/(f_I))[[v]]$ is Cohen--Macaulay of dimension $d$, so the depth lemma makes $T_I$ maximal Cohen--Macaulay and rank-one reflexive; $T_\emptyset=R$. Its valuations are $m_\alpha(I)$ at $\mathfrak P_\alpha$ and zero elsewhere, giving \eqref{eq:factor-ideal-classes} and the product formula.
\end{proof}

\begin{Cor}\label{cor:class-group-factoriality}
For $q:=\gcd(a_1,\ldots,a_{\ell})$, noncanonically,
\[
\Cl(R)\cong\mathbb Z^{\ell-1}\oplus\mathbb Z/q\mathbb Z.
\]
Its torsion subgroup is generated by $[(u,\prod_\alpha p_\alpha^{a_\alpha/q})]$, of order $q$. Moreover, $R$ is factorial exactly when $f$ is prime up to a unit, and $\Spec R$ is $\mathbb Q$-factorial exactly when $\ell=1$.
\end{Cor}
\begin{proof}
The primitive vector $a/q$ generates the torsion of $\mathbb Z^{\ell}/\mathbb Za$. Smith normal form gives the decomposition. Since $R$ is normal and local, factoriality means $\Cl(R)=0$, and $\mathbb Q$-factoriality means that $\Cl(R)$ is torsion.
\end{proof}

\subsection{Indecomposable modifying modules}

We prove that every indecomposable modifying $R$-module has rank one (Theorem \ref{modrk1}) by reducing to a normal surface. The induction step preserves rank, the modifying property and indecomposability through an isomorphism of endomorphism algebras. In the surface case, Kn\"orrer periodicity reduces the problem to the structure theorem for modules over a discrete valuation ring. We first recall the regularity criterion used to choose the parameter.

\begin{Thm}\cite[4.4]{AG}\label{AusGol}
Let $A$ be a regular local ring. Then all modifying $A$-modules are free.
\end{Thm}

The next lemma gives the induction step. Its parameter choice extends \cite[Lemma 5.23]{IW18}, and the endomorphism-algebra argument follows the proof of \cite[Proposition 5.24]{IW18}.

\begin{Lem}[Hyperplane reduction]\label{modred}
Assume $d\ge2$, and let $M$ be a nonzero modifying $R$-module. There exists $t\in\n\setminus\n^2$ with $f\notin tS$ such that, writing
\[
\overline S:=S/(t),\qquad
\overline R:=R/(t)\cong\overline S[[u,v]]/(\overline f-uv),
\qquad L:=(M/tM)^{**}_{\overline R},
\]
the natural map induces an isomorphism of $\overline R$-algebras
\begin{equation}\label{eq:modifying-hyperplane-endomorphisms}
\End_R(M)/t\End_R(M)\xrightarrow{\ \sim\ }\End_{\overline R}(L).
\end{equation}
In particular, $L$ is modifying and $\rk_{\overline R}L=\rk_RM$. If $M$ is indecomposable, then so is $L$.
\end{Lem}
\begin{proof}
Put $E:=\Ext_R^1(M,M)$ and $\La:=\End_R(M)$. Modifying modules remain modifying after localization, so Theorem \ref{AusGol} and Lemma \ref{uvreg} give
\[
\operatorname{Supp}E\subseteq\operatorname{Sing}R\subseteq V(u,v).
\]
Choose $0\to\Omega M\to P\to M\to0$ with $P$ finite free. It induces
\[
0\longrightarrow\La\longrightarrow\Hom_R(P,M)
\longrightarrow\Hom_R(\Omega M,M)\longrightarrow E\longrightarrow0.
\]
The two Hom modules are reflexive, hence have depth at least two, while $\La$ has depth $d+1\ge3$. The depth lemma gives $\depth_RE\ge1$, so $\m\notin\Ass_RE$. Together with the support inclusion, this shows that $\p\cap S\ne\n$ for every $\p\in\Ass_RE$.

The finite set
\[
\mathcal Q:=\{\p\cap S\mid\p\in\Ass_RE\}
\cup\{(f_i)\mid1\le i\le n\}
\]
consists of nonmaximal prime ideals of $S$. Prime avoidance, allowing the one nonprime ideal $\n^2$ \cite[Tag 00DS]{stacks-project}, gives
\[
t\in\n\setminus\left(\n^2\cup\bigcup_{\mathfrak q\in\mathcal Q}\mathfrak q\right).
\]
Thus $t$ acts injectively on $E$. Moreover, $S/(t)$ is regular complete local of dimension $d-1$, so $t$ is prime in $S$. If $t$ divided $f$, it would be associated to some $f_i$, contrary to its choice. Hence $\overline f\ne0$, and $\overline R$ is again a normal domain of the required form.

Put $N:=M/tM$. Since $M$ is torsionfree, $t$ is $M$-regular. Applying $\Hom_R(M,-)$ to $0\to M\xrightarrow{t}M\to N\to0$ and using injectivity on $E$ gives
\[
\La/t\La\cong\Hom_R(M,N)\cong\End_{\overline R}(N).
\]
This is a maximal Cohen--Macaulay, hence reflexive, $\overline R$-module.

We next check that $N$ is torsionfree. For a nonzero prime $\mathfrak q$ of $\overline R$, its inverse image $\p$ in $R$ strictly contains the height-one prime $(t)$, so $\height\p\ge2$. Reflexivity of $M$ gives the $S_2$ condition \cite[Tag 0AVB]{stacks-project}, and therefore
\[
\depth N_{\mathfrak q}=\depth M_{\p}-1\ge1.
\]
Thus $N$ has no nonzero associated prime and is torsionfree. It follows that $N$ and $L=N^{**}_{\overline R}$ agree in codimension one. Double dualization induces
\[
\End_{\overline R}(N)\longrightarrow\End_{\overline R}(L),
\]
which is an isomorphism in codimension one. Both sides are reflexive, so it is an isomorphism everywhere. This proves \eqref{eq:modifying-hyperplane-endomorphisms} and shows that $L$ is modifying.

Finally, $M_{(t)}$ is free over the discrete valuation ring $R_{(t)}$. Reduction modulo $t$ computes the generic rank of $N$ over $\overline R$, giving $\rk_{\overline R}L=\rk_{\overline R}N=\rk_RM$. If $M$ is indecomposable, completeness makes $\La$ local. Its quotient $\End_{\overline R}(L)\cong\La/t\La$ is then local, so $L$ is indecomposable.
\end{proof}

For $S=\mathbb C[[x_1,x_2]]$, the following theorem is \cite[Proposition 5.24]{IW18}.

\begin{Thm}\label{modrk1}
Let $(S,\n,k)$ be a regular complete local ring with $d:=\dim S\ge1$ and take $0\ne f\in\n$. Put $R:=S[[u,v]]/(f-uv)$. Then every indecomposable modifying $R$-module has rank one.
\end{Thm}
\begin{proof}
If $d\ge2$, Lemma \ref{modred} associates to an indecomposable modifying module $M$ an indecomposable modifying module $L$ over $\overline R=\overline S[[u,v]]/(\overline f-uv)$, with $\dim\overline S=d-1$ and $\rk_{\overline R}L=\rk_RM$. Thus the assertion for $d-1$ implies that for $d$. By induction, it remains to prove the case $d=1$.

In this case $S$ is a complete discrete valuation ring, and after rescaling $v$ by a unit we may write
\[
R\cong S[[u,v]]/(\pi^e-uv),
\]
where $\pi$ is a uniformizer and $e\ge1$. If $e=1$, then $R$ is regular and every reflexive $R$-module is free. Assume $e\ge2$. Since $R$ is a normal surface, reflexive modules are maximal Cohen--Macaulay. The indecomposable nonfree modules over the Artinian principal ideal ring $S/(\pi^e)$ are exactly $S/(\pi^i)$ with $1\le i<e$. Hyperbolic Kn\"orrer periodicity of Pfister and Popescu \cite[Theorem 3.12]{PP96}, as presented in \cite[Theorem 2.4 and Corollary 2.6]{Leu26}, identifies them with the indecomposable nonfree maximal Cohen--Macaulay $R$-modules. Their matrix factorizations $(\pi^i,\pi^{e-i})$ give the images
\[
\operatorname{coker}\begin{pmatrix}\pi^i&u\\v&\pi^{e-i}\end{pmatrix}
\cong(u,\pi^{e-i}),
\qquad 1\le i<e,
\]
all of rank one. The free indecomposable module $R$ also has rank one, completing the proof.
\end{proof}

\subsection{Classification by flags}\label{sec:modifying-classification}

We now combine Theorems \ref{thm:divisor-class-group} and \ref{modrk1} with a rank-one consequence of Kn\"orrer periodicity, using the factorization and multiplicity vector fixed at the start of Part \ref{part:MM}. For $b=(b_\alpha)\in\mathbb Z^\ell$ with $0\le b\le a$ coordinatewise, set
\[
g_b:=\prod_\alpha p_\alpha^{b_\alpha},\qquad I_b:=I(g_b).
\]
Thus $I_0\cong I_a\cong R$. Theorem \ref{thm:divisor-class-group} identifies $[I_b]$ with $\overline b\in\mathbb Z^\ell/\mathbb Za$.

\begin{Lem}\label{lem:rank-one-mcm-classification}
Every rank-one maximal Cohen--Macaulay $R$-module is isomorphic to $I_b$ for some $0\le b\le a$. If it is nonfree, then $b\ne0,a$ and this vector $b$ is unique.
\end{Lem}
\begin{proof}
If $f\notin\n^2$, then $R$ is regular local and the assertion follows. Otherwise let $N$ be a nonfree rank-one maximal Cohen--Macaulay module. By hyperbolic Kn\"orrer periodicity \cite[Theorem 2.4]{Leu26}, it is stably isomorphic to the cokernel of
\[
\Phi=\begin{pmatrix}\alpha&uI_s\\vI_s&\beta\end{pmatrix},
\qquad \alpha\beta=\beta\alpha=fI_s,
\]
for a reduced matrix factorization over $S$. Both $\Phi$ and its companion matrix have entries in $\m$, so this cokernel has no free summand. Krull--Schmidt cancellation therefore identifies it with $N$. Over $\Frac R$, the block $uI_s$ is invertible and the relation $\alpha\beta=uvI_s$ makes $\Phi$ have rank $s$. Its cokernel also has rank $s$, whence $s=1$. Write $\alpha=g$, $\beta=h$, with $gh=f$. The cokernel is $(v,g)$, which becomes $(u,h)$ after multiplication by $u/g$ in the fraction field. Thus $N\cong I_b$ for the exponent vector of $h$.

Finally, if $0\le b,c\le a$ and $I_b\cong I_c$, the class-group calculation gives $b-c=ta$ for some integer $t$. Within this box a nonzero $t$ is possible only for $(b,c)=(a,0)$ or $(0,a)$. These endpoints give the free module $R$, so they are excluded in the nonfree case.
\end{proof}

\begin{Lem}\label{lem:factor-hom-compatibility}
For $0\le b,c\le a$, one has
\begin{equation}\label{eq:factor-hom-compatibility}
\Hom_R(I_b,I_c)\in\CM R
\quad\Longleftrightarrow\quad b\le c\text{ or }c\le b.
\end{equation}
In particular, for divisors $g,h$ of $f$, the modules $(u,g)$ and $(u,h)$ are compatible precisely when $g\mid h$ or $h\mid g$, up to units.
\end{Lem}
\begin{proof}
The Hom module is rank-one reflexive, with class $\overline{c-b}$. By Lemma \ref{lem:rank-one-mcm-classification}, it is maximal Cohen--Macaulay exactly when
\[
c-b=z+ta\qquad\text{for some }0\le z\le a,\quad t\in\mathbb Z.
\]
For $t\ge0$ all coordinates on the right are nonnegative; for $t\le-1$ all are nonpositive. Thus $b,c$ must be comparable. Conversely, if $b\le c$ take $z=c-b$, $t=0$; if $c\le b$ take $z=c-b+a$, $t=-1$. The class-group calculation identifies the Hom module with $I_z$, proving the converse as well.
\end{proof}

A \emph{flag} is a chain $\F=(b_1,\ldots,b_m)$ of exponent vectors, with endpoints
\[
0=b_0<b_1<\cdots<b_m<b_{m+1}=a.
\]
Here $b<c$ means $b\le c$ coordinatewise and $b\ne c$. The empty flag is allowed. Let $\mathsf{Flag}(a)$ be the poset of flags ordered by inclusion, and put
\[
T^\F:=R\oplus\bigoplus_{j=1}^m I_{b_j}.
\]
A flag is \emph{maximal} if its chain is saturated in $[0,a]\cap\mathbb Z^\ell$, equivalently if $m=n-1$.

For the labelled prime factors, put $[n]:=\{1,\ldots,n\}$ and $\mu(I):=(m_\alpha(I))_\alpha$, using \eqref{eq:factor-ideal-classes}; thus $T_I=I_{\mu(I)}$. A \emph{labelled lift} of $\F$ is a chain $\emptyset\subsetneq I_1\subsetneq\cdots\subsetneq I_m\subsetneq[n]$ with $\mu(I_j)=b_j$. Every flag admits such a lift, and permuting associate factors changes only the lift. Write $s_j:=(j\ j+1)\in\mathfrak S_n$ for $1\le j<n$. For $\omega\in\mathfrak S_n$, put $I_j^\omega:=\{\omega(1),\ldots,\omega(j)\}$ for $0\le j\le n$ and write
\[
\F_\omega:=(\mu(I_1^\omega),\ldots,\mu(I_{n-1}^\omega)),\qquad
T^\omega:=T^{\F_\omega}.
\]

\begin{Thm}[Classification of modifying modules]\label{thm:modifying-classification}
Let $M$ be a nonzero basic modifying $R$-module with $r$ indecomposable summands, and distinguish one of them, $L$. Then $L$ has rank one, and there is a unique flag $\F=(b_1,\ldots,b_{r-1})$ with endpoints
\begin{equation}\label{eq:anchored-modifying-chain}
0=b_0<b_1<\cdots<b_{r-1}<b_r=a
\end{equation}
such that
\[
M\cong\bigoplus_{j=0}^{r-1}\mathsf T_L(I_{b_j}).
\]
Conversely, every rank-one reflexive module $L$ and every flag give a basic modifying module by this formula.
\end{Thm}
\begin{proof}
By Theorem \ref{modrk1}, every indecomposable summand of $M$ has rank one. Put $G:=\mathsf T_{L^*}(M)$, a basic modifying module with distinguished summand $R$. Each summand $N$ of $G$ is maximal Cohen--Macaulay, since $N=\Hom_R(R,N)$ is a direct summand of $\End_R(G)$. Lemma \ref{lem:rank-one-mcm-classification} writes every nonfree summand uniquely as $I_b$ with $0<b<a$. Lemma \ref{lem:factor-hom-compatibility} makes these vectors a chain; adjoining $b_0=0$ and $b_r=a$ gives \eqref{eq:anchored-modifying-chain} uniquely.

Conversely, the modules $I_{b_0},\ldots,I_{b_{r-1}}$ are pairwise nonisomorphic by Lemma \ref{lem:rank-one-mcm-classification}, and all their pairwise Hom modules are maximal Cohen--Macaulay by Lemma \ref{lem:factor-hom-compatibility}. Their sum is therefore basic modifying. Reflexive twisting preserves these properties, proving the converse.
\end{proof}

Thus the nonzero basic modifying modules are precisely
\begin{equation}\label{eq:modifying-flag-classification}
\mathsf T_L(T^\F),\qquad L\in\Cl(R),\quad\F\text{ a flag}.
\end{equation}
For generators, distinguishing $R$ gives a bijection $\F\mapsto T^\F$ from flags to basic modifying generators. For general modules, uniqueness is relative to the specified summand $L$; changing that summand can give a different twist presentation.

\begin{Cor}[Maximal modifying modules]\label{cor:mm-classification}
With the notation of Theorem \ref{thm:modifying-classification},
\[
M\text{ is maximal modifying}
\quad\Longleftrightarrow\quad
b_\bullet\text{ is saturated}
\quad\Longleftrightarrow\quad r=n,
\]
where saturation is taken in $[0,a]\cap\mathbb Z^\ell$. Consequently, the basic maximal modifying modules are precisely
\begin{equation}\label{eq:mm-classification}
\mathsf T_L(T^\omega)
=\bigoplus_{j=0}^{n-1}(L\otimes_RT_{I_j^\omega})^{**},
\qquad L\in\Cl(R),\quad\omega\in\mathfrak S_n,
\end{equation}
and the basic maximal modifying generators are precisely the $T^\omega$.
\end{Cor}
\begin{proof}
Reflexive twisting preserves maximality, so assume $L=R$. If the chain is not saturated, insert an intermediate vector. Theorem \ref{thm:modifying-classification} then gives a modifying module containing $M$ and a new summand, so $M$ is not maximal modifying. Conversely, if $M\oplus N$ is modifying, its basic representative has a chain containing $b_\bullet$, by the uniqueness in that theorem with distinguished summand $R$. If $b_\bullet$ is saturated, no new vector can occur, so $N\in\add M$.

Each strict step increases the sum of coordinates by at least one. Since $\sum_\alpha a_\alpha=n$, the chain is saturated exactly when $r=n$, with each step adding a standard basis vector. Such chains specify words with $a_\alpha$ occurrences of each $p_\alpha$, giving \eqref{eq:mm-classification} through their labelled lifts. Distinguishing $R$ gives the assertion for generators.
\end{proof}

For $S=k[[x,y]]$, the generator cases of Theorem \ref{thm:modifying-classification} and Corollary \ref{cor:mm-classification} recover \cite[Theorem 5.1(1), (2)]{IW18}. For $S=\mathbb C[[x,y]]$, \eqref{eq:mm-classification} agrees with \cite[Theorem 5.20]{IW18}. The arguments above use the rank-one theorem over the present base $S$ and characteristic-free hyperbolic periodicity.

\begin{Cor}\label{cor:mm-counts}
There are $n!/\prod_\alpha a_\alpha!$ basic maximal modifying generators.
\end{Cor}
\begin{proof}
By Theorem \ref{thm:modifying-classification} with distinguished summand $R$ and Corollary \ref{cor:mm-classification}, these generators are indexed by maximal flags. Such flags are words with $a_\alpha$ copies of each $p_\alpha$, giving the multinomial count.
\end{proof}

\subsection{Block exchanges and mutation connectivity}\label{sec:modifying-mutation}

Fix an ordered block factorization $f=g_1\cdots g_r$, with $r\ge2$ and $0\ne g_i\in\n$, and put
\[
G_i:=g_1\cdots g_i,\quad G_0=1,\qquad
M:=\bigoplus_{i=0}^{r-1}I(G_i).
\]
Let $\F=(b_1,\ldots,b_{r-1})$, where $b_{i,\alpha}=\operatorname{ord}_{p_\alpha}(G_i)$; thus $M=T^\F$. For $1\le j<r$, abbreviate
\[
g=G_{j-1},\quad s=g_j,\quad t=g_{j+1},\quad h=f/(gst),
\qquad N=\bigoplus_{\substack{0\le i<r\\i\ne j}}I(G_i).
\]
The following extends the exchange calculation of \cite[Lemmas 5.12--5.14]{IW18} to our regular complete local base.

\begin{Lem}[Block exchanges]\label{lem:factor-exchange-mutation}
There are exact sequences
\begin{equation}\label{eq:factor-right-exchange}
0\longrightarrow I(gt)\xrightarrow{\binom{1}{s}}I(g)\oplus I(gst)
\xrightarrow{(-s\ \ 1)}I(gs)\longrightarrow0,
\end{equation}
\begin{equation}\label{eq:factor-left-exchange}
0\longrightarrow I(gs)\xrightarrow{\binom{1}{t}}I(g)\oplus I(gst)
\xrightarrow{(-t\ \ 1)}I(gt)\longrightarrow0.
\end{equation}
The last map in \eqref{eq:factor-right-exchange} is a minimal right $(\add N)$-approximation. The first map in \eqref{eq:factor-left-exchange} is a minimal left $(\add N)$-approximation. Hence
\begin{equation}\label{eq:factor-mutation-formula}
\mu_{I(gs)}^+(M)\cong\mu_{I(gs)}^-(M)\cong N\oplus I(gt).
\end{equation}
Here $I(f)=(u)\cong R$, so the middle terms belong to $\add N$ also at the endpoint.
\end{Lem}
\begin{proof}
Surjectivity follows from $sI(g)+I(gst)=I(gs)$. The kernel is $I(gt)$ because $(I(gst):s)\cap I(g)=I(gt)$: reduction modulo $u$ and cancellation of $s$ in $S[[v]]$ give this identity, using $f=gsth$. Interchanging $s,t$ gives the second sequence.

As fractional ideals in $\Frac R$, one has
\begin{equation}\label{eq:comparable-ideal-hom}
H(x,y):=\Hom_R(I(x),I(y))=
\begin{cases}
I(y/x),&x\mid y,\\
u^{-1}I(yf/x),&y\mid x.
\end{cases}
\end{equation}
For $x\mid y$, the inclusion $I(y)\subseteq I(x)$ puts the multiplier in $\End_R(I(x))=R$; reducing $xr\in(u,y)$ modulo $u$, lifting to $S[[v]]$ and cancelling $x$ proves the first formula. The identity $I(z)=(uz/f)(v,f/z)$ gives the second by interchanging $u,v$. Every remaining prefix $c$ satisfies $c\mid g$ or $gst\mid c$, so these formulas give
\[
sH(c,g)+H(c,gst)=H(c,gs),\qquad
H(g,c)+tH(gst,c)=H(gs,c).
\]
These are the two approximation properties.

If $I(g)\not\cong I(gst)$, redundancy of a component would require $s\in I(st)$ or $u/s\in I(gh)$ for the right map, and $t\in I(st)$ or $u/t\in I(gh)$ for the left map. Reduction modulo $u$ rules out the first alternatives, since $s,t$ are nonunits; reduction modulo $s$ or $t$ rules out the second. Their endomorphism rings are local, so both approximations are minimal. If $I(g)\cong I(gst)$, Theorem \ref{thm:divisor-class-group} forces $g=h=1$ and $r=2$. The right map is then $R^2\xrightarrow{(-s\ \ u)}I(s)$, and the left map is dual to the generating map of $I(s)^*=u^{-1}I(t)$. Both are minimal because $I(s),I(t)$ are nonfree two-generated ideals, again by the class-group calculation. Finally, the kernels of the right approximation and of the dual left approximation give \eqref{eq:factor-mutation-formula} by \eqref{eq:iw-mutation-convention}.
\end{proof}

In flag coordinates, the exchange replaces $b_j$ by $b_{j-1}+b_{j+1}-b_j$. Cyclic reanchoring covers vertex zero: twisting by $I(g_1)^*$ and ordering the summands as $1,\ldots,r-1,0$ gives $R,I(g_2),\ldots,I(g_2\cdots g_r)$ with blocks $(g_2,\ldots,g_r,g_1)$. Lemma \ref{lem:twist-preserves-mutation} therefore extends the calculation to every vertex and every common reflexive rank-one twist.

\begin{Cor}[MM mutation graph]\label{cor:mm-mutation-connectivity}
The irreducible mutation graph of basic maximal modifying modules is connected. If $n\ge2$, then for every $M\in\MM R$ and every indecomposable summand $X$,
\[
\mu_X^+(M)\cong\mu_X^-(M)\in\MM R.
\]
\end{Cor}
\begin{proof}
For $n=1$, Corollary \ref{cor:mm-classification} gives $\MM R=\{R\}$. Suppose $n\ge2$. Choose a summand other than $X$ as anchor in that corollary. Lemma \ref{lem:factor-exchange-mutation} and reflexive twisting then give the mutation formula; the exchanged prime-factor word is again maximal.

For a rank-one reflexive module $L$, let $\mathcal G_L$ be the subgraph on MM modules containing $L$, using mutations away from $L$. The classification of MM generators and the adjacent exchanges above make $\mathcal G_R$ connected. Lemma \ref{lem:twist-preserves-mutation} identifies its twist by $L$ with $\mathcal G_L$, so each $\mathcal G_L$ is connected.

A prime-factor word beginning with $p_\alpha$ gives an MM generator $N_\alpha$ containing $R$ and $\mathfrak P_\alpha$. Hence $\mathsf T_L(N_\alpha)$ belongs to both $\mathcal G_L$ and $\mathcal G_{\mathsf T_L(\mathfrak P_\alpha)}$. The classes $[\mathfrak P_\alpha]$ generate $\Cl(R)$ by Theorem \ref{thm:divisor-class-group}, so all these subgraphs lie in one connected component. They cover $\MM R$ by Theorem \ref{modrk1}, proving connectivity.
\end{proof}

\section{Flag models and the Auslander--McKay correspondence}\label{sec:flag-models}

We construct a projective crepant model $Y^\F$ and a tautological tilting bundle $\E^\F$ for every flag $\F$. Theorem \ref{thm:flag-models} describes the model and its closed fibre, while Proposition \ref{prop:tautological-tilting} identifies its canonical generator, pushforward and endomorphism algebra. Corollary \ref{cor:flag-product-blowup} expresses the model as the blowup of a product of factor ideals. These constructions prove the Auslander--McKay equivalence in Section \ref{sec:relative-one-model-classification}: the flag models exhaust $\mathsf{Crep}_1(R)$. Finally, their local equations and tilting equivalences characterize the modifying algebras of finite global dimension and give the NCCR criterion in Section \ref{sec:nccr}.

The successive blowup construction and its tilting realization were established by Iyama--Wemyss for $S=k[[x,y]]$, with $k$ algebraically closed of characteristic zero; see \cite[Section 5.1 and Theorem 5.2]{IW14} and \cite[Section 5.5 and Theorem 5.35]{IW18}. Here we carry out the construction and prove its properties over an arbitrary regular complete local ring $S$, including mixed characteristic.

For a flag $\F=(b_1,\ldots,b_m)$ as in Section \ref{sec:modifying-classification}, put
\[
G_0:=1,\qquad G_j:=g_{b_j}\quad(1\le j\le m),\qquad G_{m+1}:=f,
\]
and define its blocks by
\begin{equation}\label{eq:flag-blocks}
g_i:=G_i/G_{i-1},\qquad Q_i:=\prod_{j>i}g_j
\quad(1\le i\le m+1).
\end{equation}
Thus $f=G_{i-1}g_iQ_i$ and every $g_i$ is a nonzero nonunit of $S$. The unit $\varepsilon$ is included in the last block. For a maximal flag, all blocks are prime elements.

Starting with $Y_0^\F:=\Spec R$, define successively
\begin{equation}\label{eq:flag-blowups}
\rho_j\colon Y_j^\F
:=\operatorname{Bl}_{I_{b_j}\O_{Y_{j-1}^\F}}(Y_{j-1}^\F)
\longrightarrow Y_{j-1}^\F
\qquad(1\le j\le m).
\end{equation}
Here $I_{b_j}\O_{Y_{j-1}^\F}$ means its image as an ideal of $\O_{Y_{j-1}^\F}$. Write $Y^\F:=Y_m^\F$ and $\pi:=\rho_1\circ\cdots\circ\rho_m$. For the empty flag, $\pi$ is the identity of $\Spec R$. Write $Y^\omega:=Y^{\F_\omega}$ for $\omega\in\mathfrak S_n$.

\begin{Thm}\label{thm:flag-models}
For every flag $\F=(b_1,\ldots,b_m)$, the construction \eqref{eq:flag-blowups} gives a projective birational morphism
\[
\pi\colon Y^\F\longrightarrow\Spec R
\]
from a noetherian normal Gorenstein integral scheme with $\omega_{Y^\F}\cong\O_{Y^\F}$. The following hold.
\begin{enumerate}
\item $\mathbb R\pi_*\O_{Y^\F}=\O_{\Spec R}$.
\item If $m\ge1$, the scheme-theoretic closed fibre is a reduced split nodal chain
\[
Y^\F\times_Rk=C_1\cup\cdots\cup C_m,\qquad C_j\cong\mathbb P^1_k,
\]
with consecutive components meeting at one $k$-rational ordinary node and no other intersections. For $m=0$, the closed fibre is $\Spec k$.
\item If $\F$ is maximal, then $Y^\F$ is locally factorial.
\end{enumerate}
In particular, every flag model belongs to $\mathsf{Crep}_1(R)$, and every nonempty flag also satisfies \eqref{eq:split-chain}.
\end{Thm}

We prove the theorem below by calculating the blowups before completion. This also keeps track of the scheme structure of the closed fibre. No assumption on the characteristic or perfection of $k$ is needed. The converse within $\mathsf{Crep}_1(R)$ is Theorem \ref{thm:auslander-mckay}; the split chain is a consequence of that classification, rather than an assumption on the closed fibre of the given model.

\subsection{Elementary blowups and affine charts}

\begin{Lem}\label{lem:elementary-flag-blowup}
Let $0\ne g,h\in\n$, put $A_{g,h}:=S[s,t]/(st-gh)$, and let
\[
\beta\colon Z_{g,h}:=\operatorname{Bl}_{(s,g)}(\Spec A_{g,h})
\longrightarrow\Spec A_{g,h}.
\]
Then $Z_{g,h}$ has two affine charts
\begin{equation}\label{eq:elementary-blowup-charts}
\begin{aligned}
Z_1&=\Spec S[s,z]/(sz-g),&\qquad t&=hz,\\
Z_2&=\Spec S[w,t]/(wt-h),& s&=gw.
\end{aligned}
\end{equation}
Their intersection is $\Spec S[z,z^{-1}]$, with $w=z^{-1}$. The morphism $\beta$ is projective, crepant and birational, and $Z_{g,h}$ is normal and Gorenstein.
\end{Lem}
\begin{proof}
The blowup charts are the subalgebras $A_{g,h}[g/s]$ and $A_{g,h}[s/g]$ of $\Frac A_{g,h}$; see \cite[Tag 0804]{stacks-project}. On the first put $z=g/s$, so $sz=g$ and $t=hz$; on the second put $w=s/g$, so $s=gw$ and $wt=h$. These give \eqref{eq:elementary-blowup-charts}: the displayed hypersurface rings embed in the fraction field, as is seen by inverting $s$ and $t$, respectively. On the intersection $z$ is invertible, $s=g/z$ and $t=hz$.

For any $0\ne a\in\n$, the ring $S[x,y]/(xy-a)$ is a normal domain. Indeed, it is a hypersurface domain, hence Cohen--Macaulay, and the opens $D(x)$ and $D(y)$ are regular. Their complement has codimension two, since its quotient is $S/(a)$. Serre's criterion proves normality. This applies to the base and both charts. Projectivity and birationality follow from the blowup construction, and the charts are Gorenstein hypersurfaces.

For crepancy, put $X:=\Spec A_{g,h}$ and $Z:=Z_{g,h}$. For $0\ne a\in S$, hypersurface adjunction \cite[Tag 0AA4]{stacks-project} gives the relative dualizing module of $S[x,y]/(xy-a)$ a residue generator $dx/x=-dy/y$. Here $d$ is relative to $S$, and these are rational expressions for a dualizing generator that remains regular and nonvanishing along $x=y=0$. Thus the relative dualizing sheaves on the base and the two charts have generators
\[
\begin{aligned}
\eta_X&=\frac{ds}{s}=-\frac{dt}{t},\\
\eta_1&=\frac{ds}{s}=-\frac{dz}{z},\\
\eta_2&=\frac{dw}{w}=-\frac{dt}{t}.
\end{aligned}
\]
Since $s=gw$ on $Z_2$ and $dg=0$, the rational pullback of $\eta_X$ is $\eta_i$ on each $Z_i$. These generators agree on the overlap by $w=z^{-1}$. Hence the birational identification extends to $\beta^*\omega_{X/S}\cong\omega_{Z/S}$, proving crepancy.
\end{proof}

We now iterate the elementary blowup and then pass to the complete local base. Put
\[
A:=S[u,v]/(uv-f),\qquad \mathfrak q:=(\n,u,v)A.
\]
There is a canonical identification $\widehat{A_{\mathfrak q}}=R$, so $A\to R$ is flat. Perform the blowups \eqref{eq:flag-blowups} first over $\Spec A$, using the ideals $(u,G_j)$, and denote the resulting morphism by
\[
\widetilde\pi\colon\widetilde Y^\F\longrightarrow\Spec A.
\]
Formation of the Rees algebra and of the blowup commutes with flat base change \cite[Tag 0805]{stacks-project}. Hence
\begin{equation}\label{eq:flag-flat-base-change}
Y^\F\cong\widetilde Y^\F\times_{\Spec A}\Spec R.
\end{equation}

\begin{Prop}\label{prop:flag-charts}
The scheme $\widetilde Y^\F$ is covered by the $m+1$ affine charts
\begin{equation}\label{eq:algebraic-flag-charts}
\widetilde U_i=\Spec S[x_i,y_i]/(x_iy_i-g_i),\qquad
u=G_{i-1}x_i,\quad v=Q_iy_i
\quad(1\le i\le m+1).
\end{equation}
The adjacent charts are glued on $D(y_i)\subseteq\widetilde U_i$ and $D(x_{i+1})\subseteq\widetilde U_{i+1}$ by
\begin{equation}\label{eq:flag-chart-transitions}
x_{i+1}=y_i^{-1},\qquad y_{i+1}=g_{i+1}y_i.
\end{equation}
Consequently $Y^\F$ has the affine cover $U_i=\Spec B_i$, where
\begin{equation}\label{eq:completed-base-flag-charts}
B_i=S[[u,v]][x_i,y_i]/(u-G_{i-1}x_i,\ v-Q_iy_i,\ x_iy_i-g_i).
\end{equation}
\end{Prop}
\begin{proof}
For the empty flag all assertions are immediate, so assume $m\ge1$.
Apply Lemma \ref{lem:elementary-flag-blowup} first to $f=g_1(g_2\cdots g_{m+1})$. Its first chart is $\widetilde U_1$; its second is the remaining chart
\[
\Spec S[s,v]/(sv-g_2\cdots g_{m+1}),\qquad u=g_1s.
\]
More generally, after $j-1$ steps the charts already obtained are $\widetilde U_1,\ldots,\widetilde U_{j-1}$ and the remaining chart is
\[
\Spec S[s,v]/(sv-g_j\cdots g_{m+1}),\qquad u=G_{j-1}s.
\]
On $\widetilde U_i$ for $i<j$, the ideal $(u,G_j)$ is principal, since
\[
G_j=u\,y_i\prod_{i<\ell\le j}g_\ell.
\]
The next blowup is therefore the identity there. On the remaining chart its centre is
\[
(u,G_j)=G_{j-1}(s,g_j).
\]
Multiplication of an ideal by a nonzero principal factor does not change its blowup. Lemma \ref{lem:elementary-flag-blowup}, applied to $g_j$ and $g_{j+1}\cdots g_{m+1}$, therefore produces $\widetilde U_j$ and the next remaining chart. Induction proves \eqref{eq:algebraic-flag-charts}, and the coordinate change in that lemma gives \eqref{eq:flag-chart-transitions}. Base change gives \eqref{eq:completed-base-flag-charts}.
\end{proof}

Each ring in \eqref{eq:completed-base-flag-charts} is an affine chart ring obtained by base change to $R$. Its further completion at a distinguished point of the closed fibre is calculated below.

\subsection{The closed fibre and local factoriality}

\begin{Prop}\label{prop:flag-fibre}
If $m\ge1$, the scheme-theoretic closed fibre of $\pi$ is the reduced split nodal chain of $m$ projective lines in Theorem \ref{thm:flag-models}. There are distinguished $k$-rational points $p_i\in U_i$ for $1\le i\le m+1$, given by $x_i=y_i=0$ on the fibre, such that
\begin{equation}\label{eq:flag-completed-local-rings}
\widehat{\O}_{Y^\F,p_i}\cong S[[x_i,y_i]]/(x_iy_i-g_i).
\end{equation}
The internal points $p_2,\ldots,p_m$ are its nodes, and $p_1,p_{m+1}$ are its two end points. Every other closed point of the fibre is a regular point of $Y^\F$. For $m=0$, the closed fibre consists of one point $p_1$, and \eqref{eq:flag-completed-local-rings} still holds, with $g_1=f$.
\end{Prop}
\begin{proof}
Reduce \eqref{eq:completed-base-flag-charts} modulo $\m=(\n,u,v)R$. Since all $g_i$ belong to $\n$, for $m\ge1$ this gives
\begin{equation}\label{eq:flag-fibre-charts}
\begin{split}
U_1\times_R k&=\Spec k[y_1],\qquad x_1=0,\\
U_i\times_R k&=\Spec k[x_i,y_i]/(x_iy_i)
\qquad(2\le i\le m),\\
U_{m+1}\times_R k&=\Spec k[x_{m+1}],\qquad y_{m+1}=0.
\end{split}
\end{equation}
The axis $x_i=0$ in the $i$th chart and the axis $y_{i+1}=0$ in the next chart glue by $x_{i+1}=y_i^{-1}$ to a copy $C_i$ of $\mathbb P^1_k$. Internal origins are ordinary split nodes. Nonadjacent charts have no common points on the fibre: an overlap across two successive transitions requires the intervening $g_i$ to be invertible. Thus \eqref{eq:flag-fibre-charts} describes the whole scheme-theoretic fibre, with no additional components or nilpotent structure.

At $p_i$, completion of \eqref{eq:completed-base-flag-charts} allows substitution of $u=G_{i-1}x_i$ and $v=Q_iy_i$, giving \eqref{eq:flag-completed-local-rings}. Equivalently, the completed local rings before and after the base change \eqref{eq:flag-flat-base-change} agree at points over $\mathfrak q$. A different closed point on the fibre lies in $D(x_i)\cup D(y_i)$ on an algebraic chart \eqref{eq:algebraic-flag-charts}. This open is regular, being an open of $\Spec S[x_i,x_i^{-1}]$ or $\Spec S[y_i,y_i^{-1}]$. Its completed local ring is therefore regular, and the same is true after base change; hence $Y^\F$ is regular at this point. When $m=0$, $Y^\F=\Spec R$ and the assertion follows directly.
\end{proof}

\begin{proof}[Proof of Theorem \ref{thm:flag-models}]
Each step in \eqref{eq:flag-blowups} is the blowup of a nonzero ideal on an integral scheme. Thus $Y^\F$ is a noetherian integral scheme, and $\pi$ is projective and birational. Proposition \ref{prop:flag-fibre} proves (2); in particular, the closed fibre has dimension at most one.

At each distinguished point $p_i$, the completed local ring \eqref{eq:flag-completed-local-rings} is a Gorenstein normal domain, by the hypersurface argument used for $R$ in the preceding section. At every other closed point of the fibre it is regular, by Proposition \ref{prop:flag-fibre}. Normality and the Gorenstein property descend from the completion. These properties then hold everywhere on $Y^\F$: every point specializes to a closed point on the closed fibre, since $Y^\F$ is proper over the local scheme $\Spec R$, and its local ring is a localization of the local ring at such a specialization.

Fix a generator of $\omega_S$ and suppress its tensor factor in the residue notation. On $\widetilde U_i$, Lemma \ref{lem:elementary-flag-blowup} identifies the rational pullback of the canonical generator $\eta_A=du/u$ of $\omega_{\Spec A}$ with
\[
\eta_i=\frac{dx_i}{x_i}=-\frac{dy_i}{y_i},
\]
since $u=G_{i-1}x_i$ and $G_{i-1}\in S$. These nonvanishing canonical generators agree on overlaps by $x_{i+1}=y_i^{-1}$, so they trivialize $\omega_{\widetilde Y^\F}$.

Let $\eta$ be the residue generator of $\omega_R$ obtained by completing $\eta_A$. For every closed point $y$ of the closed fibre and the corresponding point $\widetilde y\in\widetilde Y^\F$, the proof of Proposition \ref{prop:flag-fibre} gives
\[
C_y:=\widehat{\O}_{Y^\F,y}
\cong\widehat{\O}_{\widetilde Y^\F,\widetilde y},
\]
with the same $S$-algebra structure and coordinate formulas $u=G_{i-1}x_i$, $v=Q_iy_i$. Canonical modules commute with local completion \cite[Tag 0AWD]{stacks-project}, and completing hypersurface adjunction preserves the residue generators compared in Lemma \ref{lem:elementary-flag-blowup}. Hence the rational section induced by $\eta$ corresponds to the completed generator $\eta_i$ of $\omega_{C_y}$. Faithful flatness of $\O_{Y^\F,y}\to C_y$ shows that $\eta$ generates $(\omega_{Y^\F})_y$. Since every point specializes to such a $y$, the section $\eta$ trivializes $\omega_{Y^\F}$. Lemma \ref{lem:crepancy-and-duality} then gives $\pi^!\O_{\Spec R}\cong\O_{Y^\F}$.

Assertion (1) follows from Proposition \ref{prop:crepant-vanishing}, but also has a short direct proof. For the elementary blowup $\beta$ of Lemma \ref{lem:elementary-flag-blowup}, normality gives $\beta_*\O_{Z_{g,h}}=\O_{\Spec A_{g,h}}$. The two affine charts have intersection $\Spec S[z,z^{-1}]$, and their images contain $S[z]$ and $S[z^{-1}]$, whose sum is $S[z,z^{-1}]$. Thus the \v{C}ech complex has no positive cohomology, so $\mathbb R\beta_*\O_{Z_{g,h}}=\O_{\Spec A_{g,h}}$. At each step of $\widetilde\pi$, the blowup has this form on the remaining chart and is the identity on the other charts. Composition and flat base change in \eqref{eq:flag-flat-base-change} prove (1).

Finally, suppose $\F$ is maximal. Then each $g_i$ is a prime element of the regular complete local ring $S$. Corollary \ref{cor:class-group-factoriality}, applied to the completed local equations \eqref{eq:flag-completed-local-rings}, gives
\[
\Cl\bigl(\widehat{\O}_{Y^\F,p_i}\bigr)
=\Cl\bigl(S[[x_i,y_i]]/(x_iy_i-g_i)\bigr)=0.
\]
Thus these completed local rings are factorial. Factoriality descends from the completion: for a noetherian normal local ring $B$ with factorial completion, the completion of a rank-one reflexive $B$-module is reflexive and hence free, and faithful flatness descends freeness. Consequently, $\O_{Y^\F,p_i}$ is factorial. By Proposition \ref{prop:flag-fibre}, all other closed points of the fibre are regular, hence factorial. Every point of $Y^\F$ specializes to a closed point of the fibre, as above, so localization proves that $Y^\F$ is locally factorial. This proves (3) and completes the theorem.
\end{proof}

\begin{Cor}\label{cor:flag-type-a}
For any nonempty flag with $m$ terms, Theorem \ref{thm:type-a-classification} applies to $Y^\F$. In particular, every indecomposable pretilting bundle on $Y^\F$ is a line bundle, and there are $m!$ basic tilting bundles modulo tensoring by line bundles.
\end{Cor}
\begin{proof}
Theorem \ref{thm:flag-models} verifies both \eqref{eq:geometric-setup} and \eqref{eq:split-chain}.
\end{proof}

\begin{Rem}
Maximality is necessary for the local factoriality assertion. For example, for $S=k[[x,y]]$ and $f=xy$, the empty flag gives $Y^\emptyset=\Spec k[[x,y,u,v]]/(uv-xy)$. The height-one prime ideal $(u,x)$ needs two generators, so this local ring is not factorial.
\end{Rem}

\subsection{Tautological tilting bundles}\label{sec:tautological-tilting}

We identify the tautological bundle with the canonical generator of Proposition \ref{prop:canonical-tilting-general}. Its endomorphism algebra is then given by Lemma \ref{lem:pretilting-endomorphisms}(1).

For the model of any flag, write
\[
\Psi(U):=(\pi_*U)^{**}.
\]
The morphism $\pi$ is an isomorphism over $D(u)\cup D(v)$: on $D(u)$ every centre is the unit ideal, and on $D(v)$ the ideal $(u,G_j)$ is generated by $G_j$. This open contains every codimension-one point of $\Spec R$.

\begin{Prop}\label{prop:tautological-tilting}
For a flag $\F=(b_1,\ldots,b_m)$, put
\[
J_j:=I_{b_j}\O_{Y^\F}\quad(1\le j\le m),\qquad
\E^\F:=\O_{Y^\F}\oplus\bigoplus_{j=1}^mJ_j.
\]
Then $J_j$ is an invertible ideal sheaf, with
\begin{equation}\label{eq:tautological-degrees}
\deg(J_j|_{C_i})=\delta_{ij},\qquad \pi_*J_j=I_{b_j}.
\end{equation}
The bundle $\E^\F$ is the canonical basic tilting generator $\P_{Y^\F}$, and
\begin{equation}\label{eq:tautological-endomorphisms}
\pi_*\E^\F=T^\F,\qquad
\End_{Y^\F}(\E^\F)\cong\End_R(T^\F).
\end{equation}
For the empty flag these statements mean $\E^\emptyset=\O_{\Spec R}$.
\end{Prop}
\begin{proof}
The empty-flag case is immediate. Assume $m>0$. On the chart $U_i$ in \eqref{eq:completed-base-flag-charts}, the ideal $J_j$ is
\[
J_j|_{U_i}=\begin{cases}(u),&i\le j,\\(G_j),&i>j.\end{cases}
\]
These generators are non-zero divisors. The only change of generator along the closed fibre occurs between $U_j$ and $U_{j+1}$, where $G_j=y_ju$. With coordinates $y_j$ and $x_{j+1}=y_j^{-1}$ on $C_j$, this is the transition for $\O_{\mathbb P^1}(1)$, proving the degree formula.

The inclusions $I_{b_j}\subseteq\pi_*J_j\subseteq\pi_*\O_{Y^\F}=R$ become equalities with $I_{b_j}$ in codimension one. Since $I_{b_j}$ is reflexive, intersection inside $\Frac R$ over the height-one localizations gives $\pi_*J_j=I_{b_j}$.

The degree formula \eqref{eq:tautological-degrees} identifies the summand classes of $\E^\F$ with the vertices of the fundamental type $A$ alcove. Theorem \ref{thm:type-a-classification} therefore makes $\E^\F$ a basic tilting bundle. It is globally generated and contains $\O_{Y^\F}$, so Proposition \ref{prop:canonical-tilting-general} identifies it with $\P_{Y^\F}$. Lemma \ref{lem:pretilting-endomorphisms}(1) gives the endomorphism assertion.
\end{proof}

\begin{Cor}[Reconstruction of flag models]\label{cor:flag-product-blowup}
For a flag $\F=(b_1,\ldots,b_m)$ and positive integers $c_1,\ldots,c_m$, there is an isomorphism over $X=\Spec R$
\begin{equation}\label{eq:flag-product-blowup}
Y^\F\cong\operatorname{Bl}_{\prod_{j=1}^mI_{b_j}^{c_j}}X.
\end{equation}
All products here are ordinary products of ideals. For the empty flag the product is $R$ and the blowup is $X$.
\end{Cor}
\begin{proof}
Put $E_{\mathbf c}:=\O_{Y^\F}\oplus\bigoplus_{j=1}^mJ_j^{\oplus c_j}$. By Proposition \ref{prop:tautological-tilting}, this bundle is globally generated with $\add E_{\mathbf c}=\add\P_{Y^\F}$, so Proposition \ref{prop:canonical-tilting-general} makes $\det E_{\mathbf c}$ relatively ample. Moreover,
\[
\pi_*E_{\mathbf c}=R\oplus\bigoplus_{j=1}^mI_{b_j}^{\oplus c_j},
\qquad
\mathfrak d(\pi_*E_{\mathbf c})=\prod_{j=1}^mI_{b_j}^{c_j},
\]
the determinant image being the ordinary product of the rank-one ideals. Theorem \ref{thm:flag-models} puts $Y^\F$ in $\mathsf{Crep}_1(R)$, so Proposition \ref{cor:module-reconstruction}(2) gives \eqref{eq:flag-product-blowup}.
\end{proof}

\begin{Cor}[Tautological derived realization]\label{cor:flag-derived-equivalence}
For every flag $\F$, its tautological bundle gives an $R$-linear equivalence
\[
D^b(\Coh Y^\F)\simeq D^b(\mod\End_R(T^\F)),
\]
restricting to $\per Y^\F\simeq\per\End_R(T^\F)$.
\end{Cor}
\begin{proof}
Apply Lemma \ref{lem:tilting-derived-realization} and Proposition \ref{prop:tautological-tilting}.
\end{proof}

\subsection{The Auslander--McKay correspondence}\label{sec:relative-one-model-classification}

Retain $\mathsf{Crep}_1(R)$ and $\mathsf{ModGen}(R)$ from Part \ref{part:auslander-mckay}, and $\mathsf{Flag}(a)$ from Section \ref{sec:modifying-classification}.

\begin{Thm}[Auslander--McKay correspondence in relative dimension one]\label{thm:auslander-mckay}
The functor of Theorem \ref{thm:general-auslander-mckay-embedding}(1) gives equivalences
\begin{equation}\label{eq:auslander-mckay}
\begin{gathered}
\mathsf{Crep}_1(R)^{\op}\xrightarrow{\ \sim\ }\mathsf{ModGen}(R)
\xrightarrow{\ \sim\ }\mathsf{Flag}(a),\\
Y^\F\longmapsto T^\F\longmapsto\F,
\end{gathered}
\end{equation}
Every object of $\mathsf{Crep}_1(R)$ is isomorphic over $X$ to $Y^\F$ for a unique flag $\F$. For a nonempty flag of length $m$, its scheme-theoretic closed fibre is a reduced split nodal chain of $m$ projective lines; the empty flag gives $X$. The canonical generator on $Y^\F$ is $\E^\F$, with $C_j$ corresponding to the summand $I_{b_j}$ of $T^\F$. Moreover,
\begin{equation}\label{eq:flag-model-morphisms}
\Hom_X(Y^\F,Y^{\mathcal G})\ne\varnothing
\quad\Longleftrightarrow\quad
\add T^{\mathcal G}\subseteq\add T^\F,
\end{equation}
and this Hom set is a singleton whenever it is nonempty.
\end{Thm}
\begin{proof}
Theorem \ref{thm:general-auslander-mckay-embedding}(1) gives a fully faithful functor from $\mathsf{Crep}_1(R)^{\op}$ to $\mathsf{ModGen}(R)$. Applying Theorem \ref{thm:modifying-classification} with distinguished summand $R$ writes every basic modifying generator as $T^\F$. Theorem \ref{thm:flag-models} puts $Y^\F$ in $\mathsf{Crep}_1(R)$, and Proposition \ref{prop:tautological-tilting} gives
\[
\P_{Y^\F}\cong\E^\F,\qquad \pi_*\P_{Y^\F}=T^\F.
\]
Thus the functor is essentially surjective, and every model is isomorphic over $X$ to a flag model.

By Theorem \ref{thm:modifying-classification} with distinguished summand $R$, each generator determines its flag uniquely, and additive inclusions are precisely inclusions of flags. This gives the second equivalence and uniqueness. The closed-fibre assertion is Theorem \ref{thm:flag-models}, and Lemma \ref{lem:canonical-curve-summands} with \eqref{eq:tautological-degrees} gives the curve--summand correspondence. The morphism criterion of Theorem \ref{thm:general-auslander-mckay-embedding}(1) gives \eqref{eq:flag-model-morphisms}.
\end{proof}

These categories depend, up to equivalence, only on the multiplicity vector $a$; repeated prime factors are recorded by exponent vectors. Deleting entries of a flag realizes the contractions of Theorem \ref{thm:general-auslander-mckay-embedding}, with the corresponding nonfree summands omitted. The direction is reversed: $X$ is terminal in $\mathsf{Crep}_1(R)$, while $R$ is initial in $\mathsf{ModGen}(R)$. Refining a flag to a maximal one shows that every object of $\mathsf{Crep}_1(R)$ is dominated by a locally factorial flag model.

For $S=\mathbb C[[x,y]]$, the restriction to maximal modifying generators is the type $A$ case of the three-dimensional Auslander--McKay correspondence; see \cite[Theorems 1.2 and 4.3]{Wem23survey}. The formulation above also includes partial flags over an arbitrary regular complete local base, with repeated factors recorded by their multiplicities.

\begin{Rem}
Crepancy is essential even for a regular source and a split reduced closed fibre. For $S=k[[x,y]]$ and $f=x$, one has $R\cong k[[y,u,v]]$ and $n=1$, so the only flag model is $X$ itself. The nontrivial blowup $Z=\operatorname{Bl}_{(u,y)}X$ satisfies \eqref{eq:geometric-setup} and has closed fibre $C\cong\mathbb P^1_k$. For its exceptional divisor $E$, the codimension-two blowup formula gives
\[
\omega_Z\cong\O_Z(E),\qquad \omega_Z|_C\cong\O_C(-1).
\]
Thus its canonical sheaf is nontrivial, so $Z$ is not crepant and is not a flag model over $X$.
\end{Rem}

\subsection{Noncommutative crepant resolutions}\label{sec:nccr}

Let $B$ be a Gorenstein local normal domain. A \emph{modifying algebra} over $B$ is an algebra $\End_B(N)$ with $0\ne N\in\refl B$ and $\End_B(N)\in\CM B$. A \emph{noncommutative crepant resolution} (NCCR) of $B$ is a modifying algebra of finite global dimension. We say that $N$ \emph{gives an NCCR} if $\End_B(N)$ is an NCCR.

\begin{Prop}\label{prop:flag-global-dimension}
For a flag $\F=(b_1,\ldots,b_m)$ with blocks $g_1,\ldots,g_{m+1}$ as in \eqref{eq:flag-blocks}, put $\Lambda^\F:=\End_R(T^\F)$. Then
\begin{equation}\label{eq:flag-global-dimension}
\gl\Lambda^\F<\infty
\quad\Longleftrightarrow\quad Y^\F\text{ is regular}
\quad\Longleftrightarrow\quad
g_i\notin\n^2\quad(1\le i\le m+1).
\end{equation}
These conditions hold precisely when $\F$ is maximal and $p_\alpha\notin\n^2$ for every $1\le\alpha\le\ell$.
\end{Prop}
\begin{proof}
Corollary \ref{cor:flag-derived-equivalence} and Lemma \ref{lem:tilting-regularity} give the first equivalence. By Proposition \ref{prop:flag-fibre}, every closed point of the closed fibre other than the distinguished $p_i$ is regular, and
\[
\widehat{\O}_{Y^\F,p_i}\cong S[[x_i,y_i]]/(x_iy_i-g_i).
\]
Since $S[[x_i,y_i]]$ is regular local, this hypersurface is regular exactly when
\[
x_iy_i-g_i\notin(\n,x_i,y_i)^2,
\quad\text{or equivalently}\quad g_i\notin\n^2.
\]
Regularity is unchanged by completion. Every point of $Y^\F$ specializes to a closed point of the closed fibre, by properness over the local base, so regularity there implies regularity everywhere. This proves the second equivalence, including the empty flag.

Each block $g_i$ is, up to a unit, a nonempty product of prime factors of $f$. A block containing at least two factors belongs to $\n^2$. Thus all $g_i$ lie outside $\n^2$ exactly when every block is a single prime factor outside $\n^2$, proving the last assertion.
\end{proof}

\begin{Thm}[NCCR criterion]\label{thm:nccr-criterion}
The ring $R$ admits an NCCR if and only if
\begin{equation}\label{eq:nccr-existence}
p_\alpha\notin\n^2\qquad(1\le\alpha\le\ell).
\end{equation}
When this holds, the basic reflexive modules giving NCCRs are precisely the basic maximal modifying modules, namely
\[
\mathsf T_L(T^\omega)
=\bigoplus_{j=0}^{n-1}(L\otimes_RT_{I_j^\omega})^{**},
\qquad L\in\Cl(R),\quad\omega\in\mathfrak S_n.
\]
Equivalently, for any nonzero modifying module $M$, its endomorphism algebra is an NCCR exactly when $M$ is maximal modifying and \eqref{eq:nccr-existence} holds.
\end{Thm}
\begin{proof}
Passing to the basic representative of $M$ preserves its additive closure and changes its endomorphism algebra only by Morita equivalence. We may therefore assume that $M$ is basic. By \eqref{eq:modifying-flag-classification}, write $M\cong\mathsf T_L(T^\F)$. The composition-compatible isomorphisms \eqref{eq:reflexive-twist-hom} give
\[
\End_R(M)\cong\End_R(T^\F)=\Lambda^\F.
\]
Proposition \ref{prop:flag-global-dimension} shows that this algebra has finite global dimension exactly when $\F$ is maximal and \eqref{eq:nccr-existence} holds. By Corollary \ref{cor:mm-classification}, maximality of $\F$ is equivalent to $M$ being maximal modifying. Conversely, under \eqref{eq:nccr-existence}, any $T^\omega$ gives an NCCR, proving existence as well.
\end{proof}

For $S=k[[x,y]]$, the existence criterion \eqref{eq:nccr-existence} recovers Iyama--Wemyss's \cite[Theorem 5.1(3)]{IW18}: their formulation uses CT modules, whose existence is equivalent to that of an NCCR by \cite[Corollary 5.9]{IW}. Their classification of basic CT modules gives the same NCCR generators $T^\omega$. The argument above extends the criterion to arbitrary regular complete local rings $S$, including mixed characteristic, using the local equations of the flag models and their tilting equivalences.

\section{Reflexive pushforward and wall crossing}\label{sec:reconstruction-pushforward}\label{sec:mutation-geometry}

Mutation compatibility and the algebraic connectivity of MM modules first give realization on a fixed maximal-flag model, in arbitrary dimension and characteristic. We then classify the fibres by the affine Weyl symmetries of repeated prime factors, identify the quotient mutation graph, and compute the algebras and categories of its walls.

Fix a maximal flag $\F_0=(b_1^0,\ldots,b_{n-1}^0)$, with endpoints $b_0^0=0$ and $b_n^0=a$, and write $Y:=Y^{\F_0}$ and $\Psi:=\Psi_Y$. This model remains fixed through Section \ref{sec:pushforward-fibres}; Section \ref{sec:elementary-contractions} returns to arbitrary flags. If $n=1$, then $Y=\Spec R$ is factorial, $\ptilt Y=\{0,\O_Y\}$, $\modif R=\{0,R\}$ and $\MM R=\{R\}$. Statements involving mutations or components are understood for $n\ge2$.

\subsection{Mutation compatibility and realization}\label{sec:fixed-model-realization}\label{sec:mma-derived-equivalence}\label{sec:mutation-compatibility}

Lemma \ref{lem:pretilting-endomorphisms} sends basic pretilting bundles to basic modifying modules, preserving their summands. Proposition \ref{prop:tautological-tilting} identifies the canonical bundle $\P_Y=\E^{\F_0}$, which has $n$ summands. Thus Corollaries \ref{cor:number-of-summands} and \ref{cor:mm-classification} give
\[
\Psi(\tilt Y)\subseteq\MM R,\qquad
\Psi^{-1}(\MM R)=\tilt Y.
\]

\begin{Thm}[Mutation compatibility]\label{thm:pushforward-mutation}\label{prop:pushforward-mutation}
Let $\E$ be a basic tilting bundle on $Y$, with $n\ge2$, and let $V$ be an indecomposable summand. For $M=\Psi(\E)$ and $X=\Psi(V)$,
\begin{equation}\label{eq:pushforward-mutation-compatibility}
\Psi(\mu_V(\E))\cong\mu_X^+(M)\cong\mu_X^-(M).
\end{equation}
\end{Thm}
\begin{proof}
The exchange sequence of Theorem \ref{thm:facet-mutation} has the oriented form of Lemma \ref{lem:oriented-pushforward-mutation}. That lemma identifies $\Psi(\mu_V(\E))$ with one of $\mu_X^+(M)$ and $\mu_X^-(M)$, according to the orientation. The two agree by Corollary \ref{cor:mm-mutation-connectivity}.
\end{proof}

\begin{Cor}[Realization on a fixed maximal-flag model]\label{thm:tilting-modifying-surjection}
Reflexive pushforward gives a surjection $\Psi\colon\ptilt Y\twoheadrightarrow\modif R$, preserves the number of indecomposable summands and satisfies $\Psi^{-1}(\MM R)=\tilt Y$. For $n\ge2$, its restriction to tilting bundles is the chamber surjection
\begin{equation}\label{eq:tilting-modifying-surjections}
\tilt Y\cong\operatorname{Cham}(\widetilde A_{n-1},\varnothing)
\overset{\Psi}{\twoheadrightarrow}\MM R.
\end{equation}
For $n=1$, the restriction is the bijection $\{\O_Y\}\to\{R\}$.
\end{Cor}
\begin{proof}
The case $n=1$ is immediate. For $n\ge2$, Theorem \ref{thm:flag-models} identifies the closed fibre with a split ADE highest-root cycle of type $A_{n-1}$, so Theorem \ref{thm:ade-cycle-fan} gives the chamber identification. Only surjectivity remains. The initial bundle $\P_Y$ pushes forward to $T^{\F_0}$. By Corollary \ref{cor:mm-mutation-connectivity}, any $M\in\MM R$ is reached from it by a finite mutation sequence. At each step, Lemma \ref{lem:pretilting-endomorphisms}(1) matches the indecomposable summands, and Theorem \ref{thm:pushforward-mutation} lifts the mutation to a tilting bundle on $Y$. This proves surjectivity onto $\MM R$.

Every nonzero basic modifying module is a direct summand of a basic MM module, by extending its flag to a maximal flag (Corollary \ref{cor:mm-classification}). Taking the corresponding summands of a tilting lift gives a pretilting lift. Finally, $\Psi(0)=0$.
\end{proof}

For $S=\mathbb C[[x,y]]$ and $n\ge2$, the general Du Val diagram may be larger than $A_{n-1}$; restricting to the $n-1$ retained finite vertices gives the same affine $A_{n-1}$ arrangement. Thus the comparison with \cite[Theorem 0.13]{IWTC} after Corollary \ref{cor:cdv-modifying-realization} applies. Here $S$ may have arbitrary dimension and mixed characteristic.

\begin{Cor}[Derived equivalence of maximal-flag models and MMAs]\label{cor:mma-derived-equivalence}
For any two maximal-flag models $Y,Y'$ and any two maximal modifying algebras $\Lambda,\Lambda'$ of $R$, there are $R$-linear triangle equivalences
\begin{equation}\label{eq:mma-derived-equivalence}
D^b(\Coh Y)\simeq D^b(\Coh Y')\simeq D^b(\mod\Lambda)\simeq D^b(\mod\Lambda').
\end{equation}
These equivalences restrict to perfect complexes. They give a non-commutative Bondal--Orlov analogue and its commutative counterpart for maximal models in $\mathsf{Crep}_1(R)$; compare \cite[Conjectures 1.1 and 1.4]{IW13}.
\end{Cor}
\begin{proof}
Corollary \ref{thm:tilting-modifying-surjection} realizes each MMA, up to Morita equivalence, by a tilting bundle on either fixed model $Y$ or $Y'$; Lemma \ref{lem:tilting-derived-realization} then gives \eqref{eq:mma-derived-equivalence}.
\end{proof}

For $d=2$, compare \cite[Theorem 4.16]{IW}; the argument also applies for $d=1$ and $d\ge3$.

\subsection{Fibres and mutation graphs}\label{sec:pushforward-fibres}\label{sec:mutation-quotients}\label{sec:fixed-model-classification}\label{sec:cone-graph-quotients}

We now describe the ambiguity of the realization in Corollary \ref{thm:tilting-modifying-surjection}. Order the labelled factors so that $(f_i)=(p_\alpha)$ when $b_i^0-b_{i-1}^0=e_\alpha$. The prefixes $\{1,\ldots,j\}$ then form a labelled lift of $\F_0$.

Put $J_0=J_n=\O_Y$, using the $J_j$ of Proposition \ref{prop:tautological-tilting} for $1\le j<n$, and set
\[
\mathcal A_i:=J_i\otimes J_{i-1}^{-1}\quad(1\le i\le n),\qquad
\mathcal A_I:=\bigotimes_{i\in I}\mathcal A_i.
\]
Their tensor product $\bigotimes_{i=1}^n\mathcal A_i$ is trivial. To describe the pushforward, use the presentation
\[
D:=\mathbb Z^\ell/\mathbb Za\cong\Cl(R),\qquad
B_\alpha:=\{i\in[n]\mid(f_i)=(p_\alpha)\},\quad |B_\alpha|=a_\alpha.
\]
Let $\mathbf1=(1,\ldots,1)$ and put $P:=\mathbb Z^n/\mathbb Z\mathbf1$. For $q\in P$, denote by $\mathcal L_q$ the line bundle $\bigotimes_i\mathcal A_i^{q_i}$; this is independent of the representative of $q$.

\Needspace{16\baselineskip}
\begin{Lem}[Picard and class-group coordinates]\label{lem:picard-class-map}
The map $q\mapsto\mathcal L_q$ identifies $P$ with $\Pic(Y)$, and
\begin{equation}\label{eq:picard-factor-coordinates}
\deg(\mathcal L_q|_{C_j})=q_j-q_{j+1}.
\end{equation}
Under this identification, $[\Psi(-)]$ is the surjective homomorphism
\begin{equation}\label{eq:picard-class-map}
\phi\colon P\longrightarrow D,\qquad
q\longmapsto
\overline{\left(\sum_{i\in B_\alpha}q_i\right)_{\alpha=1}^\ell}.
\end{equation}
In particular, $\Psi(\mathcal A_I)\cong T_I$. Its kernel has the concrete presentation
\begin{equation}\label{eq:exceptional-picard-lattice}
K:=\ker\phi
\cong\bigoplus_{\alpha=1}^\ell
\left\{(\lambda_i)_{i\in B_\alpha}\in\mathbb Z^{B_\alpha}
\ \middle|\ \sum_{i\in B_\alpha}\lambda_i=0\right\}.
\end{equation}
The displayed direct sum embeds in $P$ and has rank $n-\ell$.
\end{Lem}
\begin{proof}
Picard restriction and Proposition \ref{prop:tautological-tilting} identify $J_1,\ldots,J_{n-1}$ as a basis of $\Pic(Y)$, with degree vectors $e_j$ and pushforward classes $[T_{\{1,\ldots,j\}}]$. Successive differences give \eqref{eq:picard-factor-coordinates}, with kernel $\mathbb Z\mathbf1$ on $\mathbb Z^n$. The pushforward class map is a homomorphism by restriction to $D(u)\cup D(v)$, so the same differences give \eqref{eq:picard-class-map}.

Every $B_\alpha$ is nonempty, so $\phi$ is surjective. Each class in $K$ has a unique representative with zero block sums: if the block-sum vector of $q$ is $ta$, replace $q$ by $q-t\mathbf1$. This gives \eqref{eq:exceptional-picard-lattice} and the rank assertion.
\end{proof}

Every indecomposable pretilting bundle is a line bundle, by Corollary \ref{cor:flag-type-a} when $n\ge2$ and trivially when $n=1$. For $I\subseteq[n]$, let $\mathbf1_I$ be its characteristic vector, viewed in $P$. The line-bundle criterion \eqref{eq:line-compatibility}, in the coordinates \eqref{eq:picard-factor-coordinates}, becomes
\begin{equation}\label{eq:factor-line-compatibility}
\Ext_Y^{>0}(\mathcal L_q,\mathcal L_{q'})=0
=\Ext_Y^{>0}(\mathcal L_{q'},\mathcal L_q)
\quad\Longleftrightarrow\quad
\max_i(q_i-q'_i)-\min_i(q_i-q'_i)\le1.
\end{equation}

For a chain $\emptyset=I_0\subsetneq I_1\subsetneq\cdots\subsetneq I_r=[n]$, write
\begin{equation}\label{eq:pretilting-flag-expression}
U(q,I_\bullet):=
\bigoplus_{j=0}^{r-1}\mathcal L_{q+\mathbf1_{I_j}}
\qquad(q\in P).
\end{equation}

Let $W_{\mathrm{rep}}:=\prod_{\alpha=1}^\ell\mathfrak S_{B_\alpha}$ act on $P$ by permuting occurrences of each prime factor. It preserves $\phi$, and therefore acts on $K$. Define
\begin{equation}\label{eq:repeated-affine-weyl-group}
\Gamma:=K\rtimes W_{\mathrm{rep}},\qquad
(\lambda,w)q:=\lambda+wq.
\end{equation}
Here $\lambda\in K$ and $w\in W_{\mathrm{rep}}$. By \eqref{eq:exceptional-picard-lattice}, this is the product of affine Weyl groups
\[
\Gamma\cong\prod_{\alpha=1}^\ell
W_{\mathrm{aff}}(A_{a_\alpha-1}),
\]
with the factor for $a_\alpha=1$ understood to be trivial. Sending $\mathcal L_q$ to $\mathcal L_{\lambda+wq}$ and extending over direct sums defines an action on pretilting bundles: translations and coordinate permutations preserve \eqref{eq:factor-line-compatibility}.

\Needspace{22\baselineskip}
\begin{Thm}[Fibres and stabilizers]\label{thm:pushforward-fibres}
The fibres of both pushforward maps in Corollary \ref{thm:tilting-modifying-surjection} are exactly the $\Gamma$-orbits. Hence
\begin{equation}\label{eq:modifying-orbit-classification}
\ptilt Y/\Gamma\xrightarrow{\ \sim\ }\modif R,
\qquad
\tilt Y/\Gamma\xrightarrow{\ \sim\ }\MM R.
\end{equation}
Let $M\ne0$ be basic modifying with $r$ summands. Choose a summand $N_0$ and its anchored chain \eqref{eq:anchored-modifying-chain}, and put
\[
d_{\alpha j}:=b_{j,\alpha}-b_{j-1,\alpha}
\qquad(1\le\alpha\le\ell,\ 1\le j\le r).
\]
For any lift $U$ of $M$, its stabilizer satisfies
\begin{equation}\label{eq:pretilting-fibre-stabilizer}
\Gamma_U\cong\prod_{\alpha=1}^\ell\prod_{j=1}^r
\mathfrak S_{d_{\alpha j}},\qquad
\Psi^{-1}(M)\cong\Gamma/\Gamma_U.
\end{equation}
The fibre over $0$ is $\{0\}$.
\end{Thm}
\begin{proof}
\emph{Flag form.}
Fix $M\ne0$ and $N_0$ as in the statement. By Corollary \ref{thm:tilting-modifying-surjection}, lifts exist and each lift has a unique summand $\mathcal L_q$ mapping to $N_0$. For any other summand $\mathcal L_{q'}$, \eqref{eq:factor-line-compatibility} gives a unique representative $\mathbf1_I$ of $q'-q$ by normalizing its minimum coordinate to zero; here $\emptyset\ne I\subsetneq[n]$. The resulting subsets are comparable, since incomparable subsets give both $1$ and $-1$ in their difference. Hence the lift has a unique form $U(q,I_\bullet)$.

Put $c_j:=(|I_j\cap B_\alpha|)_\alpha$ for $0\le j\le r$. By Lemma \ref{lem:picard-class-map}, the $r$ summands of $\Psi(U)$ have classes $\phi(q)+\overline{c_j}$ for $0\le j<r$. The vectors $c_j$ form a strict chain from $0$ to $a$, so uniqueness in Theorem \ref{thm:modifying-classification} gives $c_j=b_j$.

\emph{Transitivity.}
Since $\phi(\lambda+wq)=\phi(q)$, the action preserves the pushforward of every summand, hence preserves $\Psi$. For two lifts $U(q_0,I_\bullet^0)$ and $U(q,I_\bullet)$, their increment subsets have the same sizes $d_{\alpha j}$ inside every $B_\alpha$. Thus some $w\in W_{\mathrm{rep}}$ satisfies $wI_j^0=I_j$ for all $j$. Now $\lambda:=q-wq_0$ lies in $K$, and $(\lambda,w)$ carries the first lift to the second. This proves the orbit assertions.

\emph{Stabilizers.}
Put $B_{\alpha,j}:=B_\alpha\cap(I_j\setminus I_{j-1})$, of size $d_{\alpha j}$. A stabilizer fixes each summand separately, since their pushforward classes are distinct. Fixing $\mathcal L_q$ forces $\lambda=q-wq$; fixing the remaining summands is equivalent to preserving every $B_{\alpha,j}$. Therefore
\[
\Gamma_U=\{(q-wq,w)\mid w\in\prod_{\alpha,j}\mathfrak S_{B_{\alpha,j}}\}
\cong\prod_{\alpha,j}\mathfrak S_{d_{\alpha j}},
\]
which gives \eqref{eq:pretilting-fibre-stabilizer}. Preservation of summand counts gives the fibre over $0$.
\end{proof}

\begin{Cor}\label{cor:pushforward-bijections}
Every fibre over a maximal modifying module is a $\Gamma$-torsor. Both pushforward maps in Corollary \ref{thm:tilting-modifying-surjection} are bijections exactly when $f$ is squarefree.
\end{Cor}
\begin{proof}
For an MM module, $r=n$ and the $r$ positive column sums $\sum_\alpha d_{\alpha j}$ add to $n$. Each is one, so every factor in \eqref{eq:pretilting-fibre-stabilizer} is trivial. If $f$ is squarefree, then $K=0$ and $W_{\mathrm{rep}}=1$, giving bijectivity. Otherwise $K$ has positive rank, so an MM fibre is a nontrivial torsor and neither map is injective.
\end{proof}

\begin{Rem}[Explicit fibres]
With $N_0,d_{\alpha j}$ as above, choose a lift and let $\mathcal L_{q_0}$ be its summand mapping to $N_0$. The normal form identifies the full fibre with
\begin{equation}\label{eq:explicit-pushforward-fibre}
\Psi^{-1}(M)\cong
K\times\prod_{\alpha=1}^\ell
\operatorname{Part}(B_\alpha;d_{\alpha1},\ldots,d_{\alpha r}),
\end{equation}
recording the distinguished vertex $q_0+\lambda$ and the ordered increment partitions. Here $\operatorname{Part}$ denotes ordered partitions into subsets of the specified sizes, including empty subsets. The finite factor has cardinality
\[
\prod_{\alpha=1}^\ell\frac{a_\alpha!}{\prod_{j=1}^r d_{\alpha j}!}.
\]
If $f$ is not squarefree, $K$ has positive rank, so every nonzero fibre is countably infinite.
\end{Rem}

\begin{Ex}\label{ex:repeated-factor-fibre}
For $f=p^2q$ with $p,q$ nonassociate, $K\cong\mathbb Z$ and $\Gamma\cong W_{\mathrm{aff}}(A_1)$ is the infinite dihedral group. The module $R\oplus(u,p)\oplus(u,p^2)$ has a fibre which is a $\Gamma$-torsor, with underlying set $\mathbb Z\times\mathfrak S_2$ after choices. The module $R\oplus(u,p^2)$ has increments $(2,0),(0,1)$, so its stabilizer is $\mathfrak S_2$ and its fibre is $\Gamma/\mathfrak S_2$, with underlying set $\mathbb Z$.
\end{Ex}

\begin{Rem}[Numerical interpretation]\label{cor:pushforward-cones}
Corollary \ref{cor:rank-determinant} and Lemma \ref{lem:picard-class-map} identify $H=P_\mathbb R$ via $[\mathcal L_q]\mapsto q$, with
\[
\Delta(U(q,I_\bullet))=\Conv\{q+\mathbf1_{I_j}\mid0\le j<r\}.
\]
Thus Theorem \ref{thm:pushforward-fibres} also identifies pushforward classes with $\Gamma$-orbits of cones. For $n\ge2$, the chamber identification in Corollary \ref{thm:tilting-modifying-surjection} gives $\operatorname{Cham}(\widetilde A_{n-1},\varnothing)/\Gamma\cong\MM R$. This interpretation still uses the integral lattices $P,D,K$: the real map $P_\mathbb R\to D_\mathbb R$ alone forgets the torsion in $D$.
\end{Rem}

\begin{Cor}[Mutation graphs]\label{cor:mm-mutation-graph}
Assume $n\ge2$. Form the directed tilting and MM mutation graphs with one arrow for each indecomposable summand being mutated, including selfmutation loops on the module side. Reflexive pushforward identifies the latter graph with the $\Gamma$-quotient of the former. A projected edge is a loop exactly when the exchanged adjacent prime factors are associate. For the convention omitting selfmutations, delete these loops after taking the quotient.
\end{Cor}
\begin{proof}
Theorem \ref{thm:pushforward-fibres} identifies the quotient vertices. The $\Gamma$-action preserves compatibility by \eqref{eq:factor-line-compatibility}, hence adjacency, and preserves each summand's pushforward class. Lemma \ref{lem:pretilting-endomorphisms}(1) matches the summands and Theorem \ref{thm:pushforward-mutation} matches their mutation arrows, giving the quotient graph.

For the loop criterion, anchor at an unchanged summand. Lemma \ref{lem:factor-exchange-mutation} replaces $I(gs)$ by $I(gt)$ for the adjacent prime factors $s,t$, up to a common reflexive twist. The uniqueness of the anchored flag in Theorem \ref{thm:modifying-classification} makes the images equal exactly when $s\sim t$.
\end{proof}

\begin{Rem}\label{rem:fibre-gallery-disconnected}
A fibre need not be connected by mutations with unchanged image. For $f=p^2q$ as in Example \ref{ex:repeated-factor-fibre}, each cyclic prime-factor word has exactly one adjacent pair of associate factors. Thus every tilting bundle has exactly one neighbour with the same pushforward, so the resulting components are single edges, whereas every MM fibre is infinite.
\end{Rem}

\subsection{Contraction algebras and wall categories}\label{sec:elementary-contractions}

We now allow an arbitrary nonempty flag $\F=(b_1,\ldots,b_m)$, with prefix products $G_i$ and blocks $g_i$ as in \eqref{eq:flag-blocks}. Fix $1\le j\le m$, and write
\[
g=G_{j-1},\qquad s=g_j,\qquad t=g_{j+1},\qquad h=f/(gst).
\]
Deleting the $j$th entry gives $c_{\F,j}\colon Y^\F\to Y^{\mathcal G}$. Its local model identifies the exceptional centre; the same ring will describe its contraction algebra and wall category.

\begin{Prop}[Geometry of elementary contractions]\label{prop:elementary-contractions}
Before the flat base change \eqref{eq:flag-flat-base-change}, the contraction $c_{\F,j}$ is the blowup of $(x,s)$ on the merged-block chart $\Spec S[x,y]/(xy-st)$, and is an isomorphism on every other chart. There is a closed subscheme $Z_{\F,j}\subseteq Y^{\mathcal G}$ such that
\[
Z_{\F,j}\cong\Spec S/(s,t),\qquad
c_{\F,j}^{-1}(Z_{\F,j})\cong\mathbb P^1_{S/(s,t)}.
\]
The morphism is an isomorphism off $Z_{\F,j}$. Interchanging $s,t$ gives the blowup of $(x,t)$ over the same target.
\end{Prop}
\begin{proof}
Corollary \ref{cor:flag-product-blowup} and \cite[Tag 080A]{stacks-project} identify $c_{\F,j}$ with the blowup of $I(gs)\O_{Y^{\mathcal G}}$. Before \eqref{eq:flag-flat-base-change}, this ideal is $g(x,s)$ on the merged chart, where $u=gx,v=hy$, and is principal elsewhere by Proposition \ref{prop:flag-charts}. Cancelling $g$ and applying Lemma \ref{lem:elementary-flag-blowup} gives the local model and the isomorphism off $V(x,y,s,t)$.

Over this subscheme, the two charts \eqref{eq:elementary-blowup-charts} reduce to $\Spec\bigl((S/(s,t))[z]\bigr)$ and $\Spec\bigl((S/(s,t))[w]\bigr)$, glued by $w=z^{-1}$. It avoids the other target charts and has $u=v=0$, so it and its projective-line preimage are unchanged by base change to $R$. Interchanging $s,t$ proves the last assertion.
\end{proof}

Put
\[
M=T^\F=\bigoplus_{i=0}^m I(G_i),\qquad \Lambda=\End_R(M).
\]
Let $e_j\in\Lambda$ correspond to $I(G_j)$. By Proposition \ref{prop:tautological-tilting}, the pullback of the target canonical generator is
\[
U_{\F,j}:=\O_{Y^\F}\oplus\bigoplus_{\substack{1\le i\le m\\i\ne j}}J_i.
\]
Its contraction algebra is $A_{\F,j}:=\Lambda/\Lambda(1-e_j)\Lambda$. For $\F=\F_\omega$, also write $A_{\omega,j}:=A_{\F_\omega,j}$.

\begin{Thm}[Elementary contraction algebras]\label{thm:elementary-contraction-algebra}
For adjacent blocks $s=g_j,t=g_{j+1}$, one has
\begin{equation}\label{eq:elementary-contraction-algebra}
A_{\F,j}\cong R/(u,v,s,t)\cong S/(s,t).
\end{equation}
\end{Thm}
\begin{proof}
Put $N:=\bigoplus_{i\ne j}I(G_i)$. Since $\End_R(I(gs))=R$, the quotient $A_{\F,j}$ is $R$ modulo the endomorphisms of $I(gs)$ factoring through $\add N$. The right approximation in \eqref{eq:factor-right-exchange} shows that this ideal is
\[
sH(gs,g)+H(gs,gst)
=\frac{s}{u}I(f/s)+I(t)
=(s,v)+(u,t)=(u,v,s,t),
\]
by \eqref{eq:comparable-ideal-hom} and $f=uv$.
\end{proof}

\begin{Cor}[Smallness and mutation]\label{cor:small-contractions-mutation}
Suppose $\F$ is maximal and choose $\omega$ with $\F=\F_\omega$. Deleting the $j$th entry from $\F_\omega$ and $\F_{\omega s_j}$ gives a common flag $\mathcal G$ and the diagram
\begin{equation}\label{eq:elementary-contraction-diagram}
Y^\omega\xrightarrow{\,c_j\,}Y^{\mathcal G}
\xleftarrow{\,c'_j\,}Y^{\omega s_j},
\end{equation}
in which both exceptional loci have pure codimension $\height_S(s,t)$. Moreover,
\begin{equation}\label{eq:small-mutation-criterion}
c_j\text{ is small}
\ \Longleftrightarrow\ s\not\sim t
\ \Longleftrightarrow\ \mu_{I(G_j)}^\pm(T^\omega)\not\cong T^\omega.
\end{equation}
If $s\sim t$, the contractions are divisorial and $Y^\omega\cong_RY^{\omega s_j}$.
\end{Cor}
\begin{proof}
Proposition \ref{prop:elementary-contractions} gives the common target and the projective-line exceptional loci. For prime $s,t$, the quotient $S/(s,t)$ is pure of dimension $d-\height_S(s,t)$, with $\height_S(s,t)=1$ if $s\sim t$ and $2$ otherwise. This gives the codimension and smallness assertions. Formula \eqref{eq:factor-mutation-formula} and uniqueness in Theorem \ref{thm:modifying-classification} give the selfmutation criterion. In the associate case the flags coincide, so the models are isomorphic.
\end{proof}

\begin{Cor}[Null and wall categories]\label{cor:explicit-wall-category}
One has
\[
\mathcal N_{c_{\F,j}}\simeq\mod S/(g_j,g_{j+1}),\qquad
\W_{U_{\F,j}}=\mathcal N_{c_{\F,j}}\cap\Coh_0Y^\F
\simeq\fl S/(g_j,g_{j+1}).
\]
On a fixed maximal-flag model with $n\ge2$, let $\E=U\oplus V$ be basic tilting and $V$ indecomposable. Anchor the normal form \eqref{eq:pretilting-flag-expression} at a summand other than $V$. Its chain is maximal, so write $I_k=I_k^\omega$ for $\omega\in\mathfrak S_n$, and let $V$ correspond to the $j$th proper prefix. Put $s=f_{\omega(j)},t=f_{\omega(j+1)}$. Then
\begin{equation}\label{eq:explicit-wall-category}
\W_U=U^{\perp_0}\cap U^{\perp_1}\cap\Coh_0Y\simeq\fl S/(s,t),
\end{equation}
also when module mutation is a selfmutation.
\end{Cor}
\begin{proof}
The first assertion is Theorem \ref{thm:contraction-null-category} applied to \eqref{eq:elementary-contraction-algebra}. For the second, the flag form identifies $\Psi(\E)$ with a reflexive twist of $T^\omega$. Lemma \ref{lem:pretilting-endomorphisms}(1) and reflexive twisting then identify the marked endomorphism algebra with $\End_R(T^\omega)$. Apply the same quotient calculation and Theorem \ref{thm:semistable-reduction}.
\end{proof}

When $S=k[[x,y]]$ and $s,t$ are nonassociate primes, the contraction algebra is an Artinian complete intersection of dimension
\[
\dim_kS/(s,t)=i_0(s,t).
\]
This is the threefold criterion of Corollary \ref{cor:contraction-algebra-dimension}, with the intersection multiplicity retained; compare \cite[Theorem 5.5 and Section 6]{Wem23survey}. Over a general $S$, nonassociate prime blocks give a quotient of dimension $d-2$.

\begin{Ex}\label{ex:wall-intersection-multiplicity}
For $S=k[[x,y]]$ and $f=x(x+y^r)$, $r\ge1$, the two maximal-flag models are distinct and their elementary contractions are small, with
\[
A_{\omega,1}\cong k[[y]]/(y^r),\qquad
\W_U\simeq\mod k[y]/(y^r).
\]
For $f=x^2$ there is one maximal-flag model, the contraction is divisorial, and
\[
A_{\omega,1}\cong k[[y]],\qquad \W_U\simeq\fl k[[y]].
\]
The fixed models have the same $\widetilde A_1$ tilting fan; their wall categories distinguish these cases and detect $r$.
\end{Ex}


\bibliographystyle{alpha} 
\bibliography{reference}
\enlargethispage{4\baselineskip}

\end{document}